\documentclass[preprint,3p,12pt]{elsarticle}
\usepackage{mathrsfs}
\usepackage{amsmath}
\usepackage{stmaryrd}
\usepackage{bbding}
\usepackage{dcolumn}
\usepackage{graphicx}
\usepackage{amsfonts}
\usepackage{amssymb}
\usepackage{psfrag}
\usepackage{wrapfig}
\usepackage{subfigure}
\usepackage{makeidx}
\usepackage{bm}
\usepackage{epsf}
\usepackage{color}
\usepackage{epsfig}
\usepackage{setspace}
\usepackage{graphicx}
\usepackage{epstopdf}
\usepackage{psfrag}
\usepackage{subfigure}
\usepackage{multirow}
\usepackage{diagbox}
\usepackage{verbatim}

\usepackage{makecell}
\usepackage{float}
\usepackage{esint}
\usepackage{comment}
\usepackage{movie15}
\usepackage{hyperref}
\usepackage[ruled,linesnumbered]{algorithm2e}

\newcommand{\mathsym}[1]{{}}
\newcommand{\unicode}[1]{{}}
\begin{document}
	
\title{A Memory Reduction Compact Gas Kinetic Scheme on 3D Unstructured Meshes}

	\author[XJTU,HKUST1]{Hongyu Liu}
    \ead{hliudv@connect.ust.hk}
	
	\author[XJTU]{Xing Ji\corref{cor1}}
	\ead{jixing@xjtu.edu.cn}

	\author[XJTU]{Yunpeng Mao}
	\ead{m691810014@stu.xjtu.edu.cn}
		
	\author[XJTU]{Zhe Qian}
    \ead{qzhexj@stu.xjtu.edu.cn}
	
	\author[HKUST1,HKUST2,HKUST3]{Kun Xu}
	\ead{makxu@ust.hk}

	\address[XJTU]{Shaanxi Key Laboratory of Environment and Control for Flight Vehicle, Xi'an Jiaotong University, Xi'an, China}
	\address[HKUST1]{Department of Mathematics, Hong Kong University of Science and Technology, Clear Water Bay, Kowloon, Hong Kong}
	\address[HKUST2]{Department of Mechanical and Aerospace Engineering, Hong Kong University of Science and Technology, Clear Water Bay, Kowloon, Hong Kong}
	\address[HKUST3]{Shenzhen Research Institute, Hong Kong University of Science and Technology, Shenzhen, China}
	\cortext[cor1]{Corresponding author}

\begin{abstract}

This paper presents a memory-reduction third-order compact gas-kinetic scheme (CGKS) designed to solve compressible Euler and Navier-Stokes equations on 3D unstructured meshes. 
Utilizing the time-accurate gas distribution function, the gas kinetic solver provides a time-evolution solution at the cell interface, distinguishable from the Riemann solver with a constant solution. 
With the time evolution solution at the cell interface, evolving both the cell-averaged flow variables and the cell-averaged slopes of flow variables becomes feasible.
Therefore, with the cell-averaged flow variables and their slopes inside each cell, the Hermite WENO (HWENO) techniques can be naturally implemented for the compact high-order reconstruction at the beginning of the next time step. 
However, the HWENO reconstruction method requires the storage of a reconstruction-coefficients matrix for the quadratic polynomial to achieve third-order accuracy, leading to substantial memory consumption. This memory overhead limits both computational efficiency and the scale of simulations.

A novel reconstruction method, built upon HWENO reconstruction, has been designed to enhance computational efficiency and reduce memory usage compared to the original CGKS.
The simple idea is that the first-order and second-order terms of the quadratic polynomials are determined in a two-step way. In the first step, the second-order terms are obtained from the reconstruction of a linear polynomial of the first-order derivatives by only using the cell-averaged slopes, since the second-order derivatives are nothing but the "derivatives of derivatives". Subsequently, the first-order terms left can be determined by the linear reconstruction only using cell-averaged values. Thus, we successfully split one quadratic least-square regression into several linear least-square regressions, which are commonly used in a second-order finite volume code. Since only a $3\times3$ matrix inversion is needed in a 3-D linear least-square regression, the computational cost for the new reconstruction is dramatically reduced and the storage of the reconstruction-coefficient matrix is no longer necessary.
The proposed memory reduction CGKS is tested for both inviscid and viscous flow at low and high speeds on hybrid unstructured meshes.
The proposed new reconstruction technique can reduce the overall computational cost by about $20\%$ to $30\%$. In the meantime, it also simplifies the algorithm.
The challenging large-scale unsteady numerical simulation is performed, which demonstrates that the current improvement brings the CGKS to a new level for industrial applications.

\end{abstract}


\maketitle

\section{Introduction}

The second-order finite volume method (FVM) enjoys widespread adoption in commercial CFD software owing to its computational efficiency and robustness \cite{kim2000second}. 
The cell-based Green-Gauss method is commonly employed for slope reconstruction. However, it may suffer from decreased spatial accuracy on skewed meshes and may exhibit over-dissipation when simulating flow with discontinuities \cite{shima2013green}. 
Another slope reconstruction method utilized in second-order FVM is the least-square reconstruction technique, which employs cell-averaged variables and considers von Neumann neighbors.
While this method can achieve strict second-order accuracy, it tends to exhibit linear instability on tetrahedral grids \cite{haider2009stability}. 
To ensure linear stability, extended stencils are needed. 
In recent decades, high-order numerical methods in computational fluid dynamics (CFD) have garnered significant success \cite{HUYNH2014209,abgrall2018high}. These methods, such as Weighted Essentially Non-Oscillatory (WENO) \cite{gao2012high}, Discontinuous Galerkin (DG) \cite{shu2003fd-fv-weno-dg-review,shu2016weno-dg-review}, Correction Procedure Via Reconstruction (CPR) \cite{CPR2011,wang2017towards}, and Variational Finite Volume (VFV) \cite{cances2020variational}, have played a crucial role in improving the accuracy of numerical simulations \cite{bhatti2020recent}. They are well-suited for handling complex flow phenomena, including turbulent flows \cite{li2010direct}, shock propagation \cite{romick2017high}, and multi-physics coupling \cite{ferrer2023high}.
The traditional high-order FVMs have been developed for decades and used in some aeronautical simulations \cite{ollivier2009obtaining}. 
However, in the pursuit of high-order accuracy, high-order FVM methods often encounter a significant rise in memory consumption \cite{hu1999weighted}. 
This increase in memory not only disrupts the continuity of memory storage but also requires memory allocation on the heap that should ideally be stored temporarily on the stack.
These memory-intensive methods not only bring more challenges in programming but also increase difficulties in application to massively parallel computing architectures.

Compact high-order schemes do not need to store too much neighbor cell's topology and geometry information, as a result, it is not only more friendly to memory access but also reduce some memory usage \cite{ferrer2023high,wen2020improved}. 
Hence, in recent years, research on compact high-order methods has emerged as a prominent area of interest.
The compact methods with the updating of multiple degrees of freedom (DOFs) for each cell have been developed extensively in the past decades, such as DG and CPR.
These methods demonstrate the capability to achieve arbitrary spatial order on arbitrary cells, highlighting their exceptional mesh adaptability. Additionally, both DG and CPR only require updates for the targeted cell's degrees of freedom (DOFs).
Thus, these methods are naturally suitable for parallel computation, which indicates high scalability.
Based on the above advantages, the mesh moving and deformation techniques can be applied directly to DG and CPR without loss of accuracy.
Successful examples have been demonstrated in the DG Arbitrary Lagrangian Eulerian (ALE) method \cite{ren2016multi} and overset mesh using the CPR method \cite{duan2020high}.
Large-scale simulations like large eddy simulation (LES) \cite{wang2017towards} can be done easily due to the high scalability of these methods.
Though these methods behave very well in smooth regions, when facing discontinuities, they have less robustness compared with traditional FVM methods. 
In addition, their explicit time step is bounded by the order of spatial accuracy. 
The $P_{N}P_{M}$ \cite{dumbser2010arbitrary} and reconstructed-DG (rDG) \cite{xuan2014rdg-gks} methods are targeting to solve the above problem. 
Large time steps and less memory requirement can be achieved in the rDG methods in comparison with the same order DG ones.
However, rDG will release the compactness of the original DG methods to do reconstruction. 

In recent years, the high-order compact gas-kinetic scheme(CGKS) \cite{ji2018compact,zhao2023high,zhao2019compact} has been developed based on gas-kinetic theory.
CGKS uses the time-dependent distribution function which has an accurate analytical integral solution at the cell interface of the Bhatnagar--Gross--Krook equation \cite{BGK}.
By using the accurate time-dependent distribution function at the cell interface, not only the Navier-Stokes flux functions can be obtained, but also the time-accurate macroscopic flow variables will be evaluated. 
This implies as cell-averaged flow variables are updated within the finite volume framework, we can also update the cell-averaged slopes as another DOF in each cell to get as much information as possible with the smallest stencil possible.
With the cell-averaged flow variables and their slopes, a Hermite Weighted Essentially Non-Oscillatory (HWENO) \cite{li2021multi} method can be employed for the reconstruction, which can be found in our previous work of third-order CGKS \cite{zhang2023high}. 
As for temporal discretization, explicit two-stage fourth-order and other multi-stage multi-derivative time marching schemes can be used in CGKS for high-order temporal discretization \cite{li2016twostage}.
Due to the benefits of the more reliable evolution process based on mesoscopic gas kinetic theory, the CGKS has excellent performance on both smooth and discontinuous flow regimes. 
The CGKS also shows good performance in the regime of unsteady compressible flow, such as in computational aeroacoustics \cite{zhao2019acoustic} and implicit large eddy simulation \cite{ji2021compact}.
In CGKS, the discontinuity feedback factor (DF) extends beyond the assumptions of the finite volume framework by determining the presence of discontinuities within the cell for the upcoming time step. 
Thus, it further improves the robustness of the CGKS when facing strong discontinuity.
 As a result, CGKS also performs well in supersonic and hypersonic flow simulations on the three-dimensional hybrid unstructured mesh, such as YF-17 fighter jet and X-38 type spaceship \cite{ji2021gradient}. 
 In summary, the CGKS demonstrates strong grid adaptability and robustness, effectively handling the complexities of low-quality meshes and pronounced discontinuities.

Nevertheless, high-order numerical methods cause substantial memory consumption, which can constrain computational efficiency and the scalability of numerical simulations.
DG and CPR methods need to arrange Gaussian points in each cell and the face of the grid. 
For model DG, the hexahedron-type cell needs $N^{3}$ Gaussian points and the quadrilateral-type face needs $N^{2}$ Gaussian points, where $N$ is the order of the numerical method \cite{luo2010reconstructed}. 
For the nodal DG or CPR method, each hexahedron-type cell needs $3N^{2}(2N+1)$ Gaussian points and each quadrilateral-type face needs $2N(N+1)$ Gaussian points.
The original HWENO reconstruction used by CGKS requires a large memory consumption for the reconstruction matrix \cite{ji2021compact}. 
Each hexahedron-type cell needs to store a matrix of 24 rows and 9 columns and $6(N-1)^2$ Gaussian points.

This paper proposes a memory-reduction CGKS method aimed at further diminishing the memory overhead of the original CGKS \cite{ji2021gradient,ji2021compact,zhang2023high} and enhancing computational efficiency.
The primary strategy for memory reduction involves replacing the original linear reconstruction method of CGKS with a two-step third-order reconstruction approach.
Firstly, using cell-averaged slopes to do the least square reconstruction once.
Then, the coefficients of quadratic terms of the polynomial can be obtained.
Secondly, moving the quadratic terms to the right-hand side (RHS) of the original HWENO linear equations, means only linear terms need to be solved. 
Then, using the least square reconstruction again to obtain the first-order terms’ coefficients.
As described above, the new reconstruction method is matrix-free and only needs to compute the small coefficient matrix of the least square in the subroutine. A similar idea has been adopted in the hybrid DG/FVM method where the high-order term is reconstructed in the same manner. However, the treatments for the low-order term are different: the DG evolution is adopted in hybrid DG/FVM, while the reconstruction is still applied in the current CGKS \cite{maltsev2023hybrid}. 
The proposed memory reduction CGKS is tested for both inviscid and viscous flow at low and high speeds on hybrid unstructured meshes, demonstrating the current method's accuracy, robustness, and efficiency improvement.

The paper is organized as follows. 
In Section 2, the 3D BGK equation, the finite volume framework, and the construction of CGKS on three-dimensional hybrid unstructured meshes will be introduced.
In Section 3, the memory reduction two-step third-order linear spatial reconstruction and the nonlinear limiting procedure will be introduced.  
In Section 4, numerical examples including both inviscid and viscous flow computations will be given. The last section is the conclusion.

\section{Gas kinetic scheme under the finite volume framework}
\subsection{3-D BGK equation}
The Boltzmann equation \cite{cercignani1988boltzmann} describes the evolution of molecules at the mesoscopic scale. It indicates that each particle will either transport with a constant velocity or encounter a two-body collision. The BGK \cite{BGK} model simplifies the Boltzmann equation by replacing the full collision term with a relaxation model.
The 3-D gas-kinetic BGK equation \cite{BGK} is
\begin{equation}\label{bgk}
	f_t+\textbf{u}\cdot\nabla f=\frac{g-f}{\tau},
\end{equation}
where $f=f(\textbf{x},t,\textbf{u},\xi)$ is the gas distribution function, which is a function of space $\textbf{x}$, time $t$, phase space velocity $\textbf{u}$, and internal variable $\xi$.
$g$ is the equilibrium state
and $\tau$ is the collision time, which means an averaged time interval between two sequential collisions.  $g$ is expressed as a Maxwellian distribution function.
Meanwhile, the collision term on the right-hand side (RHS) of Eq.~\eqref{bgk} should satisfy the compatibility condition
\begin{equation*}\label{compatibility}
	\int \frac{g-f}{\tau} \pmb{\psi} \text{d}\Xi=0,
\end{equation*}
where $\pmb{\psi}=(1,\textbf{u},\displaystyle \frac{1}{2}(\textbf{u}^2+\xi^2))^T$,
$\text{d}\Xi=\text{d}u_1\text{d}u_2\text{d}u_3\text{d}\xi_1...\text{d}\xi_{K}$,
$K$ is the number of internal degrees of freedom, i.e.
$K=(5-3\gamma)/(\gamma-1)$ in the 3-D case, and $\gamma$
is the specific heat ratio. The details of the BGK equation can be found in \cite{xu2014directchapter2}.

In the continuous flow regime, distribution function $f$ can be taken as a small-scale expansion of Maxwellian distribution. Based on the Chapman-Enskog expansion \cite{CE-expansion}, the gas distribution function can be expressed as \cite{xu2014directchapter2},
\begin{align*}
	f=g-\tau D_{\textbf{u}}g+\tau D_{\textbf{u}}(\tau
	D_{\textbf{u}})g-\tau D_{\textbf{u}}[\tau D_{\textbf{u}}(\tau
	D_{\textbf{u}})g]+...,
\end{align*}
where $D_{\textbf{u}}={\partial}/{\partial t}+\textbf{u}\cdot \nabla$.
Through zeroth-order truncation $f=g$, the Euler equation can be obtained. The Navier-Stokes (NS) equations,
\begin{equation*}\label{ns-conservation}
	\begin{split}
		\textbf{W}_t+ \nabla \cdot \textbf{F}(\textbf{W},\nabla \textbf{W} )=0,
	\end{split}
\end{equation*}
 can be obtained by taking first-order truncation, i.e.,
\begin{align} \label{ce-ns}
	f=g-\tau (\textbf{u} \cdot \nabla g + g_t),
\end{align}
with $\tau = \mu / p$ and $Pr=1$.

Benefiting from the time-accurate gas distribution function, a time evolution solution at the cell interface is provided by the gas kinetic solver, which is distinguishable from the Riemann solver with a constant solution \cite{yang2022comparison}. This is a crucial point to construct the compact high-order gas kinetic scheme.
\begin{align}\label{point}
	\textbf{W}(\textbf{x},t)=\int \pmb{\psi} f(\textbf{x},t,\textbf{u},\xi)\text{d}\Xi,
\end{align}
and the flux at the cell interface can also be obtained
\begin{equation}\label{f-to-flux}
	\textbf{F}(\textbf{x},t)=
	\int \textbf{u} \pmb{\psi} f(\textbf{x},t,\textbf{u},\xi)\text{d}\Xi.
\end{equation}

\subsection{Finite volume framework}
The boundary of a three-dimensional arbitrary polyhedral cell $\Omega_i$ can be expressed as
\begin{equation*}
	\partial \Omega_i=\bigcup_{p=1}^{N_f}\Gamma_{ip},
\end{equation*}
where $N_f$ is the number of cell interfaces for cell $\Omega_i$.
$N_f=4$ for tetrahedron, $N_f=5$ for prism and pyramid, $N_f=6$ for hexahedron.
The semi-discretized form of the finite volume method for conservation laws can be written as
\begin{equation}\label{semidiscrete}
	\frac{\text{d} \textbf{W}_{i}}{\text{d}t}=\mathcal{L}(\textbf{W}_i)=-\frac{1}{\left| \Omega_i \right|} \sum_{p=1}^{N_f} \int_{\Gamma_{ip}}
	\textbf{F}(\textbf{W}(\textbf{x},t))\cdot\textbf{n}_p \text{d}s,
\end{equation}
with
\begin{equation*}\label{f-to-flux-in-normal-direction}
	\textbf{F}(\textbf{W}(\textbf{x},t))\cdot \textbf{n}_p=\int\pmb{\psi}  f(\textbf{x},t,\textbf{u},\xi) \textbf{u}\cdot \textbf{n}_p \text{d}\Xi,
\end{equation*}
where $\textbf{W}_{i}$ is the cell averaged values over cell $\Omega_i$, $\left|
\Omega_i \right|$ is the volume of $\Omega_i$, $\textbf{F}$ is the interface fluxes, and $\textbf{n}_p=(n_1,n_2,n_3)^T$ is the unit vector representing the outer normal direction of $\Gamma_{ip}$.
Through the iso-parametric transformation,
the Gaussian quadrature points can be determined and $\textbf{F}_{ip}(t)$ can be approximated by the numerical quadrature
\begin{equation*}\label{fv-3d-general-quadrature}
	\int_{\Gamma_{ip}}
	\textbf{F}(\textbf{W}(\textbf{x},t))\cdot\textbf{n}_p \text{d}s =  S_{i,p} \sum_{k=1}^{M} \omega_k
	\textbf{F}(\textbf{x}_{p,k},t)\cdot\textbf{n}_p,
\end{equation*}
where $S_{i,p}$ is the area of $\Gamma_{ip}$. Through the iso-parametric transformation,
in the current study, the linear element is considered.
To meet the requirement of a third-order spatial accuracy,
three Gaussian points are used for a triangular face and four Gaussian points are used for a quadrilateral face.
In the computation, the fluxes are obtained under the local coordinates.
The details can be found in \cite{ji2021gradient,pan2020high,JI2024112590}.

\subsection{Gas kinetic solver}
To obtain the numerical flux at the cell interface, the integration solution based on the BGK equation is used
\begin{equation}\label{integral1}
	f(\textbf{x},t,\textbf{u},\xi)=\frac{1}{\tau}\int_0^t g(\textbf{x}',t',\textbf{u},\xi)e^{-(t-t')/\tau}\text{d}t'
	+e^{-t/\tau}f_0(\textbf{x}-\textbf{u}t,\textbf{u},\xi),
\end{equation}
where $\textbf{x}=\textbf{x}'+\textbf{u}(t-t')$ is the particle trajectory. $f_0$ is the initial gas distribution function, $g$ is the corresponding
equilibrium state in the local space and time.
This integration solution describes the physical picture of the particle evolution. Starting with an initial state $f_0$, the particle will take free transport with a probability of $e^{-\Delta t/\tau}$. Otherwise, it will suffer a series of collisions. The effect of collisions is driving the system to the local Maxwellian distribution $g$, and the particles from the equilibrium propagate along the characteristics with a surviving probability of $e^{-(t-t')/\tau}$.
The components of the numerical fluxes at the cell interface can be categorized as the upwinding free transport from $f_0$ and the central difference from the integration of the equilibrium distribution.

To construct a time-evolving gas distribution function at a cell interface,
the following notations are introduced first
\begin{align*}
	a_{x_i} \equiv  (\partial g/\partial x_i)/g=g_{x_i}/g,
	A \equiv (\partial g/\partial t)/g=g_t/g,
\end{align*}
where $g$ is the equilibrium state.  The partial derivatives $a_{x_i}$ and $A$, denoted by $s$,
have the form of
\begin{align*}
	s=s_j\psi_j =s_{1}+s_{2}u_1+s_{3}u_2+s_{4}u_3
	+s_{5}\displaystyle \frac{1}{2}(u_1^2+u_2^2+u_3^2+\xi^2).
\end{align*}
The initial gas distribution function in Eq.~\eqref{integral1} can be modeled as
\begin{equation*}
	f_0=f_0^l(\textbf{x},\textbf{u})\mathbb{H} (x_1)+f_0^r(\textbf{x},\textbf{u})(1- \mathbb{H}(x_1)),
\end{equation*}
where $\mathbb{H}(x_1)$ is the Heaviside function. Here $f_0^l$ and $f_0^r$ are the
initial gas distribution functions on the left and right sides of a cell
interface, which can be fully determined by the initially
reconstructed macroscopic variables. The first-order
Taylor expansion for the gas distribution function in space around
$\textbf{x}=\textbf{0}$ can be expressed as
\begin{align}\label{flux-3d-1}
	f_0^k(\textbf{x})=f_G^k(\textbf{0})+\frac{\partial f_G^k}{\partial x_i}(\textbf{0})x_i
	=f_G^k(\textbf{0})+\frac{\partial f_G^k}{\partial x_1}(\textbf{0})x_1
	+\frac{\partial f_G^k}{\partial x_2}(\textbf{0})x_2
	+\frac{\partial f_G^k}{\partial x_3}(\textbf{0})x_3,
\end{align}
for $k=l,r$.
According to Eq.~\eqref{ce-ns}, $f_{G}^k$ has the form
\begin{align}\label{flux-3d-2}
	f_{G}^k(\textbf{0})=g^k(\textbf{0})-\tau(u_ig_{x_i}^{k}(\textbf{0})+g_t^k(\textbf{0})),
\end{align}
where $g^k$ is the equilibrium state with the form of a Maxwell distribution.
$g^k$ can be fully determined from the
reconstructed macroscopic variables $\textbf{W}
^l, \textbf{W}
^r$ at the left and right sides of a cell interface
\begin{align}\label{get-glr}
	\int\pmb{\psi} g^{l}\text{d}\Xi=\textbf{W}
	^l,\int\pmb{\psi} g^{r}\text{d}\Xi=\textbf{W}
	^r.
\end{align}
Substituting Eq.~\eqref{flux-3d-1} and Eq.~\eqref{flux-3d-2} into Eq.~\eqref{integral1},
the kinetic part for the integral solution can be written as
\begin{equation}\label{dis1}
	\begin{aligned}
		e^{-t/\tau}f_0^k(-\textbf{u}t,\textbf{u},\xi)
		=e^{-t/\tau}g^k[1-\tau(a_{x_i}^{k}u_i+A^k)-ta^{k}_{x_i}u_i],
	\end{aligned}
\end{equation}
where the coefficients $a_{x_1}^{k},...,A^k, k=l,r$ are defined according
to the expansion of $g^{k}$.
After determining the kinetic part
$f_0$, the equilibrium state $g$ in the integral solution
Eq.~\eqref{integral1} can be expanded in space and time as follows
\begin{align}\label{equli}
	g(\textbf{x},t)= g^{c}(\textbf{0},0)+\frac{\partial  g^{c}}{\partial x_i}(\textbf{0},0)x_i+\frac{\partial  g^{c}}{\partial t}(\textbf{0},0)t,
\end{align}
where $ g^{c}$ is the Maxwellian equilibrium state located at an interface.
Similarly, $\textbf{W}^c$ are the macroscopic flow variables for the determination of the
equilibrium state $ g^{c}$
\begin{align}\label{compatibility2}
	\int\pmb{\psi} g^{c}\text{d}\Xi=
	\int_{u>0}\pmb{\psi} g^{l}\text{d}\Xi+
	\int_{u<0}\pmb{\psi} g^{r}\text{d}\Xi=\textbf{W}^c.
\end{align}
Substituting Eq.~\eqref{equli} into Eq.~\eqref{integral1}, the collision part in the integral solution
can be written as
\begin{equation}\label{dis2}
	\begin{aligned}
		\frac{1}{\tau}\int_0^t
		g&(\textbf{x}',t',\textbf{u},\xi)e^{-(t-t')/\tau}\text{d}t'
		=C_1 g^{c}+C_2 a_{x_i}^{c} u_i g^{c} +C_3 A^{c} g^{c} ,
	\end{aligned}
\end{equation}
where the coefficients
$a_{x_i}^{c},A^{c}$ are
defined from the expansion of the equilibrium state $ g^{c}$. The
coefficients $C_m, m=1,2,3$ in Eq.~\eqref{dis2}
are given by
\begin{align*}
	C_1=1-&e^{-t/\tau}, C_2=(t+\tau)e^{-t/\tau}-\tau, C_3=t-\tau+\tau e^{-t/\tau}.
\end{align*}
The coefficients in Eq.~\eqref{dis1} and Eq.~\eqref{dis2}
can be determined by the spatial derivatives of macroscopic flow
variables and the compatibility condition as follows
\begin{align}\label{co}
	&\langle a_{x_1}\rangle =\frac{\partial \textbf{W} }{\partial x_1}=\textbf{W}_{x_1},
	\langle a_{x_2}\rangle =\frac{\partial \textbf{W} }{\partial x_2}=\textbf{W}_{x_2},
	\langle a_{x_3}\rangle =\frac{\partial \textbf{W} }{\partial x_3}=\textbf{W}_{x_3},\nonumber\\
	&\langle A+a_{x_1}u_1+a_{x_2}u_2+a_{x_3}u_3\rangle=0,
\end{align}
where $\left\langle ... \right\rangle$ are the moments of a gas distribution function defined by
\begin{align}\label{co-moment}
	\langle (...) \rangle  = \int \pmb{\psi} (...) g \text{d} \Xi .
\end{align}
The specific details of constructing the second-order flux on the interfaces and the formula of numerical dissipation can be found in \cite{ji2021compact}

\subsection{Direct evolution of the cell averaged slopes} \label{slope-section}

The time evolution solution at a cell interface is provided by the gas-kinetic solver, which is distinguished from the Riemann solvers with a constant solution.
By recalling Eq.~\eqref{point}, the conservative variables at the Gaussian point  $\textbf{x}_{p,k}$ can be updated through the moments $\pmb{\psi}$
of the gas distribution function,
\begin{equation*}\label{point-interface}
	\begin{aligned}
		\textbf{W}_{p,k}(t^{n+1})=\int \pmb{\psi} f^n(\textbf{x}_{p,k},t^{n+1},\textbf{u},\xi) \text{d}\Xi,~ k=1,...,M.
	\end{aligned}
\end{equation*}

Then, the cell-averaged slopes within each element at $t^{n+1}$ can be evaluated based on the divergence theorem,

\begin{equation*}\label{gauss-formula}
	\begin{aligned}
		\nabla \overline{W}^{n+1} \left| \Omega \right|
		&=\int_{\Omega} \nabla \overline{W}(t^{n+1}) \text{d}V
		=\int_{\partial \Omega} \overline{W}(t^{n+1}) \textbf{n} \text{d}S
		= \sum_{p=1}^{N_f}\sum_{k=1}^{M_p} \omega_{p,k} W^{n+1}_{p,k} \textbf{n}_{p,k} \Delta S_p,
	\end{aligned}
\end{equation*}
where $\textbf{n}_{p,k}=((n_{1})_{p,k},(n_{2})_{p,k},(n_{3})_{p,k})$ is the outer unit normal direction at each Gaussian point $\textbf{x}_{p,k}$.

\section{Spatial reconstruction}

In this section, the new memory reduction 3rd-order compact reconstruction is presented with cell-averaged values and cell-averaged slopes. To maintain both shock-capturing ability and robustness, WENO weights and DF are used \cite{ji2021gradient}. Further improvement has been made to make reconstruction simple and more robust \cite{zhang2023high}. Only one large stencil and one sub-stencil are involved in the WENO procedure.

Firstly, we will introduce the HWENO method used in the original CGKS  \cite{ji2020hweno}  for the large stencil.
Secondly, the memory reduction two-step third-order reconstruction for the large stencil will be introduced.
Thirdly, the nonlinear WENO weights and DF will be introduced.

\subsection{Original 3rd-order compact reconstruction for large stencil}
Firstly, a linear reconstruction approach is presented. To achieve a third-order accuracy in space, a quadratic polynomial $p^2$ is constructed as follows

\begin{equation}
p^2(\mathbf{x})=\overline{Q}_0+\sum\limits_{|k|=1}^2a_k[(x-x_0)^{k_1} (y-y_0)^{k_2} (z-z_0)^{k_3}-\frac{1}{\left|\Omega_0\right|} \iiint_{\Omega_0} (x-x_0)^{k_1} (y-y_0)^{k_2} (z-z_0)^{k_3}] \mathrm{~d} V,
\end{equation}
where $k=\left(k_1, k_2, k_3\right)$ is the multi-index, $k|=k_1+k_2+k_3$, $(x_0, y_0, z_0)$ is the geometric center coordinate. The $p^2$ on $\Omega_0$  is constructed on the compact stencil $S$ including $\Omega_0$ and all its von Neumann neighbors $\Omega_m$ ($m=1,\cdots,N_f$, $N_f$ is the number of $\Omega_0$'s faces). The cell averages $\overline{Q}$ on $\Omega_0$ and $\Omega_m$ together with cell-averaged slopes $\overline{Q}_x,\overline{Q}_y $ and $\overline{Q}_z$ on $\Omega_m$ are used to obtain $p^2$.

The polynomial $p^2$ naturally satisfies cell averages over $\Omega_0$ 
\begin{equation}\label{self-constrain}
	\iiint_{\Omega_0} p^2 \text{d}V = \overline{Q}_0|\Omega_0|,
\end{equation}
Meanwhile, the $p^2$ is also required to exactly satisfy cell averages over the target cell's neighbors $\Omega_m$ 
\begin{equation}\label{neighbor-constrain}
 \iiint_{\Omega_m} p^2 \text{d}V = \overline{Q}_m|\Omega_m|\\
\end{equation}
Then, the following conditions are satisfied in a least-square sense
\begin{equation*}
	\begin{aligned}
		\iiint_{\Omega_m}\frac{\partial}{\partial x} p^2 \text{d}V = \left(\overline{Q}_x\right)_m|\Omega_m|\\
		\iiint_{\Omega_m}\frac{\partial}{\partial y} p^2 \text{d}V = \left(\overline{Q}_y\right)_m|\Omega_m|\\
		\iiint_{\Omega_m}\frac{\partial}{\partial z} p^2 \text{d}V = \left(\overline{Q}_z\right)_m|\Omega_m|.\\
	\end{aligned}
\end{equation*}
The constrained least-square method is used to meet the above requirements.

\subsection{The memory reduction 3rd-order compact reconstruction for large stencil}
The new memory reduction third-order reconstruction for $p^2$ of the large stencil consists of two steps: \\
\textbf{Step 1}\\
Having the cell-averaged slopes, we can reconstruct the distribution of these slopes in space, which means the coefficients of the quadratic terms of $p^2$ can be obtained.
$L_{x}^{1}$ is the linear polynomial of the x-direction cell-averaged slope;
$L_{y}^{1}$ is the linear polynomial of the y-direction cell-averaged slope;
$L_{z}^{1}$ is the linear polynomial of the z-direction cell-averaged slope;
\begin{equation}
    \begin{aligned}
    L_{x}^{1}=b_0+b_1(x-x_0)+b_2(y-y_0)+b_3(z-z_0),\\
    L_{y}^{1}=c_0+c_1(x-x_0)+c_2(y-y_0)+c_3(z-z_0),\\
    L_{z}^{1}=d_0+d_1(x-x_0)+d_2(y-y_0)+d_3(z-z_0).\\
    \end{aligned}
\end{equation}
Substituting the first-order terms' coefficients into the $p^2$ of the large stencil, the quadratic terms can be written as:

\begin{equation}
    a_4=b_1, a_5=c_2, a_6=d_3
\end{equation}

\begin{equation}
    a_7=\frac{b_1+c_1}{2}, a_8=\frac{c_3+d_2}{2}
    , a_9=\frac{b_3+d_1}{2}\\
\end{equation}
\textbf{Step2}\\
After step 1, only the linear term in $p^2$ will be determined.
Moving the quadratic terms to the RHS of the $p^2$, a new linear system of the linear terms of $p^2$ can be obtained.
For convenience, we denote the integration of the quadratic terms as
\begin{equation}
\begin{aligned}
    	R_{m}^2 &= \iiint_{\Omega_m} [a_4(x-x_0)^2+a_5(y-y_0)^2+a_6(z-z_0)^2 \\
		&+a_7(x-x_0)(y-y_0)+a_8(y-y_0)(z-z_0) \\
        &+a_9(x-x_0)(z-z_0)] \text{d}V-R_{0}^2\Omega_m,
\end{aligned}
\end{equation}
in which $R_{0}^2$ is
\begin{equation}
\begin{aligned}
    	R_{0}^2 &= \frac{1}{\left|\Omega_0\right|}\iiint_{\Omega_0} [a_4(x-x_0)^2+a_5(y-y_0)^2+a_6(z-z_0)^2 \\
		&+a_7(x-x_0)(y-y_0)+a_8(y-y_0)(z-z_0) \\
        &+a_9(x-x_0)(z-z_0)] \text{d}V
\end{aligned}
\end{equation}

Based on Eq.~\eqref{self-constrain} and Eq.~\eqref{neighbor-constrain}, the new linear system for the linear terms of $p^2$ can be written as: 
\begin{equation}
\begin{aligned}
          \iiint_{\Omega_m} a_1(x-x_0)+a_2(y-y_0)+a_3(z-z_0) \text{d}V = \overline{Q}_m|\Omega_m|
     -R_m^2\\
\end{aligned}
\end{equation}

\begin{figure}[htp]
	\centering
	\includegraphics[width=1.0\textwidth]
	{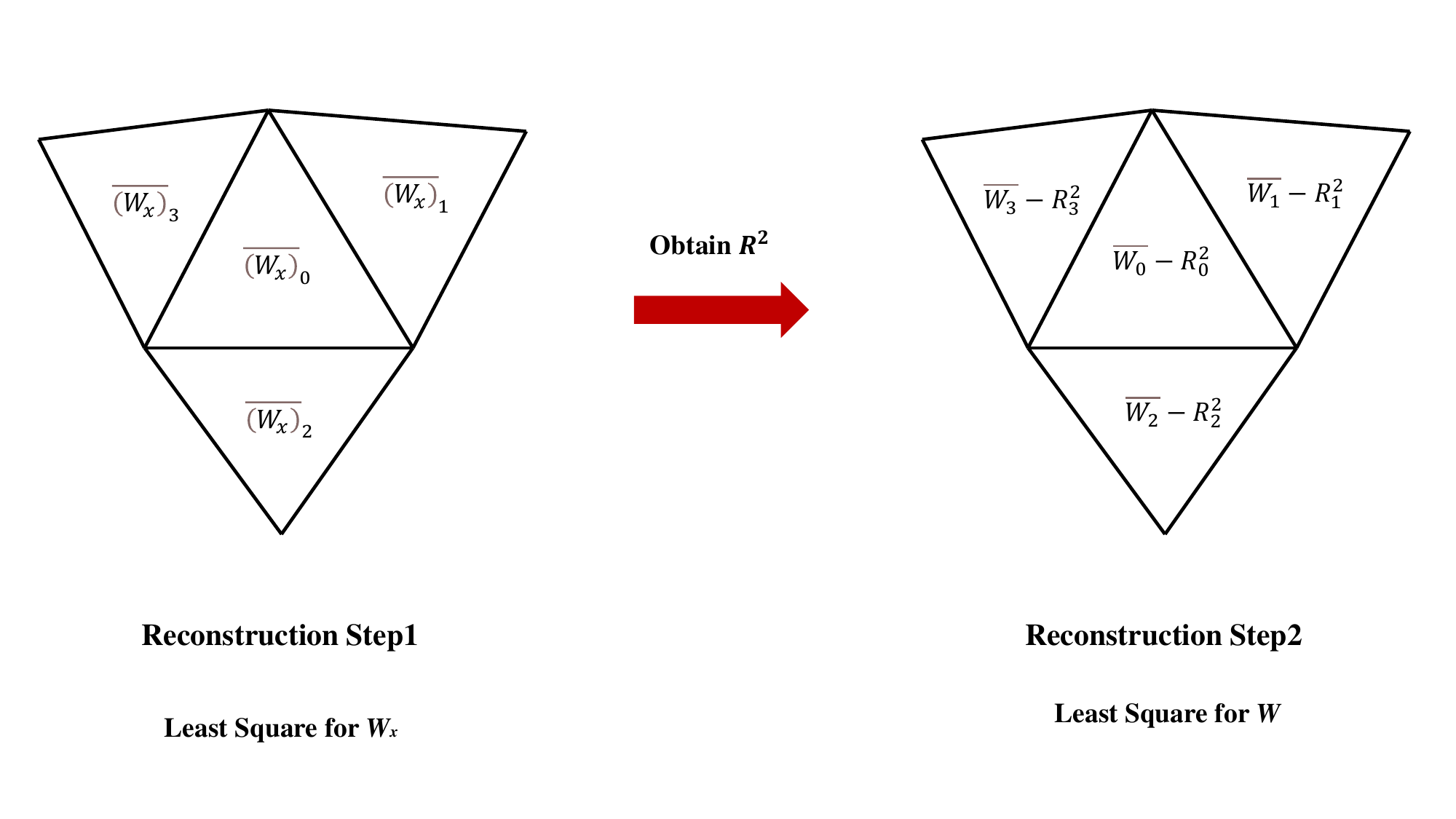}
	\caption{\label{two-step-reconstruction}
		The procedure of memory reduction reconstruction for $p^2$.}
\end{figure}

After completing the above steps, all the coefficients of $p^2$ of the large stencil are solved. Fig.~\ref{two-step-reconstruction} illustrates the procedure of the memory reduction reconstruction.

As can be seen from the above, the main difference between the memory reduction reconstruction and the original reconstruction of $p^2$ is how to use the evolved cell-averaged slopes. 
In the original reconstruction, the cell-averaged slopes are put on the linear system's left-hand side (LFS) and satisfied in a least-square sense. 
In memory reduction reconstruction, the cell-averaged slopes are first used to obtain the coefficients of the quadratic terms of $p^2$. Then, the original linear system's quadratic terms are put on the new linear system's RHS. 

Since the small coefficient matrix with the dimension of $3\times6$ can be constructed in the subroutine temporarily, the new memory reduction reconstruction method can be matrix-free.
Thus, it's very cache-friendly and can speed up the reconstruction. 
In addition to the advantage of compaction and memory reduction, the method is easy to program like traditional second-order FVM, which uses second-order least-square reconstruction.

The memory reduction reconstruction compared to the original reconstruction is shown in Table \ref{reconstruction-compare-table}. 
In this table, the memory consumption is the number of double-precision floating point numbers in each cell, and the matrix assembly time is the time consumption used to assemble the matrix in least-square reconstruction, which is tested using a hexahedron mesh with 50688 cells.

\begin{table}[htp]
	\small
	\begin{center}
		\def\temptablewidth{1.0\textwidth}
		{\rule{\temptablewidth}{1pt}}
		\begin{tabular*}{\temptablewidth}{@{\extracolsep{\fill}}c|c|c|c}
			Reconstruction method & Space order &  Memory consumption  & Matrix assembly time \\
			\hline
			Original & $3$ &  $276$ & 1.3 \\ 	
			Memory reduction reconstruction & $3$ & $60$ & $0$ (matrix-free)\\ 	
		\end{tabular*}
		{\rule{\temptablewidth}{0.1pt}}
	\end{center}
	\vspace{-4mm} \caption{\label{reconstruction-compare-table} 
 Comparison of two reconstruction methods. memory consumption for the reconstruction module(Sixty double-precision floating point numbers in each cell are needed for storing zero-mean basis and the coefficients of the polynomial).}
\end{table}

\subsection{Green-Gauss reconstruction for the sub stencil}
The classical Green-Gauss reconstruction \cite{shima2013green} with only cell-averaged values is adopted to provide the linear polynomial $p^1$ for the sub stencil.
\begin{equation*}
	p^1=\overline{Q}+\boldsymbol{x}\cdot \sum_{m=1}^{N_f}\frac{\overline{Q}_m+\overline{Q}_0}{2}S_m\boldsymbol{n} _m,
\end{equation*}
where $S_m$ is the area of the cell's surface and $\boldsymbol{n}_m$ is the surface's normal vector. In most cases, Green-Gaussian reconstruction has only first-order precision.
\subsection{Discontinuity Feedback}
The DF was first proposed in \cite{ji2021gradient}. Here several improvements have been made in \cite{zhang2023high}: there is no $\epsilon$ in the improved expression of DF; the difference of Mach number is added to improve the robustness under strong rarefaction waves. Denote $\alpha_i \in[0,1]$ as DF at targeted cell $\Omega_i$
\begin{equation*}
	\alpha_i = \prod_{p=1}^m\prod_{k=0}^{M_p}\alpha_{p,k},
\end{equation*}
where $\alpha_{p,k}$ is the CF obtained by the $k$th Gaussian point at the interface $p$ around cell $\Omega_i$, which can be calculated by
\begin{equation*}
	\begin{aligned}
		&\alpha_{p,k}=\frac{1}{1+D^2}, \\
		&D=\frac{|p^l-p^r|}{p^l} +\frac{|p^l-p^r|}{p^r}+(\text{Ma}^{l}_n-\text{Ma}^{r}_n)^2+(\text{Ma}^{l}_t-\text{Ma}^{r}_t)^2,
	\end{aligned}
\end{equation*}
where $p$ is pressure, $\text{Ma}_n$ and $\text{Ma}_t$ are the Mach numbers defined by normal and tangential velocity, and superscripts $l,r$ denote the left and right values of the Gaussian points.

Then, the updated slope is modified by
\begin{equation*}
	\widetilde{\overline{\nabla \boldsymbol{W}}}_i^{n+1} = \alpha_i\overline{\nabla \boldsymbol{W}}_i^{n+1},
\end{equation*}
and the Green-Gauss reconstruction is modified as
\begin{equation*}
	p^1=\overline{Q}+\alpha \boldsymbol{x}\cdot \sum_{m=1}^{N_f}\frac{\overline{Q}_m+\overline{Q}_0}{2}S_m\boldsymbol{n} _m.
\end{equation*}
\subsection{Non-linear WENO weights}
To deal with discontinuity, the idea of multi-resolution WENO reconstruction is adopted \cite{ji2021gradient,zhu2020new}.
Here only two polynomials are chosen
\begin{equation*}
	P_2=\frac{1}{\gamma_2}p^2-\frac{\gamma_1}{\gamma_2}p^1  ,P_1 =p^1.
\end{equation*}
where $\gamma_1=\gamma_2=0.5$. So the quadratic polynomial $p^2$ can be written as
\begin{equation}\label{weno-linear}
	p^2 = \gamma_1P_1 + \gamma_2 P_2.
\end{equation}
Then, we can define the smoothness indicators
\begin{equation*}
	\beta_{j}=\sum_{|\alpha|=1}^{r_{j}}\Omega^{\frac{2}{3}|\alpha|-1} \iiint_{\Omega}\left(D^{\alpha} p^{j}(\mathbf{x})\right)^{2} \mathrm{~d} V,
\end{equation*}
where $\alpha$ is a multi-index and $D$ is the derivative operator, $r_1=1,r_2=2$. Special care is given for $\beta_1$ for better robustness
\begin{equation*}
	\beta_1=\min(\beta_{1,\text{Green-Gauss}},\beta_{1,\text{least-square}}),
\end{equation*}
where $\beta_{1,\text{Green-Gauss}}$ is the smoothness indicator defined by Green-Gauss reconstruction, and $\beta_{1,\text{least-square}}$ is the smoothness indicator defined by second-order least-square reconstruction. Then, the smoothness indicators $\beta_i$ are non-dimensionalized by
\begin{equation*}
	\tilde{\beta}_i=\frac{\beta_i}{Q_0^2+\beta_1+10^{-40}}.
\end{equation*}
The nondimensionalized global smoothness indicator $\tilde{\sigma}$ can be defined as
\begin{equation*}
	\tilde{\sigma}=\left|\tilde{\beta}_{1}-\tilde{\beta}_{0}\right|^{}.
\end{equation*}
Therefore, the corresponding non-linear weights are given by
\begin{equation*}
	\tilde{\omega}_{m}=\gamma_{m}\left(1+\left(\frac{\tilde{\sigma}}{\epsilon+\tilde{\beta}_{m}}\right)^{2}\right),\epsilon=10^{-5},
\end{equation*}
\begin{equation*}
	\bar{\omega}_{m}=\frac{\tilde{\omega}_{m}}{\sum \tilde{\omega}_{m}}, m=1,2.
\end{equation*}
Replacing $\gamma_m$ in equation (\ref{weno-linear}) by $\bar{\omega}_{m}$ , the final non-linear reconstruction can be obtained
\begin{equation*}
	R(\boldsymbol{x})=\bar{\omega}_2P_2+\bar{\omega}_1 P_1.
\end{equation*}
The desired non-equilibrium states at Gaussian points become
\begin{equation*}
	Q_{p, k}^{l, r}=R^{l, r}\left(\boldsymbol{x}_{p, k}\right), \left(Q_{x_{i}}^{l, r}\right)_{p, k}=\frac{\partial R^{l, r}}{\partial x_{i}}\left(\boldsymbol{x}_{p, k}\right).
\end{equation*}

\section{Numerical examples} \label{test-case}
In this section, we present numerical tests to validate the proposed scheme. 
To achieve high-order accuracy in time advance, the two-stage fourth-order (S2O4) time discretization is adopted, the details can be found in \cite{zhao2021direct}.
All simulations are conducted on a three-dimensional hybrid unstructured mesh. 
The hybrid meshes ensure flexibility and high resolution, allowing accurate capture of complex geometries and flow features.
The tests include benchmark problems and more practical applications, focusing on metrics such as computational cost, and error analysis. 
The results demonstrate the scheme's effectiveness and robustness, providing a comprehensive validation against existing methods.
The simulations are conducted by our in-house C++ solver, where MPI is used for parallel computation, and METIS is used for mesh partitioning.

\begin{figure}[htp]
	\centering
	\includegraphics[width=0.8\textwidth]
	{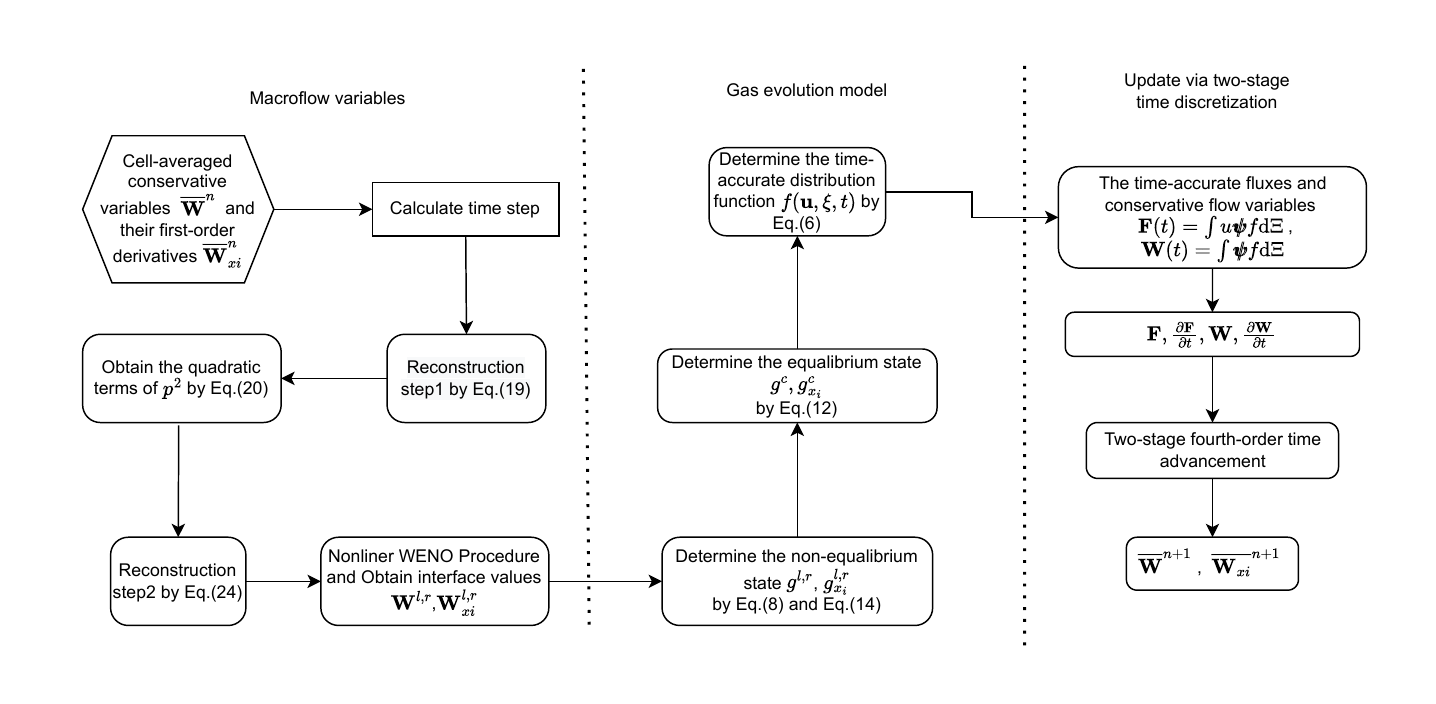}
	\caption{\label{memory-reduction-procedure}
		The framework of memory reduction CGKS. }
\end{figure}

A brief flowchart of the whole memory-reduction CGKS is shown in Fig.~\ref{memory-reduction-procedure}.

\subsection{Accuracy Test}
In this test case, 3-D sinusoidal wave propagation is calculated to verify the accuracy of the scheme.
The initial condition for the advection of density per
perturbation is given as
\begin{equation*}
\begin{aligned}
& \rho(x, y, z)=1+0.2 \sin (\pi(x+y+z)), \\
& \mathbf{U}(x, y, z)=(1,1,1), \quad p(x, y, z)=1.
\end{aligned}
\end{equation*}

\begin{figure}[htp]
	\centering	
	\includegraphics[height=0.35\textwidth]{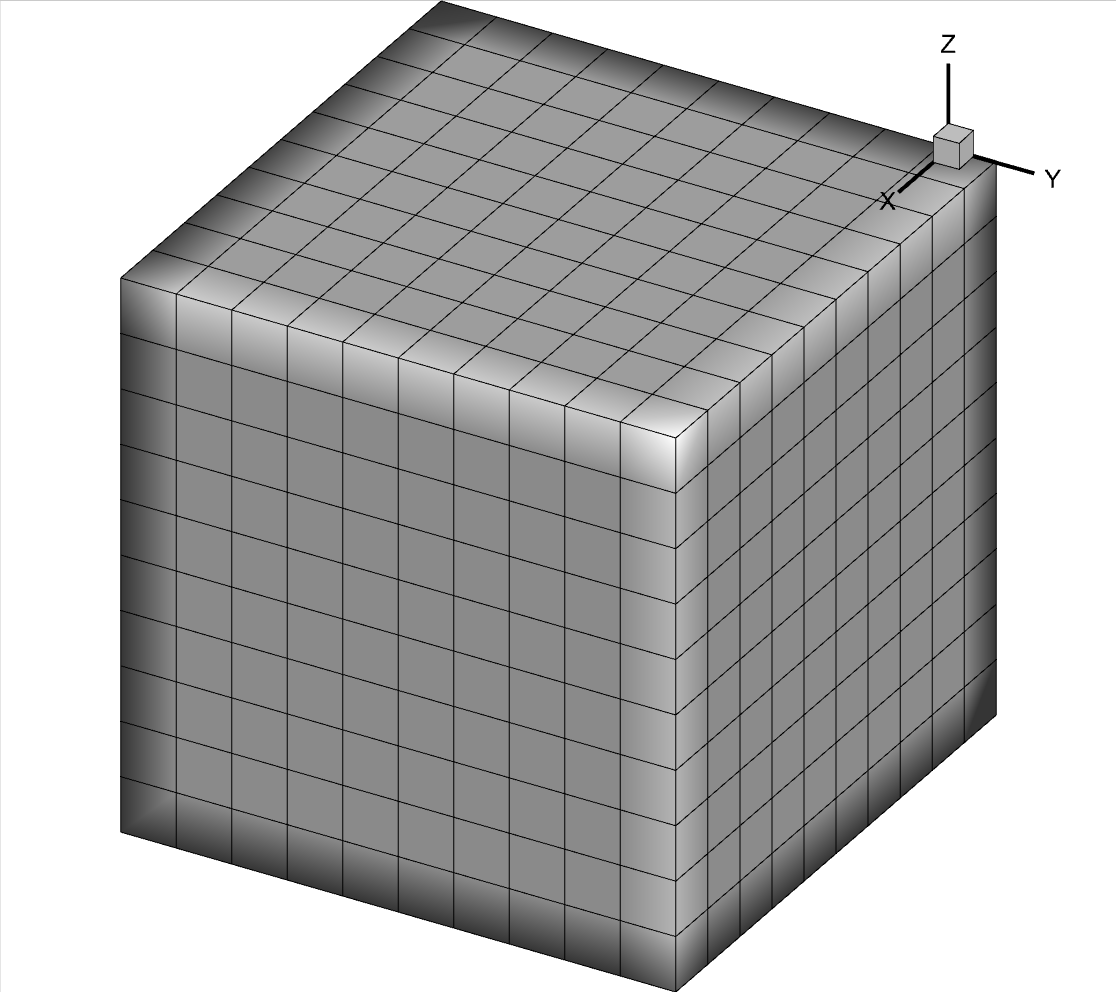}
	\includegraphics[height=0.35\textwidth]{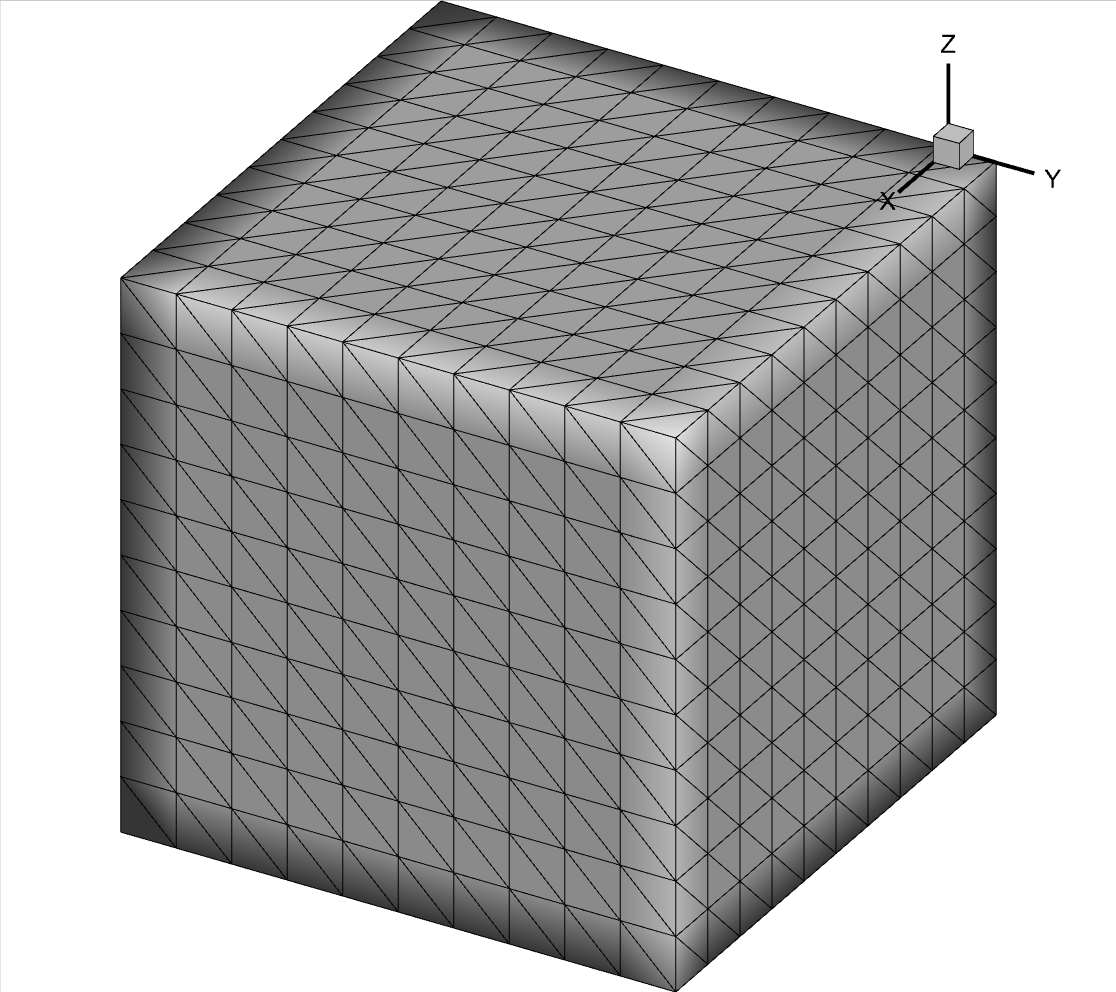}
	\caption{\label{sin-wave-mesh}
		Meshes used in accuracy test. Left: hexahedron mesh. Right: tetrahedron mesh.}
\end{figure}

The domain is cubic and its size is $[0,2] \times[0,2] \times[0,2]$.  
A series of sequentially refined hexahedron meshes and tetrahedron meshes are used in this test case, as shown in Fig.~\ref{sin-wave-mesh}. 
With the periodic boundary condition in all directions, the analytic solution is

\begin{equation*}
\begin{aligned}
& \rho(x, y, z, t)=1+0.2 \sin (\pi(x+y+z-t)) \\
& \mathbf{U}(x, y, z)=(1,1,1), \quad p(x, y, z, t)=1 .
\end{aligned}
\end{equation*}

\begin{figure}[htp]
	\centering	
	\includegraphics[height=0.35\textwidth]{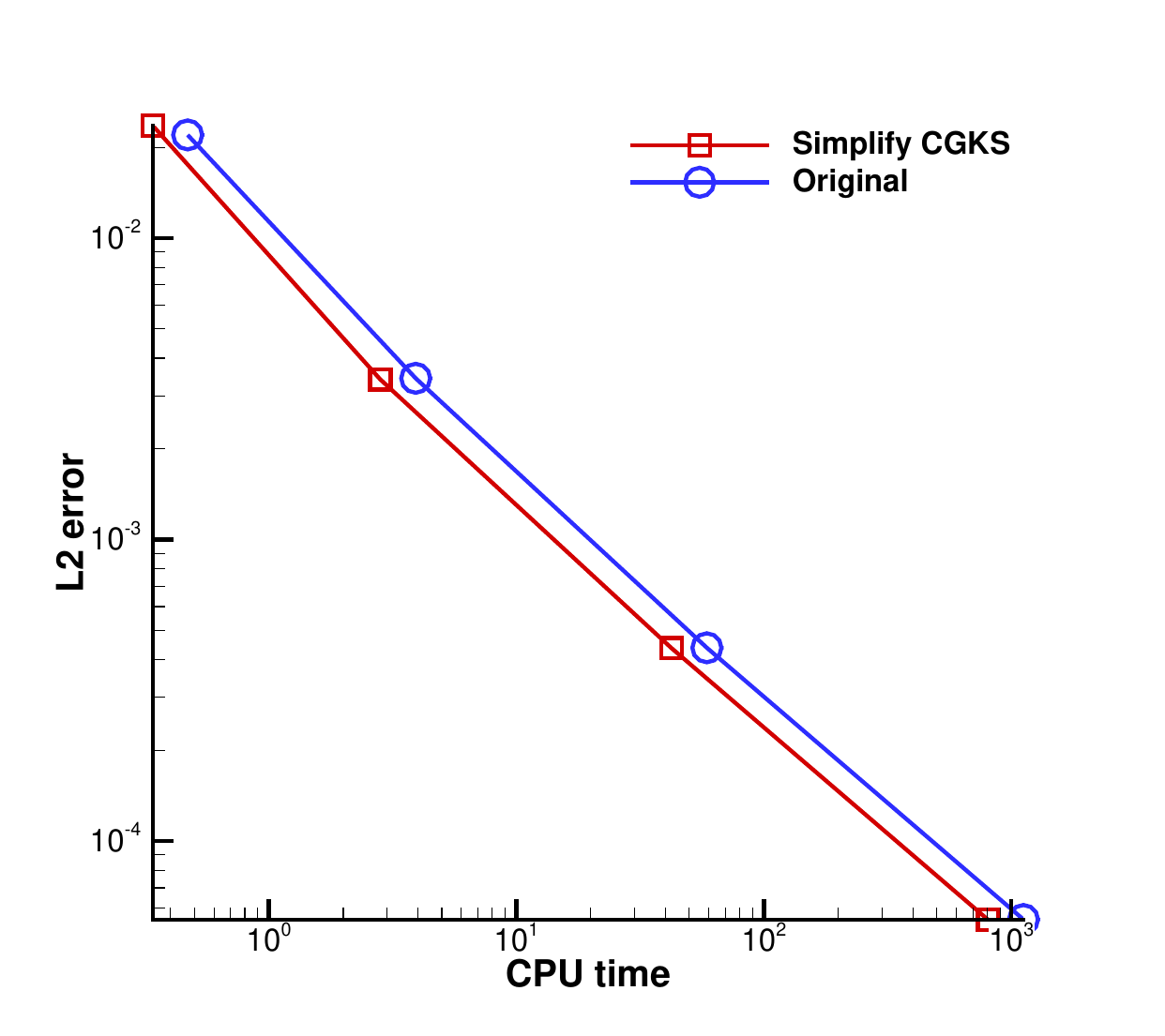}
	\includegraphics[height=0.35\textwidth]{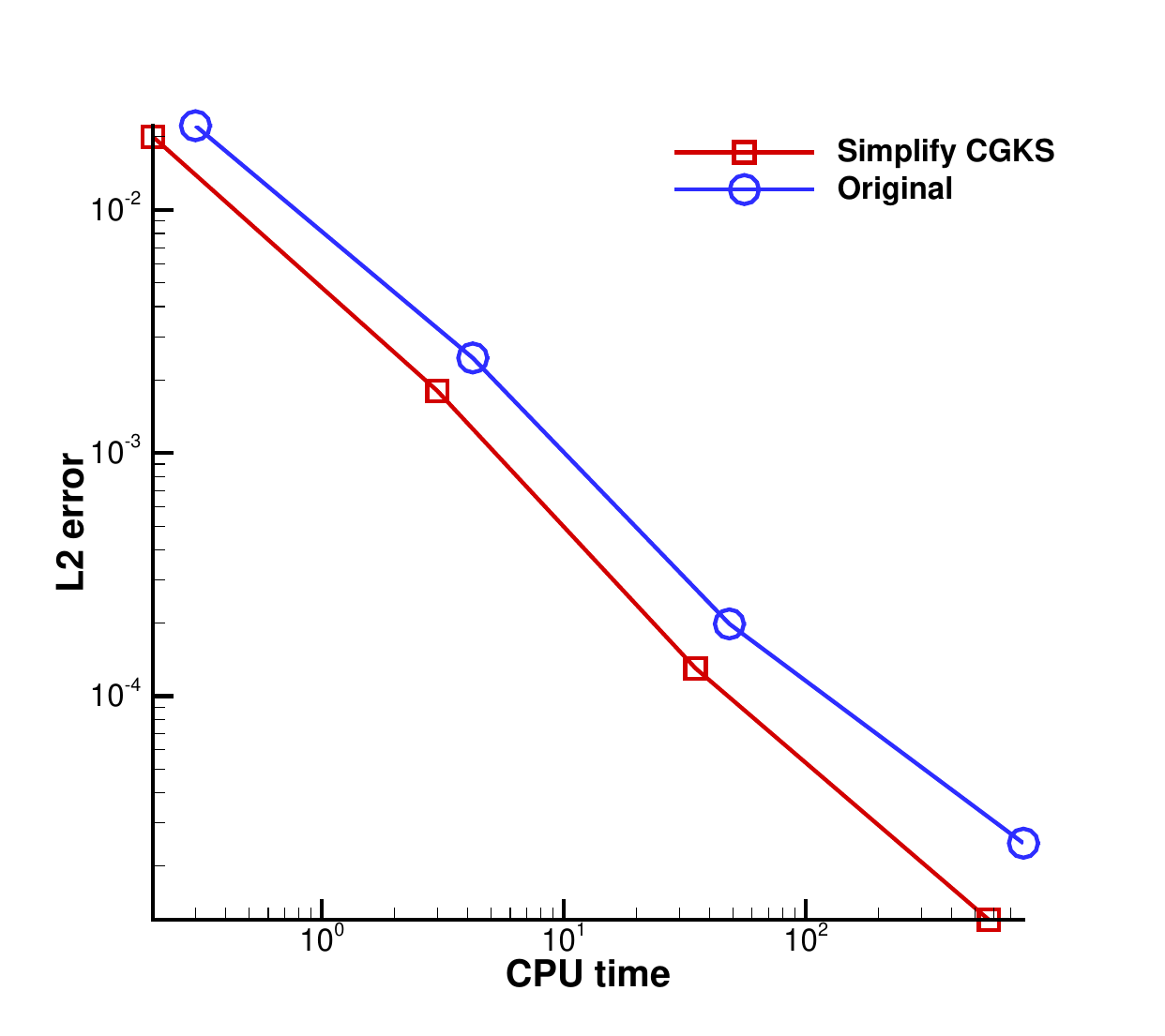}
	\caption{\label{sin-wave-cputime-compare}
		Original CGKS versus Simplify CGKS of CPU time. Left: hexahedron mesh. Right: tetrahedron mesh.}
\end{figure}

The flow is inviscid and the collision time $\tau$ is 0. The $L^1, L^2$ and $L^{\infty}$ errors and the corresponding orders with linear weights at $t=2$ under both meshes are given in Table \ref{sin-hexa} and Table \ref{sin-tet}. 
Expected accuracy is achieved for all cases.
Meanwhile, the original CGKS versus Simplify CGKS of CPU time is shown in Fig.~\ref{sin-wave-cputime-compare}. The overall $28\%$ efficiency improvements have been achieved in all cases.

\begin{table}[htp]
	\small
	\begin{center}
		\def\temptablewidth{1.0\textwidth}
		{\rule{\temptablewidth}{1pt}}
		\begin{tabular*}{\temptablewidth}{@{\extracolsep{\fill}}c|c|c|c|c|c|c}
			Mesh number & $L^{1}$ error &  order & $L^{2}$ error & order & $L^{\infty}$ error & order\\
			\hline
			$10^3$ & $2.147907e^{-2}$ & -- & $2.365047e^{-2}$ & -- &  $3.353870e^{-2}$ & --\\ 	
			$20^3$ & $3.064556e^{-3}$ & 2.81 & $3.401127e^{-3}$ & 2.80 &  $4.998815e^{-3}$ & 2.75\\ 	
			$40^3$ & $3.933715e^{-4}$ & 2.96 & $4.367541e^{-4}$ & 2.96 &  $6.561719e^{-4}$ & 2.93\\ 	
			$80^3$ & $4.952024e^{-5}$ & 2.99 & $5.486601e^{-5}$ & 2.99 &  $8.282316e^{-5}$ & 2.99\\ 	
		\end{tabular*}
		{\rule{\temptablewidth}{0.1pt}}
	\end{center}
	\vspace{-4mm} \caption{\label{sin-hexa} Accuracy test using hexahedron mesh.}
\end{table}

\begin{table}[htp]
	\small
	\begin{center}
		\def\temptablewidth{1.0\textwidth}
		{\rule{\temptablewidth}{1pt}}
		\begin{tabular*}{\temptablewidth}{@{\extracolsep{\fill}}c|c|c|c|c|c|c}
			Mesh number & $L^{1}$ error &  order & $L^{2}$ error & order & $L^{\infty}$ error & order\\
			\hline
			$6 \times 5^3$ & $1.794963e^{-2}$ & -- & $2.017067e^{-2}$ & -- &  $3.288013e^{-2}$ & --\\ 	
			$6 \times 10^3$ & $1.560746e^{-3}$ & 3.52 & $1.766208e^{-3}$ & 3.51 &  $3.297513e^{-3}$ & 3.32\\ 	
			$6 \times 20^3$ & $1.135311e^{-4}$ & 3.78 & $1.318874e^{-4}$ & 3.74 &  $2.839491e^{-4}$ & 3.54\\ 	
			$6 \times 40^3$ & $9.508264e^{-6}$ & 3.58 & $1.165079e^{-5}$ & 3.50 &  $2.797772e^{-5}$ & 3.34\\ 	
		\end{tabular*}
		{\rule{\temptablewidth}{0.1pt}}
	\end{center}
	\vspace{-4mm} \caption{\label{sin-tet} Accuracy test using tetrahedron mesh.}
\end{table}

\subsection{Subsonic flow around a cylinder}
\begin{figure}[htp]
	\centering
	\includegraphics[width=0.45\textwidth]
	{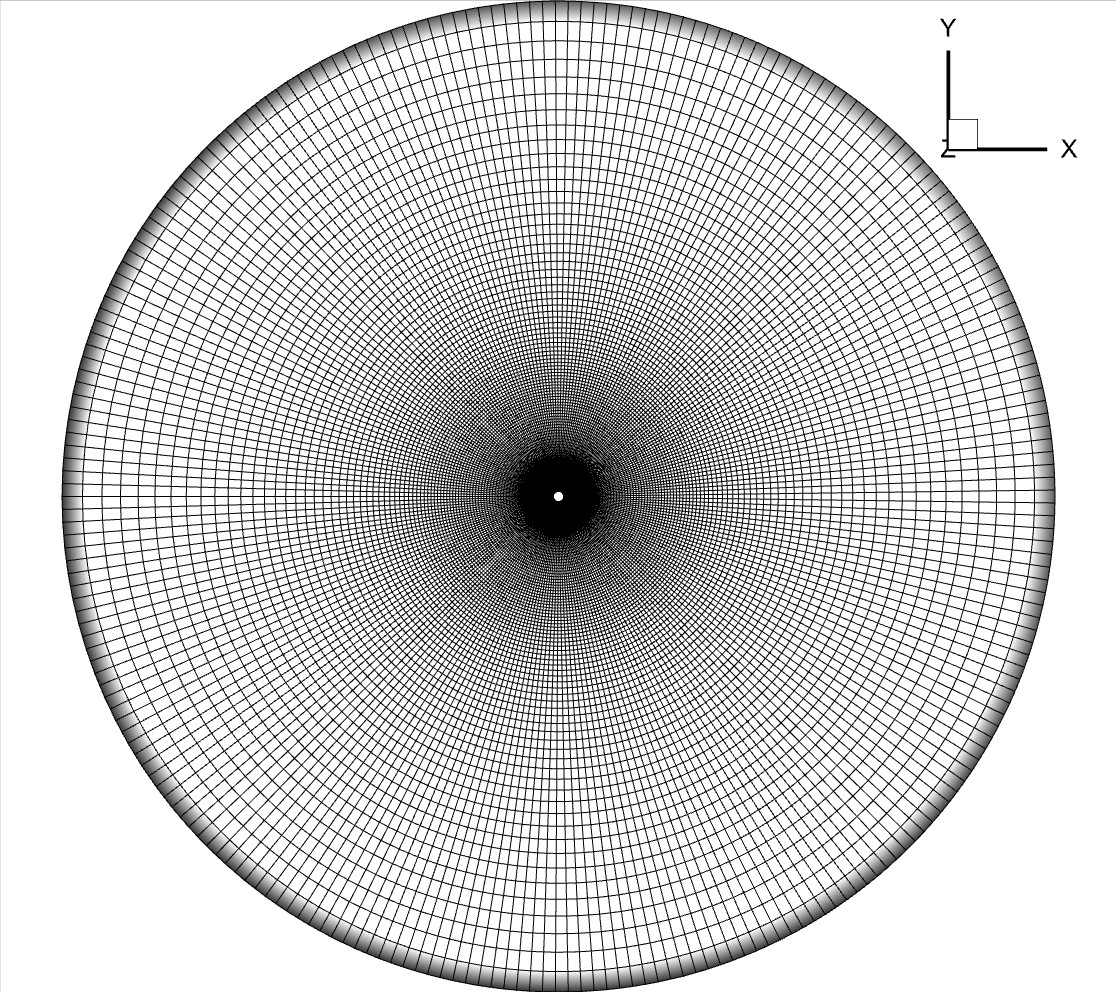}
        \includegraphics[width=0.45\textwidth]
	{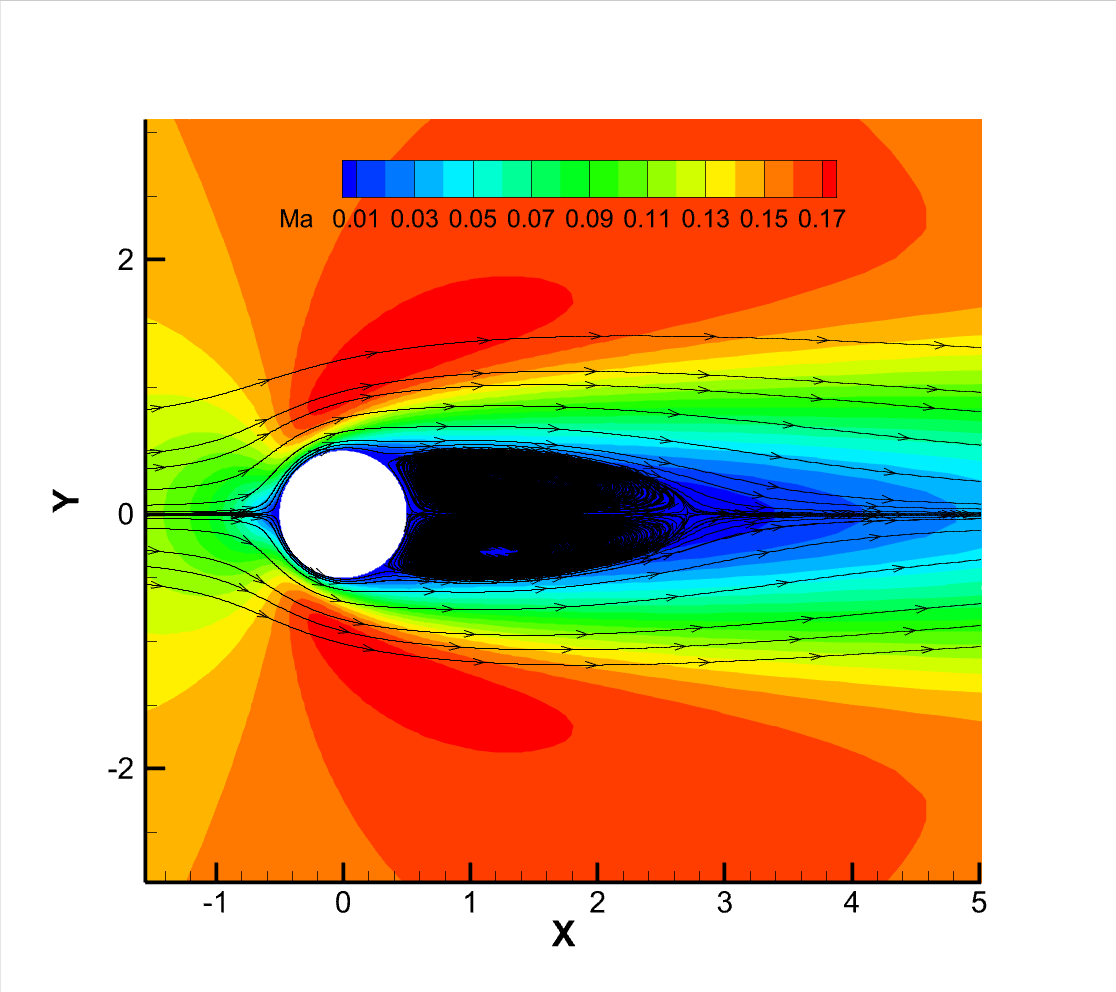}
	\caption{\label{cylinder-result}
		Subsonic flow around a cylinder. Left: hexahedron mesh. Right: Mach number contour with streamlines.}
\end{figure}

In this case, subsonic flow around a cylinder is simulated. 
The far-field incoming flow condition is Mach number equals to 0.15 and Reynolds number equals to 40.
The diameter of the cylinder is 1, and the diameter of the whole computational domain is 96.
A total hexahedron-type mesh with a mesh number equal to 9450 is used in the simulation, and the near-wall size is $h=1/96$. 
Two stable and symmetrical vortices appear at the cylindrical tail.
The mesh and the Mach number contour with streamlines are shown in Fig.~\ref{cylinder-result}.

Quantitative results include the drag coefficient $C_D$, the lift coefficient $C_L$, the wake length $L$, and the separation angle $\theta$ are listed in Table \ref{cylinder-table}.
The results above show the current scheme agrees well with the experimental and numerical references.

\begin{table}[htp]
	\small
	\begin{center}
		\def\temptablewidth{1.0\textwidth}
		{\rule{\temptablewidth}{1pt}}
		\begin{tabular*}{\temptablewidth}{@{\extracolsep{\fill}}c|c|c|c|c|c|c}
			Method & $C_D$ &  $C_L$ & $L$ &  Vortex Height & Vortex Width & $\theta$\\
			\hline
			Experiment \cite{tritton1959experiments} & $1.46-1.56$ & -- & -- & -- &  -- & --\\ 	
			Experiment \cite{coutanceau1977experimental} & -- & -- & 2.12 & 0.297 &  0.751 & $53.5^{\circ}$\\ 	
			DDG \cite{zhang2019direct} & 1.529 & -- & 2.31 & -- & -- & --\\ 	
			Current & 1.527 & $5.8e^{-14}$ & 2.22 & 0.296 &  0.714 & $53.5^{\circ}$\\ 	
		\end{tabular*}
		{\rule{\temptablewidth}{0.1pt}}
	\end{center}
	\vspace{-4mm} \caption{\label{cylinder-table} Comparison of the quantitative results of subsonic flow around a cylinder.}
\end{table}

\subsection{Subsonic flow around a NACA0012 airfoil}

\begin{figure}[htp]	
	\centering	
	\includegraphics[height=0.35\textwidth]{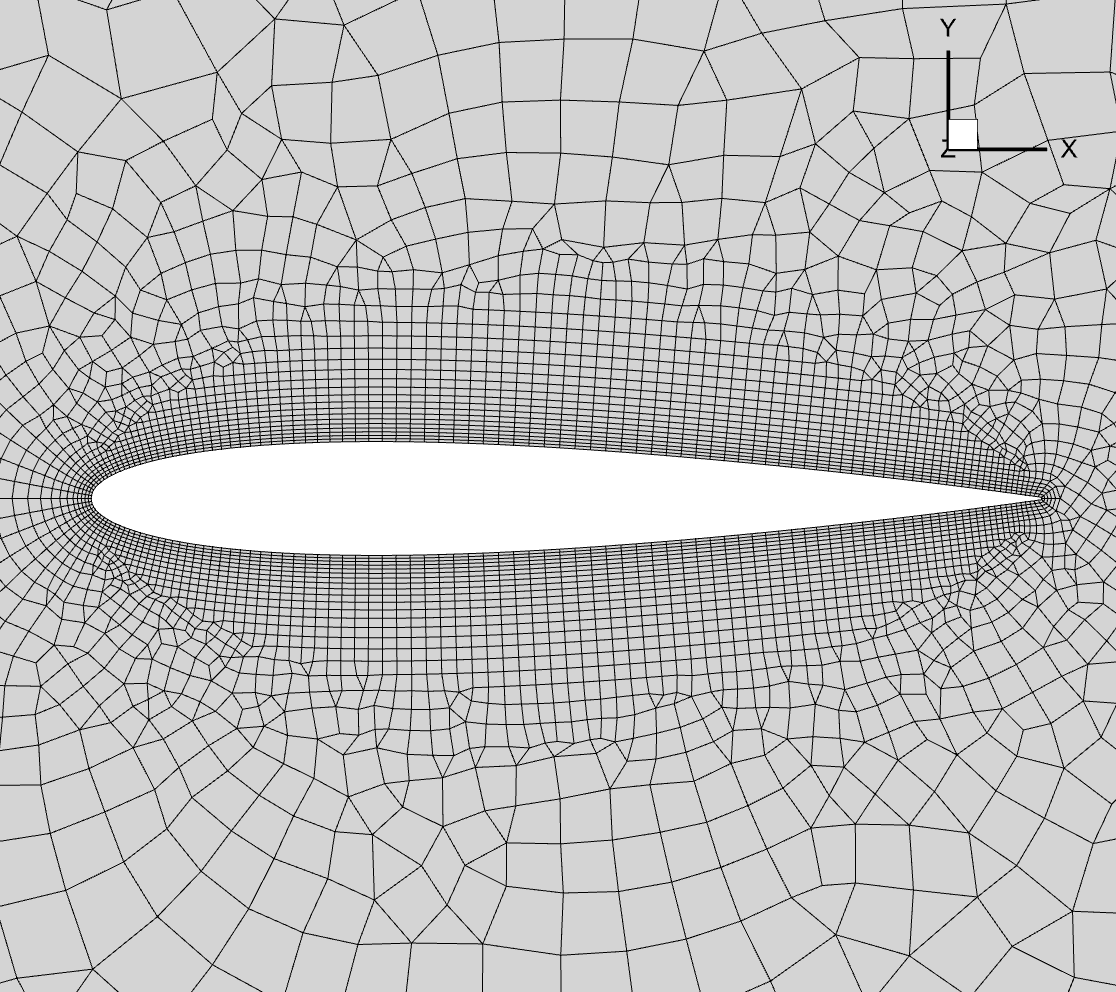}
	\includegraphics[height=0.35\textwidth]{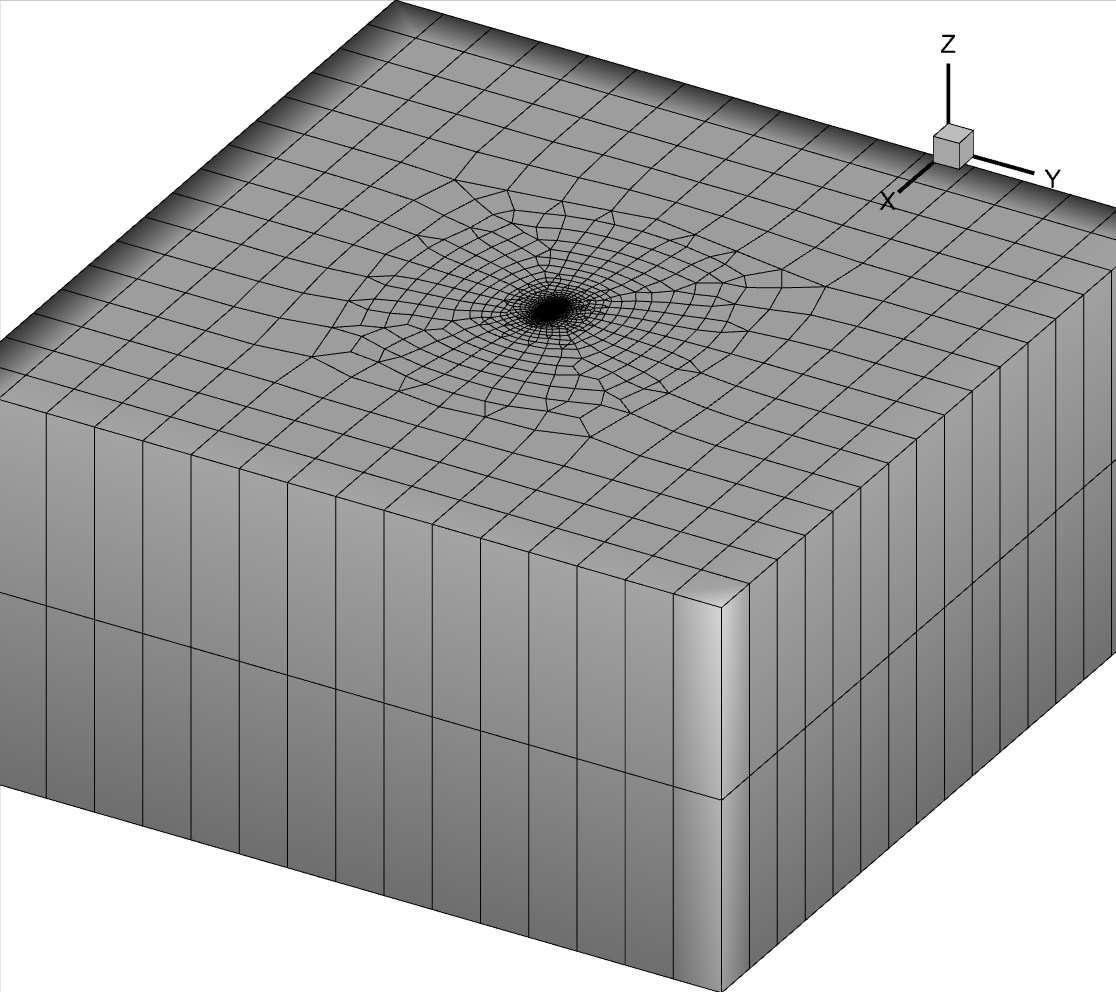}
	\caption{\label{naca0012 mesh}
		NACA0012 airfoil Mesh.}
\end{figure}

\begin{figure}[htp]	
	\centering	
	\includegraphics[height=0.35\textwidth]{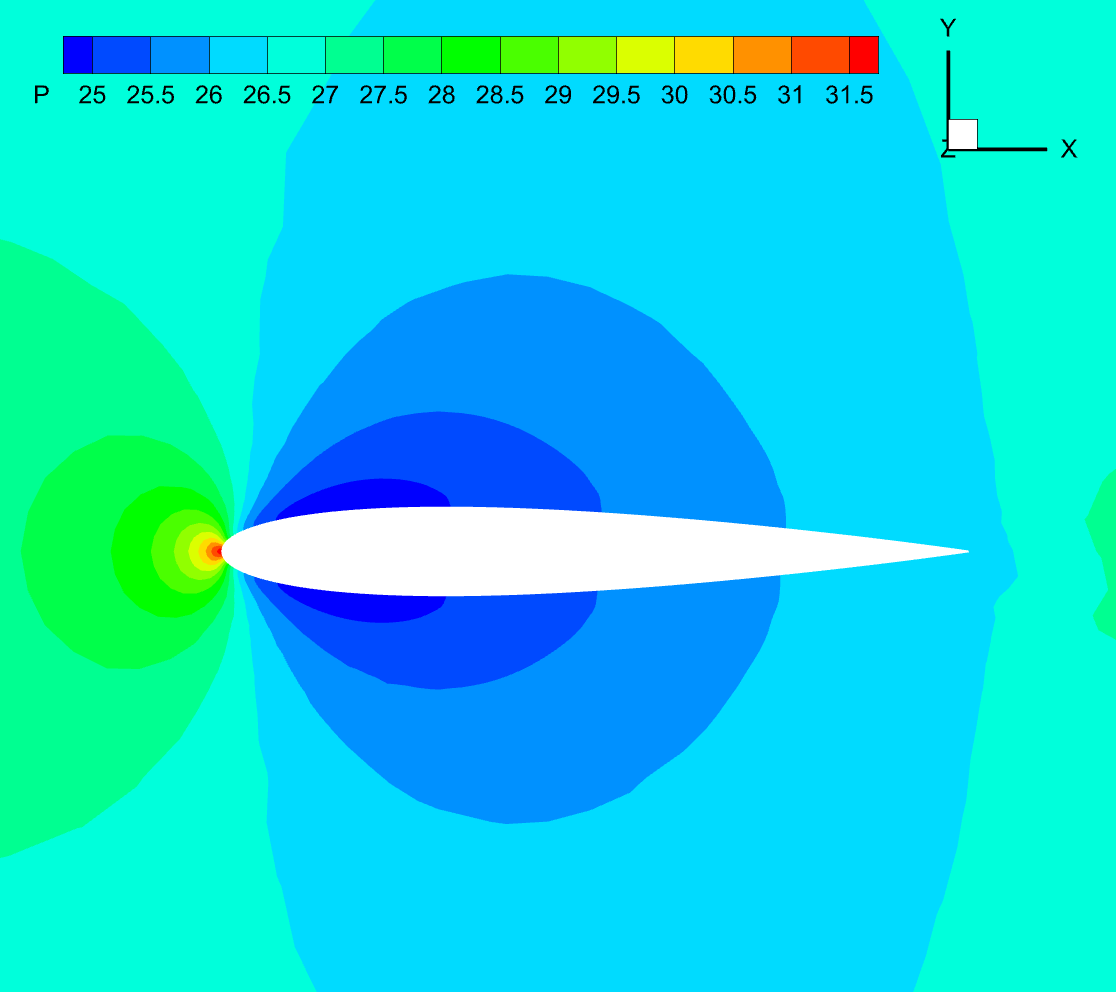}
	\includegraphics[height=0.35\textwidth]{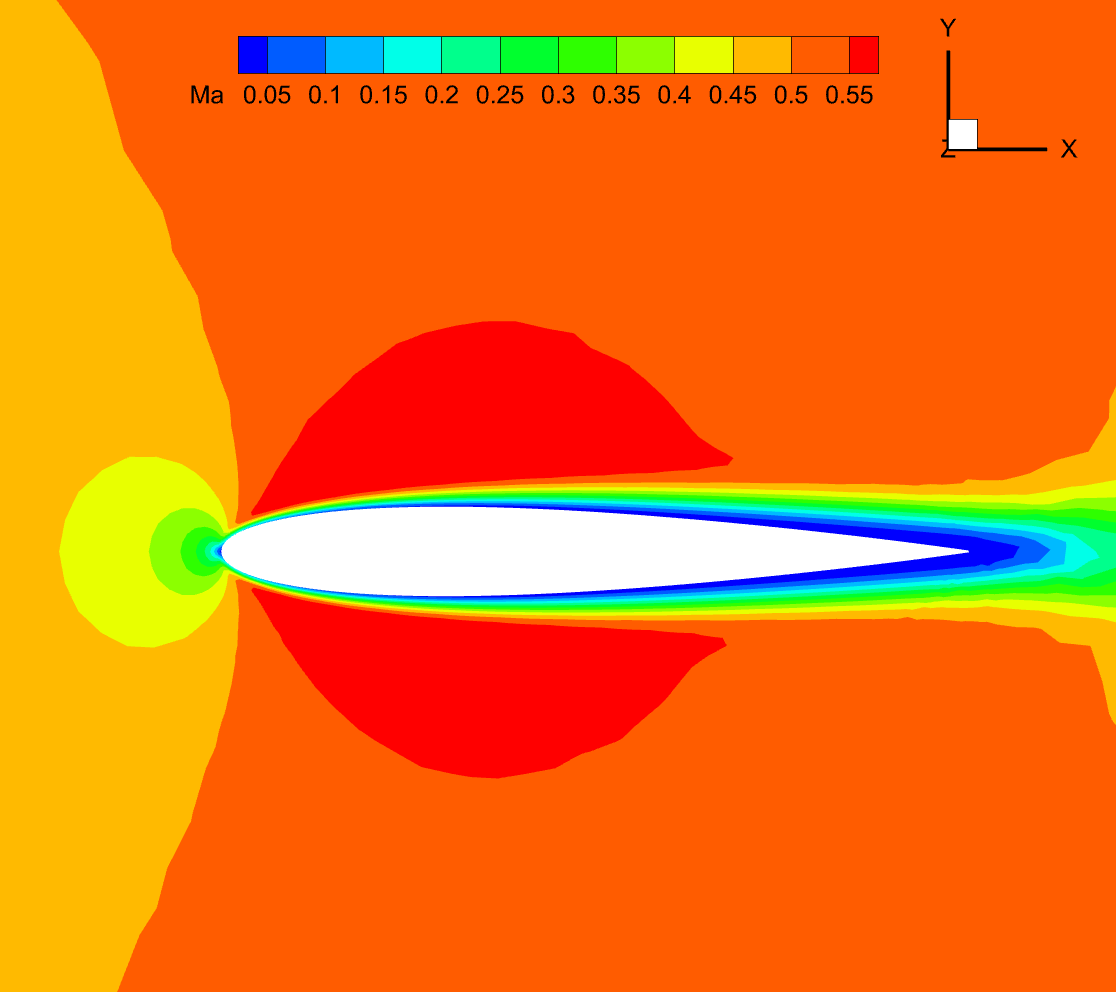}
	\caption{\label{naca0012 result}
		Subsonic flow around a NACA0012 airfoil. Left: Pressure contour. Right: Mach number contour.}
\end{figure}

\begin{figure}[htp]	
	\centering	
	\includegraphics[height=0.35\textwidth]{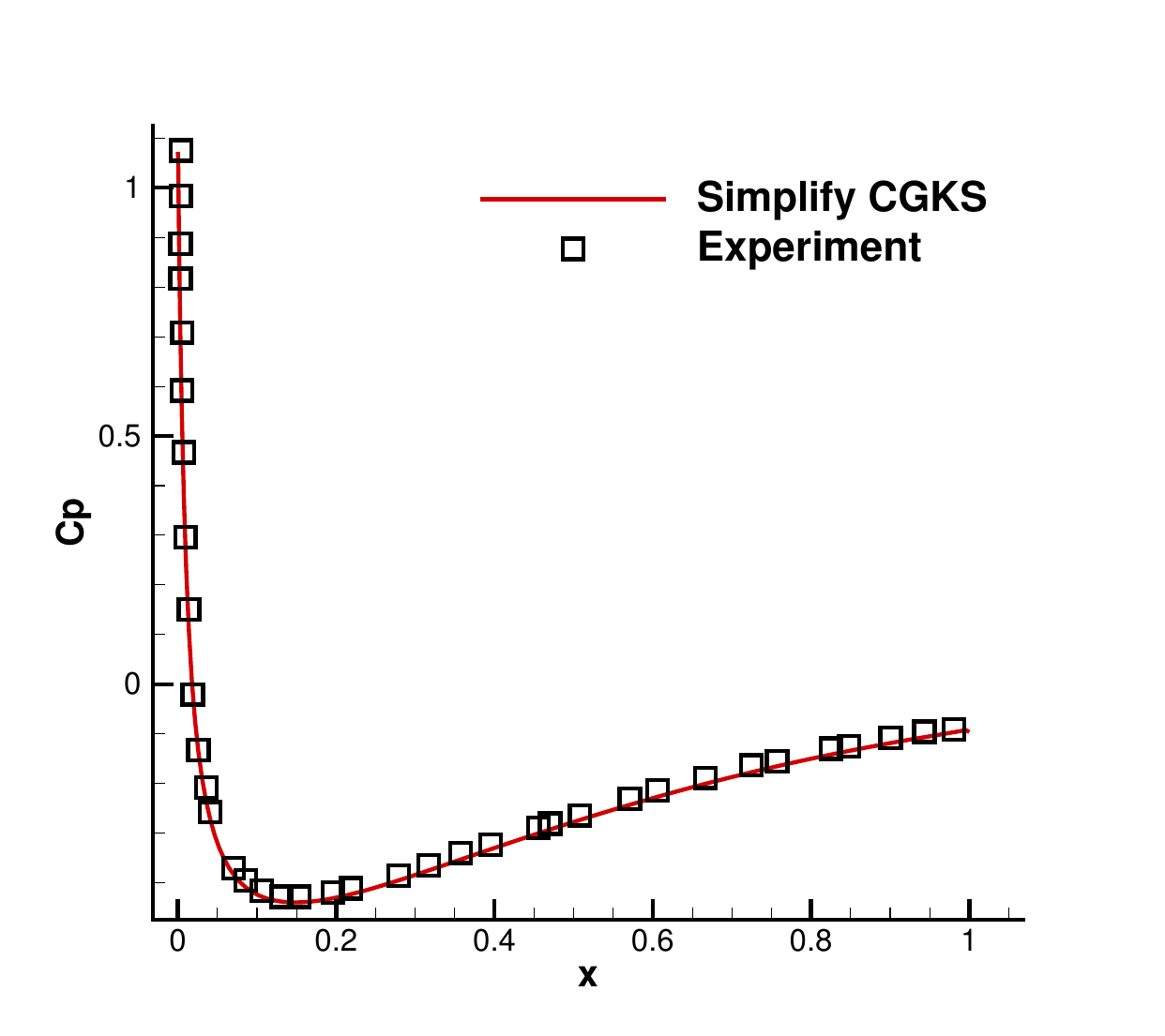}
	\caption{\label{naca0012 cp}
		Subsonic flow around a NACA0012 airfoil. Surface pressure coefficient distribution.}
\end{figure}

In this section, viscous flow around a NACA0012 airfoil is simulated. 
The incoming Mach number is set to 0.5 and the incoming Reynolds number is set to 5000 based on the chord length L=1. 
The subsonic far-field is calculated by Riemann invariants and the solid wall of the airfoil is set to be an adiabatic non-slip wall. 
Total 6538 $\times$ 2 hybrid prismatic cells are used in a cuboid domain [-15, 15] $\times$ [15, 15] $\times$ [0, 0.1]. 
The hybrid unstructured mesh is shown in Fig.~\ref{naca0012 mesh}.  
The Mach number contour is shown in Fig.~\ref{naca0012 result}. Quantitative result including the surface pressure coefficient is extracted and plotted in Fig.~\ref{naca0012 cp}, which highly agrees with the Ref \cite{bassi1997ns}.

\subsection{Transonic flow around dual NACA0012 airfoils}
\begin{figure}[htp]	
	\centering	
	\includegraphics[height=0.35\textwidth]{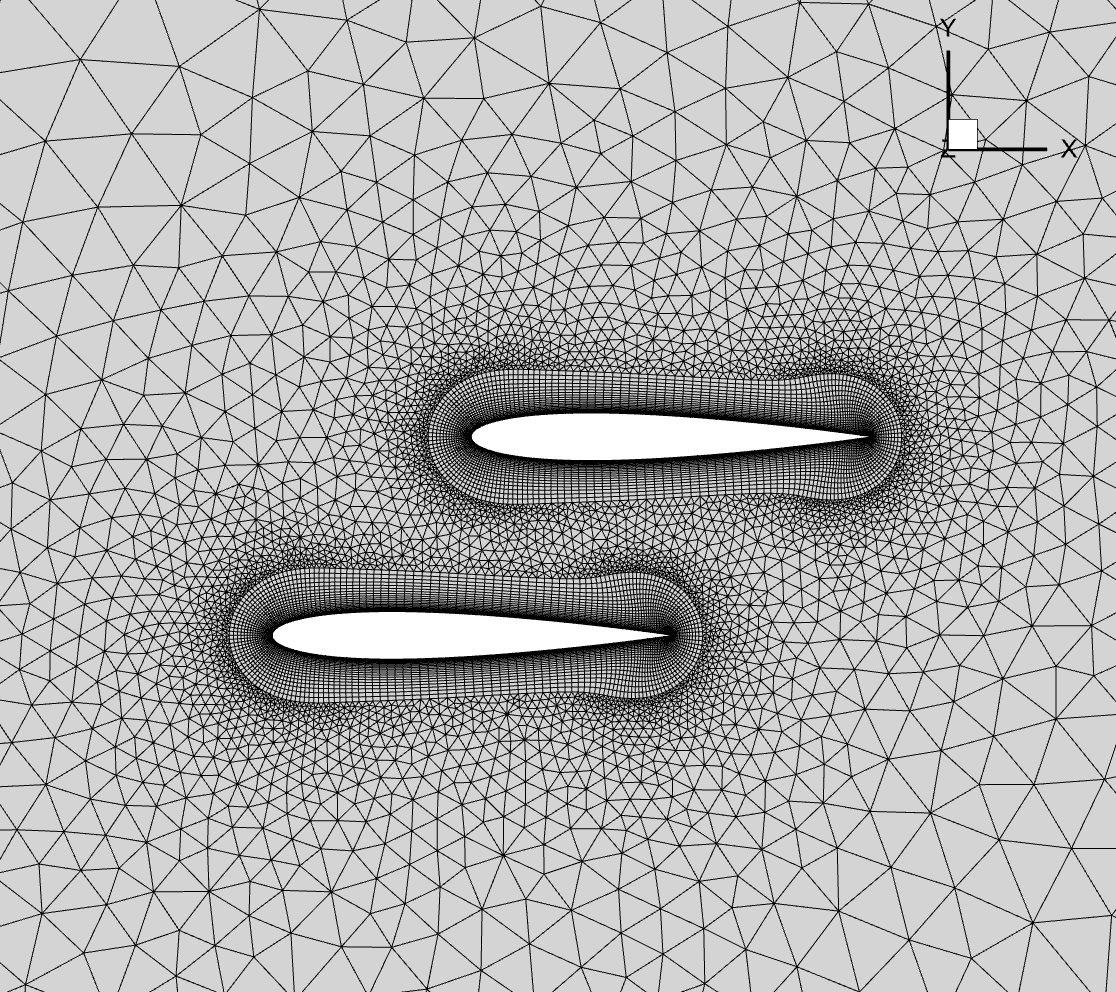}
	\includegraphics[height=0.35\textwidth]{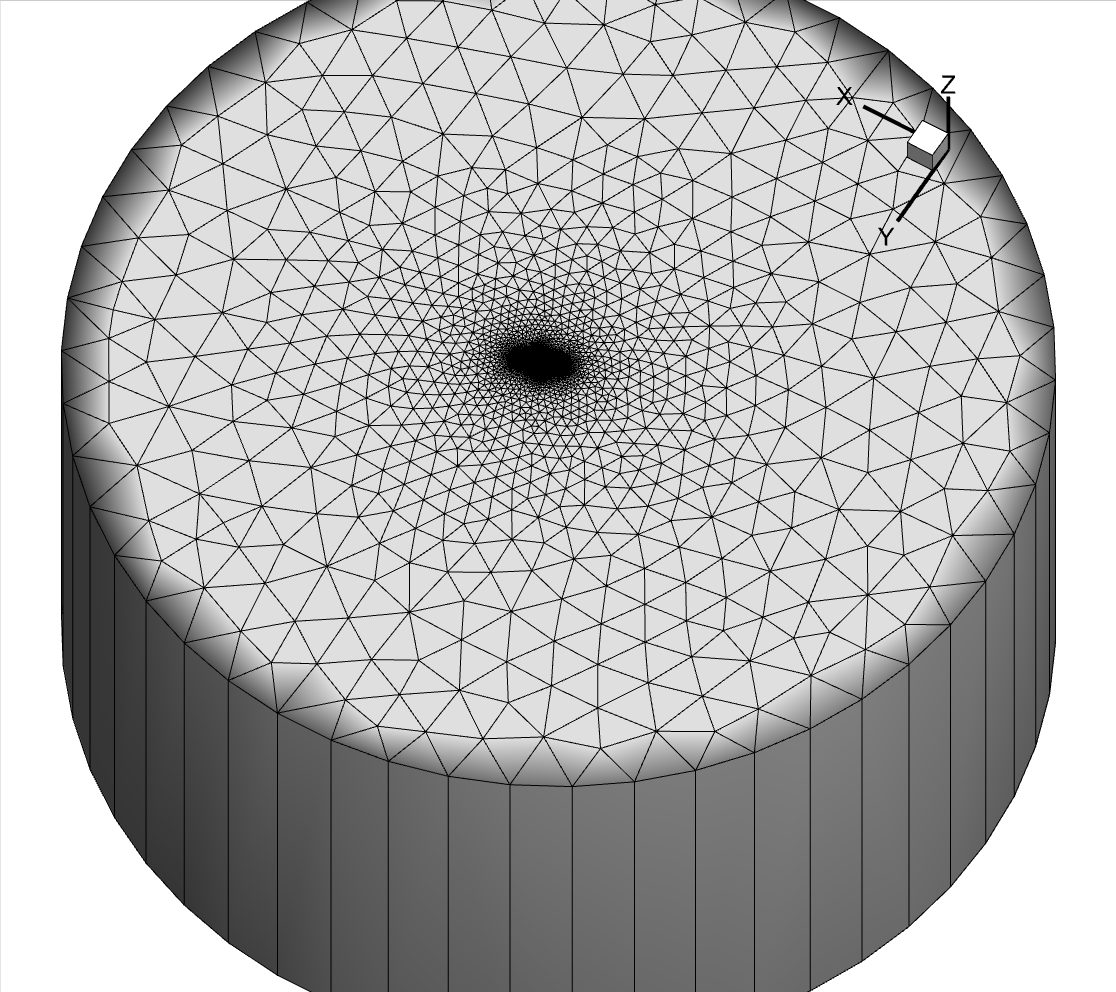}
	\caption{\label{dual naca0012 mesh}
		Transonic flow around dual NACA0012 airfoils. Mesh.}
\end{figure}

\begin{figure}[htp]	
	\centering	
	\includegraphics[height=0.35\textwidth]{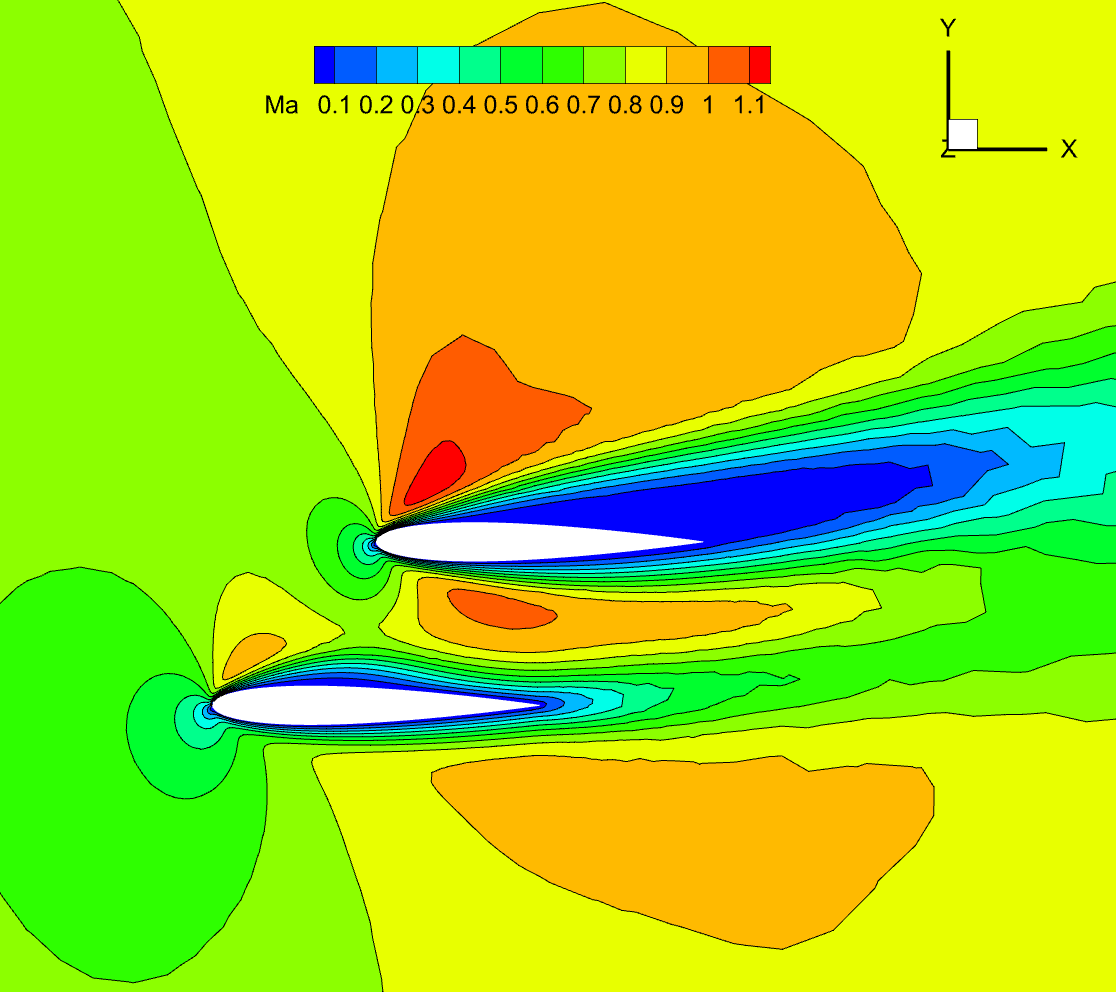}
	\includegraphics[height=0.35\textwidth]{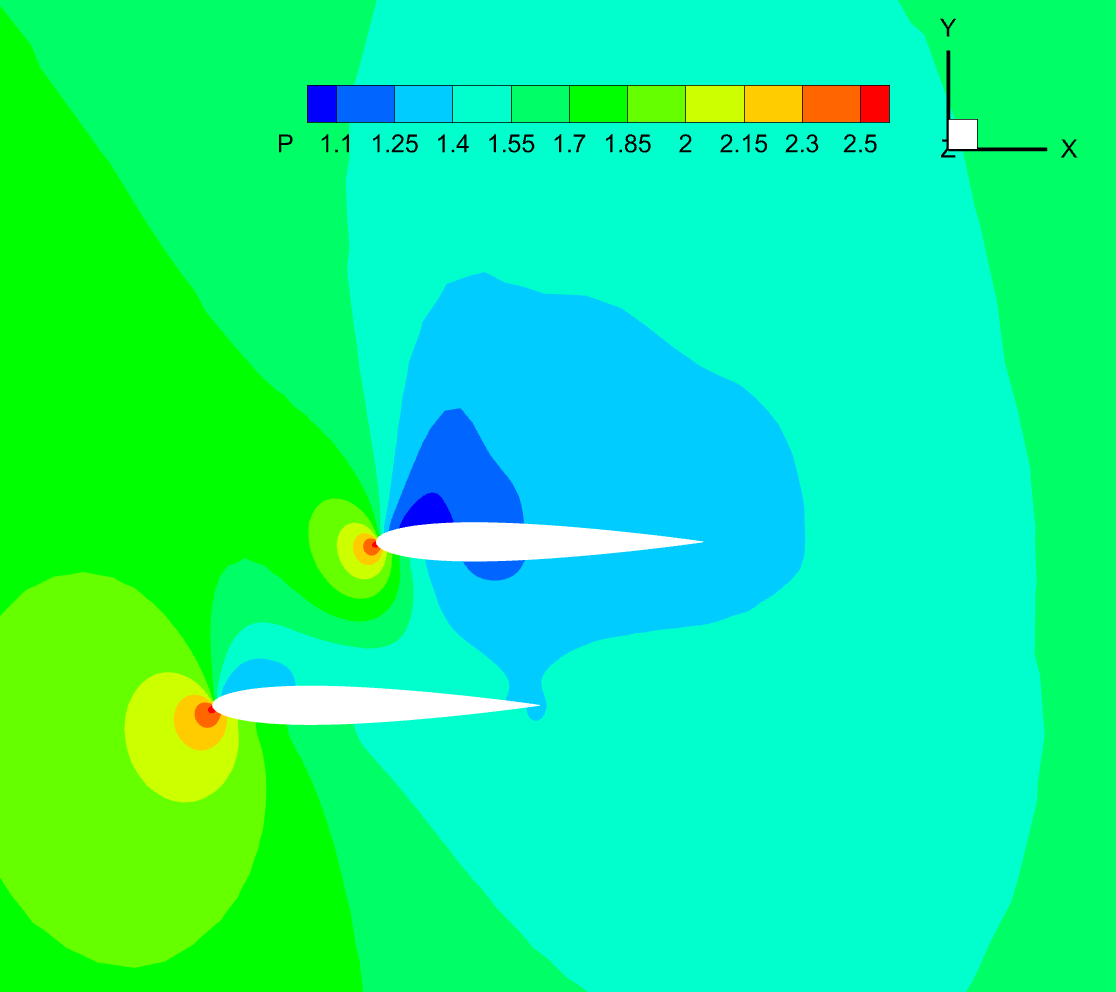}
	\caption{\label{dual naca0012 result}
		Transonic flow around dual NACA0012 airfoils. Left: Mach number distribution. Right: Pressure distribution.}
\end{figure}

\begin{figure}[htp]	
	\centering	
	\includegraphics[height=0.35\textwidth]{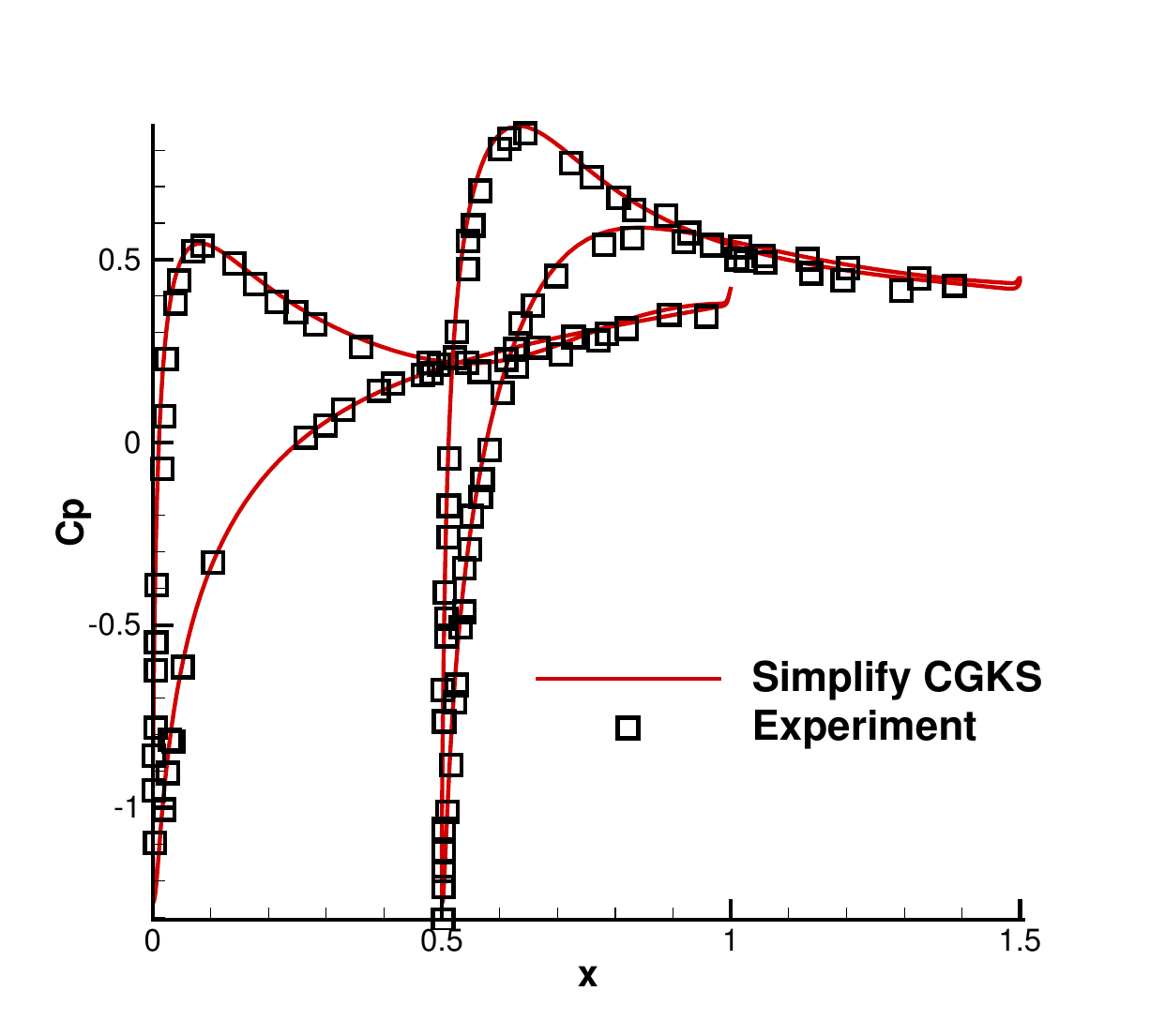}
	\caption{\label{dual naca0012 cp}
		Transonic flow around dual NACA0012 airfoils. Surface pressure coefficient distribution.}
\end{figure}

To verify the memory reduction CGKS on a more complicate case, transonic flow around dual NACA0012 airfoils is simulated. The head of the first airfoil is located at (0, 0) and the second one is located at (0.5, 0.5). 
Both airfoils are put in parallel with the x-axis. 
The incoming Mach number is set to be 0.8 with an angle of attack AOA = 10 $^{\circ}$ and the Reynolds number is set to be 500 based on the chord length L =1. The mesh consists of 28678 mixed elements. The near wall size of the mesh is set to be h = 2 $\times$ $10^{-3}$, which indicates that the grid Reynolds number is 2.5 $\times$ $10^{5}$. 
The far-field boundary condition is set to be subsonic inflow using Riemann invariants and the wall is set to be a non-slip adiabatic wall. 
The mesh is presented in Fig.~\ref{dual naca0012 mesh}.  The Mach number distribution and the pressure distribution are shown in Fig.~\ref{dual naca0012 result}. 
The oblique shock wave can be observed at the front of the top airfoil. 
The surface pressure coefficient is also extracted and compared with the reference data \cite{jawahar2000high}, as shown in Fig.~\ref{dual naca0012 cp}. The results obtained by memory reduction CGKS agrees well with the experimental data.


\subsection{Flow around a sphere}
\noindent{\sl{(a) subsonic viscous flow around a sphere}}

\begin{figure}[htp]
	\centering	
	\includegraphics[height=0.35\textwidth]{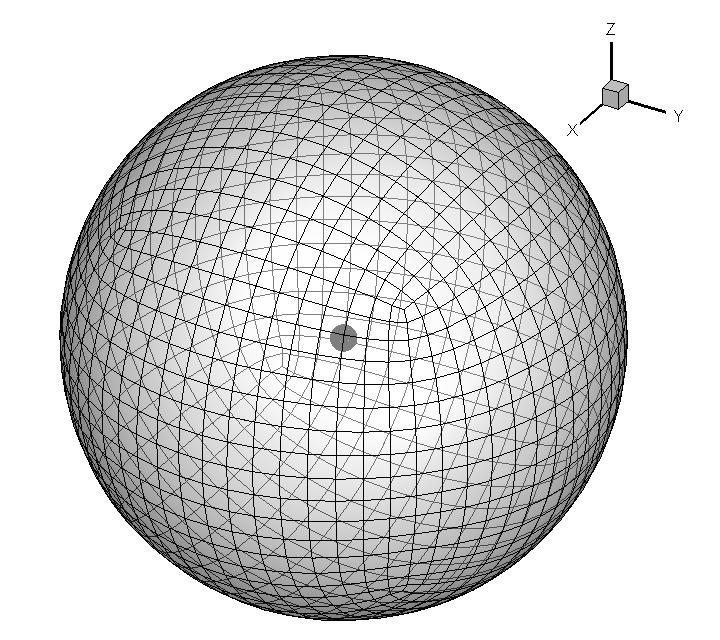}
	\includegraphics[height=0.35\textwidth]{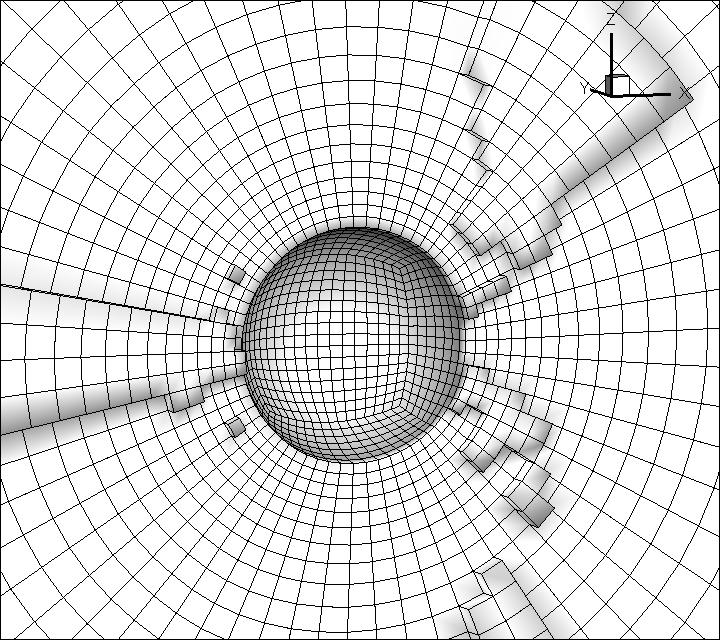}
	\caption{\label{viscous sphere mesh}
		Subsonic flow around a sphere. Mesh sample.}
\end{figure}

A subsonic flow around a sphere is simulated in this case. The Mach number is set to be 0.2535 and the Reynolds number is set to be 118.0. 
The surface of the sphere is set as non-slip and adiabatic.
The first mesh off the wall has the size $h = 4.5 \times 10^{-2}D$, and the total cell number is 50688. 
 The mesh is shown in Fig.~\ref{viscous sphere mesh}.

\begin{figure}[htp]	
	\centering	
	\includegraphics[height=0.35\textwidth]{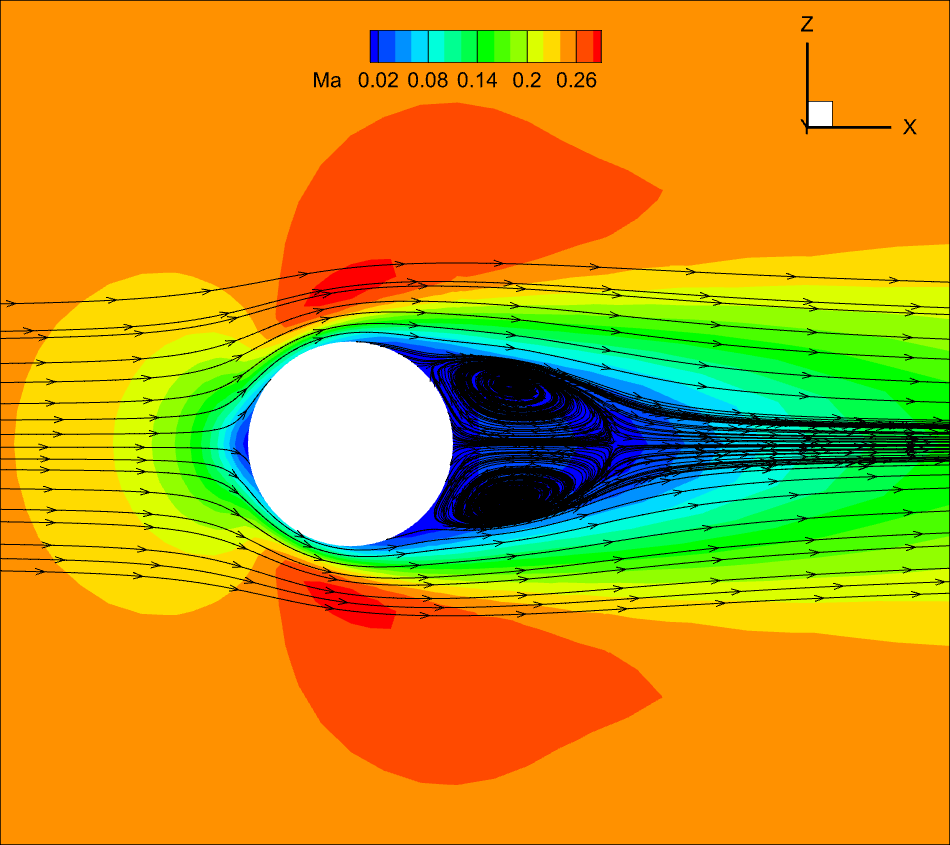}
	\includegraphics[height=0.35\textwidth]{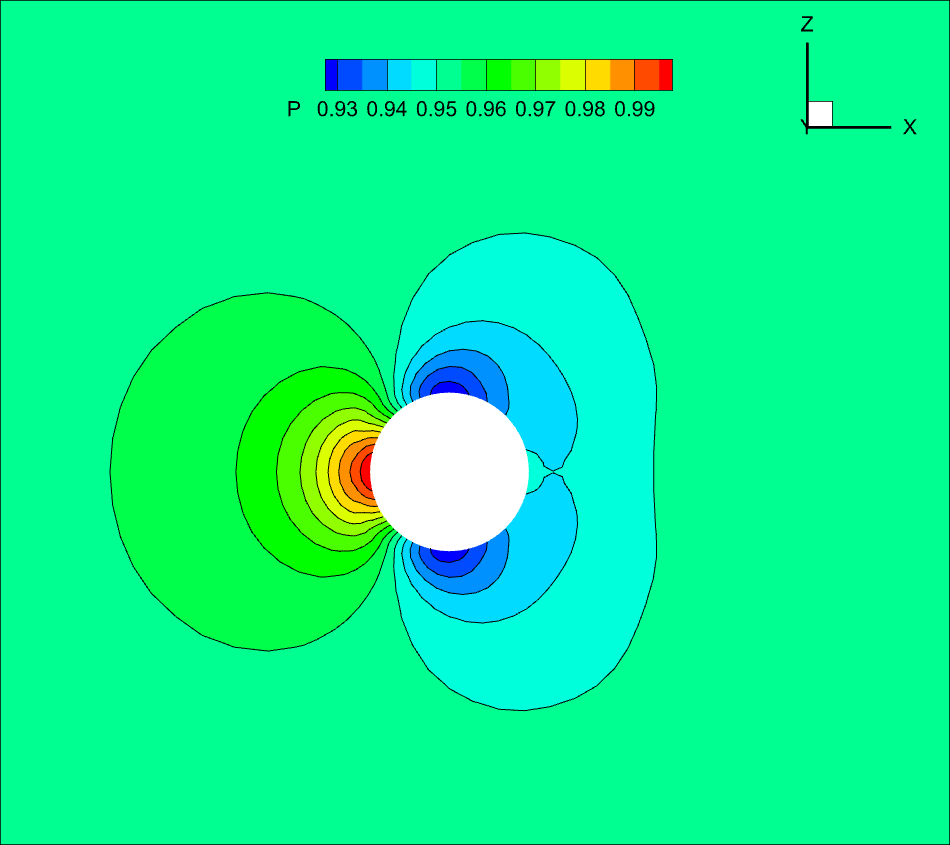}
	\caption{\label{viscous-sphere-density-contour}
		Subsonic flow around a sphere. Left: Mach contour and streamlines. Right: Pressure contour}
\end{figure}

 The Mach number contour and streamline are presented in Fig.~\ref{viscous-sphere-density-contour} to show the high resolution of the memory reduction CGKS. 
 Quantitative results are given in Table \ref{viscous subsonic sphere}, including the drag coefficient $C_D$, the separation angle $\theta$, and the closed wake length L, as defined in \cite{ji2021compact}.
 The results above show the current scheme agrees well with the experimental and numerical references.

\begin{table}[htp]
	\small
	\begin{center}
		\def\temptablewidth{1.0\textwidth}
		{\rule{\temptablewidth}{1pt}}
		\begin{tabular*}{\temptablewidth}{@{\extracolsep{\fill}}c|c|c|c|c|c}
			Scheme & Mesh number & $C_D$  & $\theta$  &L &Cl\\
			\hline
			Experiment \cite{taneda1956experimental}	&-- & 1.0  & 151 & 1.07 & -- \\ 	
			Third-order DDG \cite{cheng2017parallel} & 160,868 & 1.016 & 123.7 & 0.96 & --\\
			Fourth-order VFV \cite{wang2017thesis}  & 458,915 & 1.014 & --& -- & 2.0e-5\\
			Current & 50688 & 1.023  & 126.9 & 0.96 & 2.26e-5\\
		\end{tabular*}
		{\rule{\temptablewidth}{0.1pt}}
	\end{center}
	\vspace{-4mm} \caption{\label{viscous subsonic sphere} Quantitative comparisons among different compact schemes for the subsonic flow around a sphere.}
\end{table}

\noindent{\sl{(b) transonic viscous flow around a sphere}}

\begin{figure}[htp]	
	\centering	
	\includegraphics[height=0.35\textwidth]{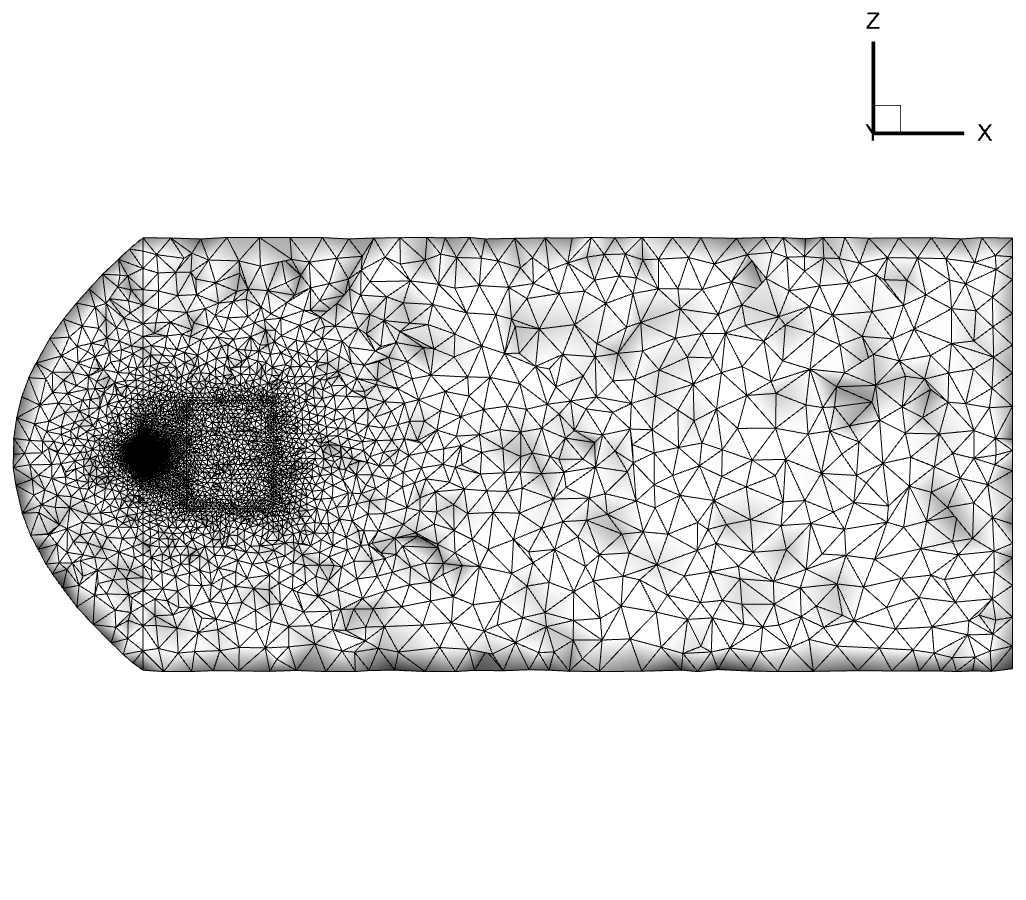}
	\includegraphics[height=0.35\textwidth]{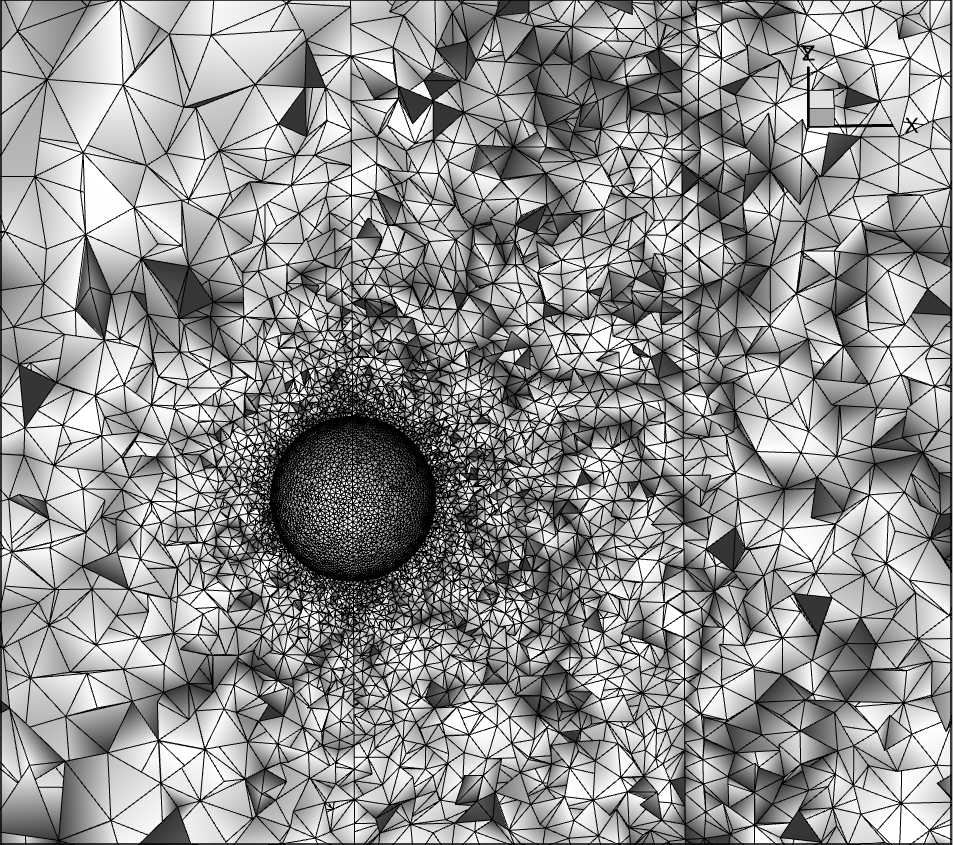}
	\caption{\label{tetrahedron-sphere-mesh}
		The mesh of Transonic flow around a sphere. Left: Global mesh. Right: Local mesh. }
\end{figure}

\begin{figure}[htp]	
	\centering	
	\includegraphics[height=0.35\textwidth]{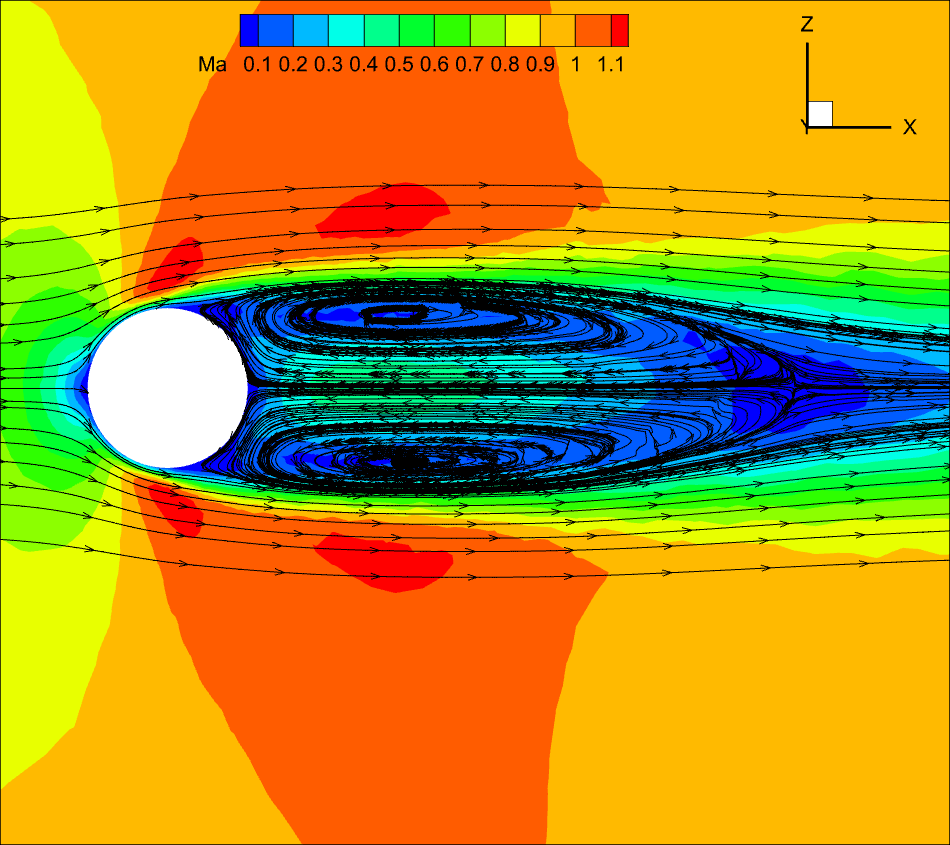}
	\includegraphics[height=0.35\textwidth]{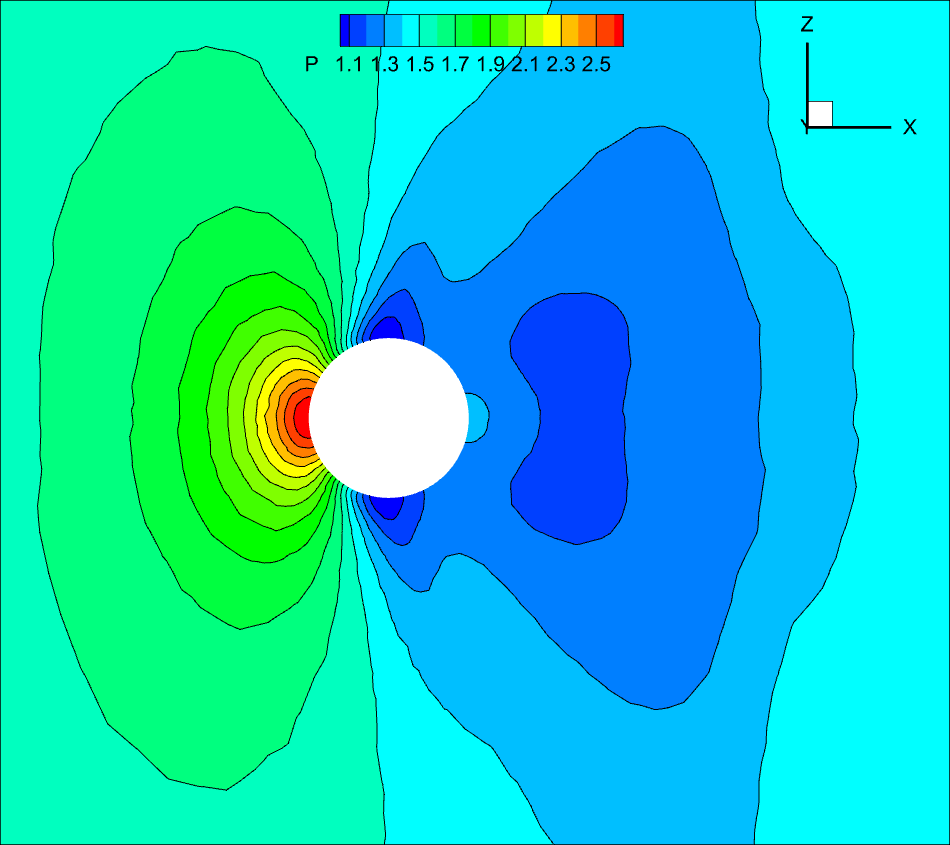}
	\caption{\label{viscous-sphere-ma0.95-contour}
		Transonic flow around a sphere. Left: Mach number contour with streamline through the sphere. Right: Pressure contour. }
\end{figure}

A transonic viscous flow around a sphere is simulated to show the performance of the memory reduction CGKS for transonic viscous flow. 
The Mach number is set to be 0.95 and the Reynolds number is set to be 300.0. 
In this case, we use the pure tetrahedron mesh with a mesh number equal to 665914, and the wake part of the sphere is refined to capture the vortex.
The mesh used in this case is shown in Fig.~\ref{tetrahedron-sphere-mesh}. 
The numerical results of the Mach number contour and streamline around a sphere are shown in Fig.~\ref{viscous-sphere-ma0.95-contour}, which indicates the high resolution of the memory reduction CGKS.
Quantitative results include the drag coefficient $C_D$,  the wake length $L$, and the separation angle $\theta$ are listed in Table \ref{transsonic-sphere}. 
The results above show the current scheme agrees well with the numerical references even using the higher order.

\begin{table}[htp]
	\small
	\begin{center}
		\def\temptablewidth{1.0\textwidth}
		{\rule{\temptablewidth}{1pt}}
		\begin{tabular*}{\temptablewidth}{@{\extracolsep{\fill}}c|c|c|c|c}
			Scheme & Mesh Number & $C_D$  & $\theta$  &L\\
			\hline
			WENO6 \cite{Nagata2016sphere} 	&909,072 & 0.968  & 111.5 & 3.48\\ 	
			Original CGKS \cite{ji2021compact}  & 515,453 & 0.950  & 112.7 & 3.30\\
   Current  & 665,914 & 0.974  & 110.0 & 3.44\\
		\end{tabular*}
		{\rule{\temptablewidth}{0.1pt}}
	\end{center}
	\vspace{-4mm} \caption{\label{transsonic-sphere} Quantitative comparisons between the current scheme and the reference solution for the transonic flow around a sphere.}
\end{table}

\noindent{\sl{(c) supersonic viscous flow around a sphere}}

\begin{figure}[htp]	
	\centering	
	\includegraphics[height=0.35\textwidth]{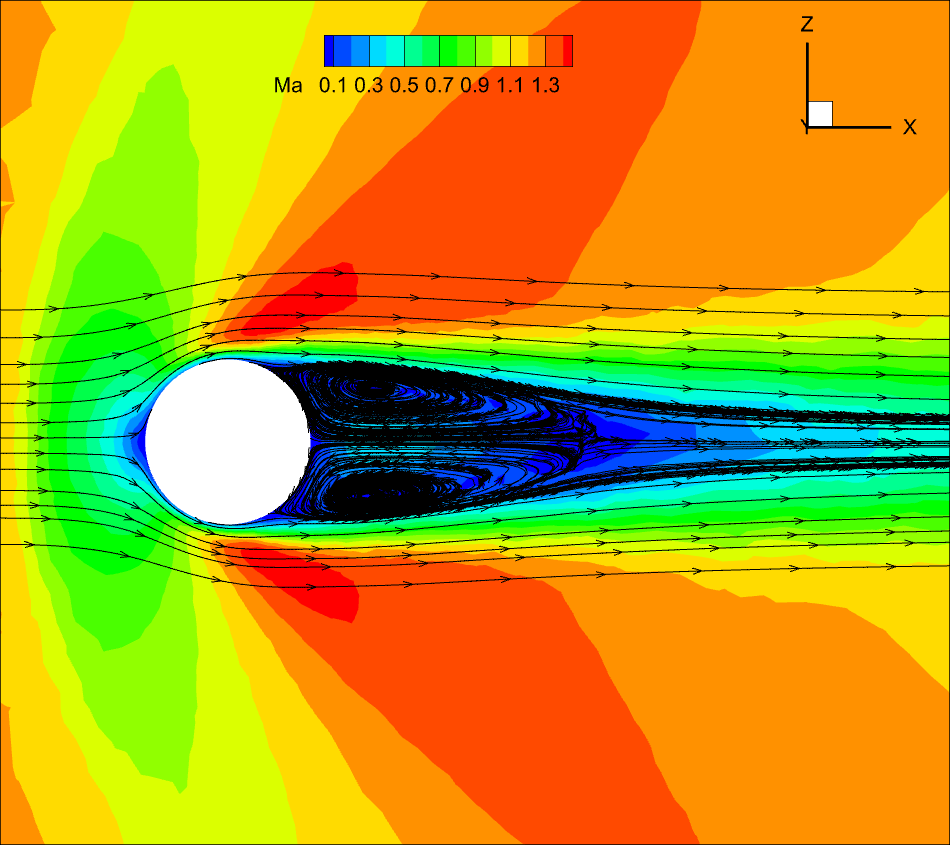}
	\includegraphics[height=0.35\textwidth]{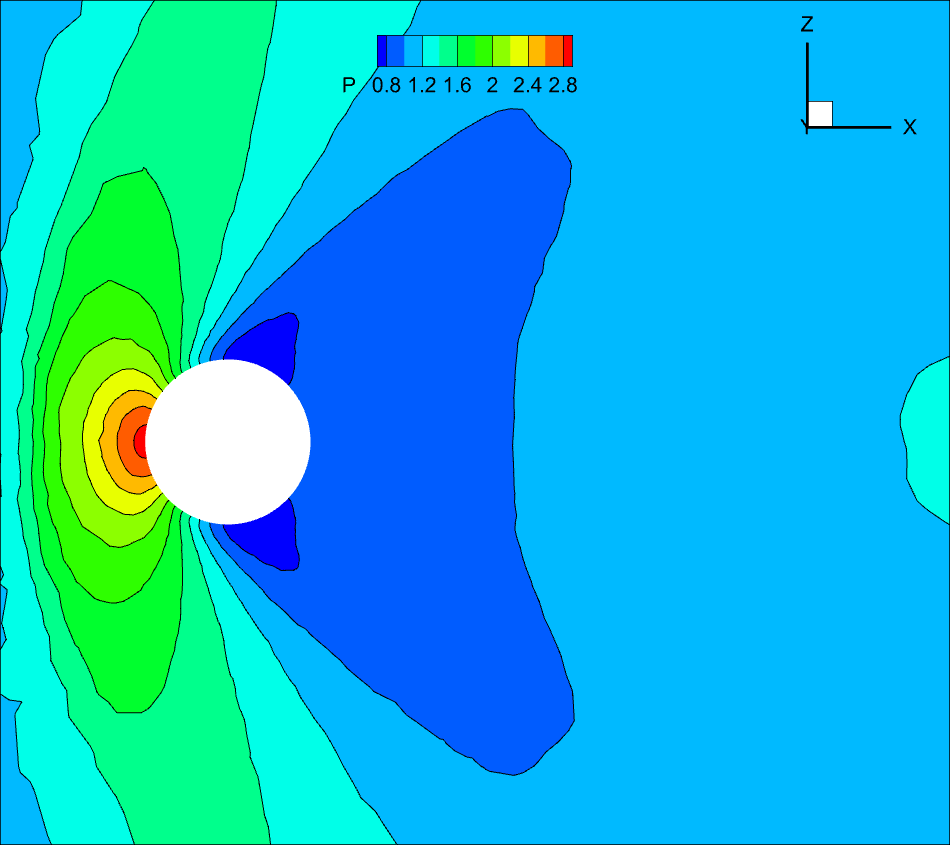}
	\caption{\label{viscous-sphere-ma1x2-contour}
		Supersonic flow around a sphere. Left: Mach number contour with streamline through the sphere. Right: Pressure contour. }
\end{figure}

To verify that the memory reduction CGKS can also have good performance in the supersonic flow region, a supersonic flow around a sphere is simulated. 
The Mach number is set to be 1.2 and the Reynolds number is set to be 300.  
The mesh used in this case is the same as the transonic case. The upstream length is 5 and the downstream length is 40. 
The first layer mesh at the wall has a thickness $2.3 \times 10^{-2}D$. 
The result of Mach number with streamline around the sphere is also shown in Fig.~\ref{viscous-sphere-ma1x2-contour}, which indicates the high resolution of the memory reduction CGKS. 
Quantitative results are listed in Table \ref{supersonic-sphere}, which agrees well with those given by Ref.\cite{Nagata2016sphere}.

\begin{table}[htp]
	\small
	\begin{center}
		\def\temptablewidth{1.0\textwidth}
		{\rule{\temptablewidth}{1pt}}
		\begin{tabular*}{\temptablewidth}{@{\extracolsep{\fill}}c|c|c|c|c|c}
			Scheme & Mesh Number & Cd  & $\theta$  &L & Shock stand-off\\
			\hline
			WENO6 \cite{Nagata2016sphere} 	&909,072 & 1.281  & 126.9 & 1.61 & 0.69 \\ 	
			Original CGKS \cite{JI2024112590}  & 665,914 & 1.274  & 126.3 & 1.64 & 0.72 \\
                Current  & 665,914 & 1.303  & 126.3 & 1.60 & 0.72 \\
		\end{tabular*}
		{\rule{\temptablewidth}{0.1pt}}
	\end{center}
	\vspace{-4mm} \caption{\label{supersonic-sphere} Quantitative comparisons between the current scheme and the reference solution for the supersonic flow around a sphere.}
\end{table}

\subsection{M6-wing}
\begin{figure*}[htp]	
	\centering	
	\includegraphics[height=0.35\textwidth]{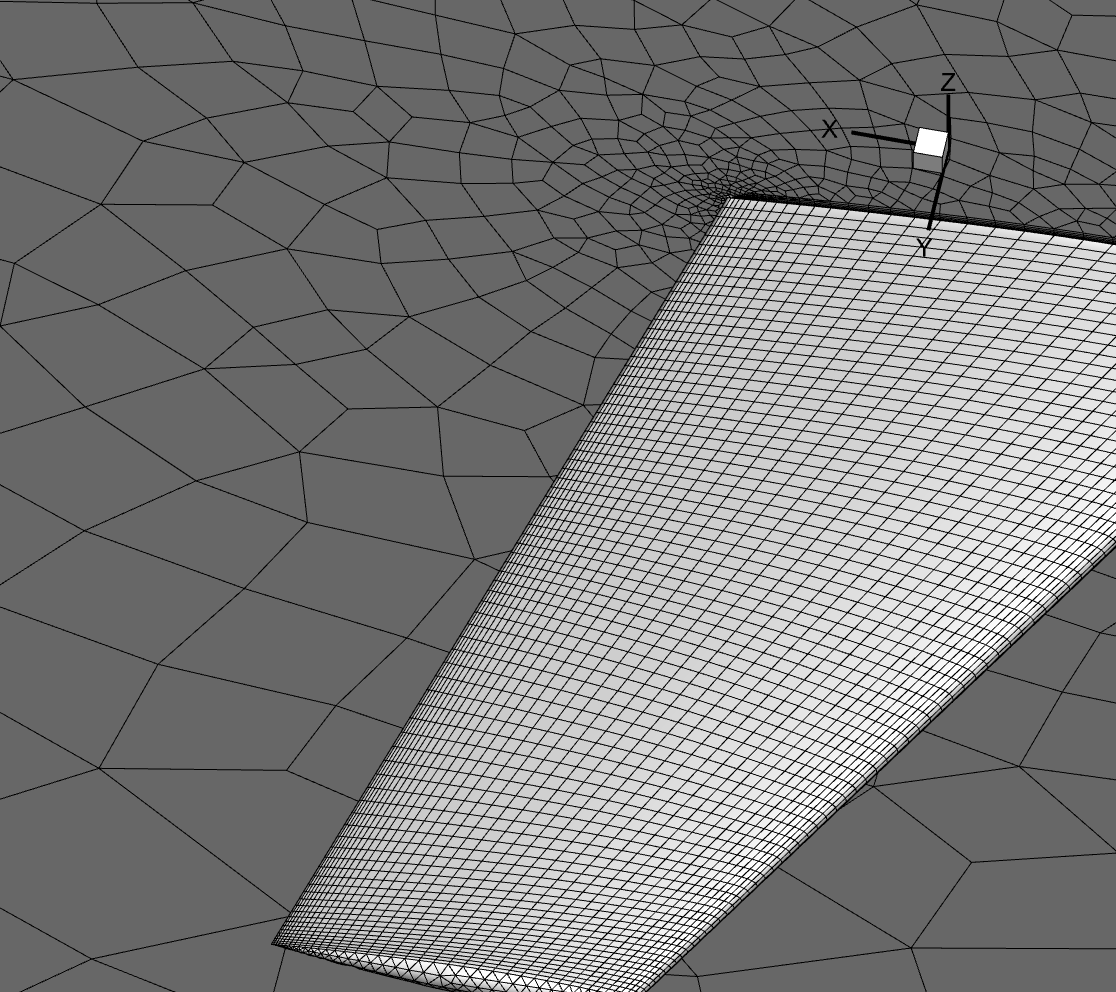}
	\includegraphics[height=0.35\textwidth]{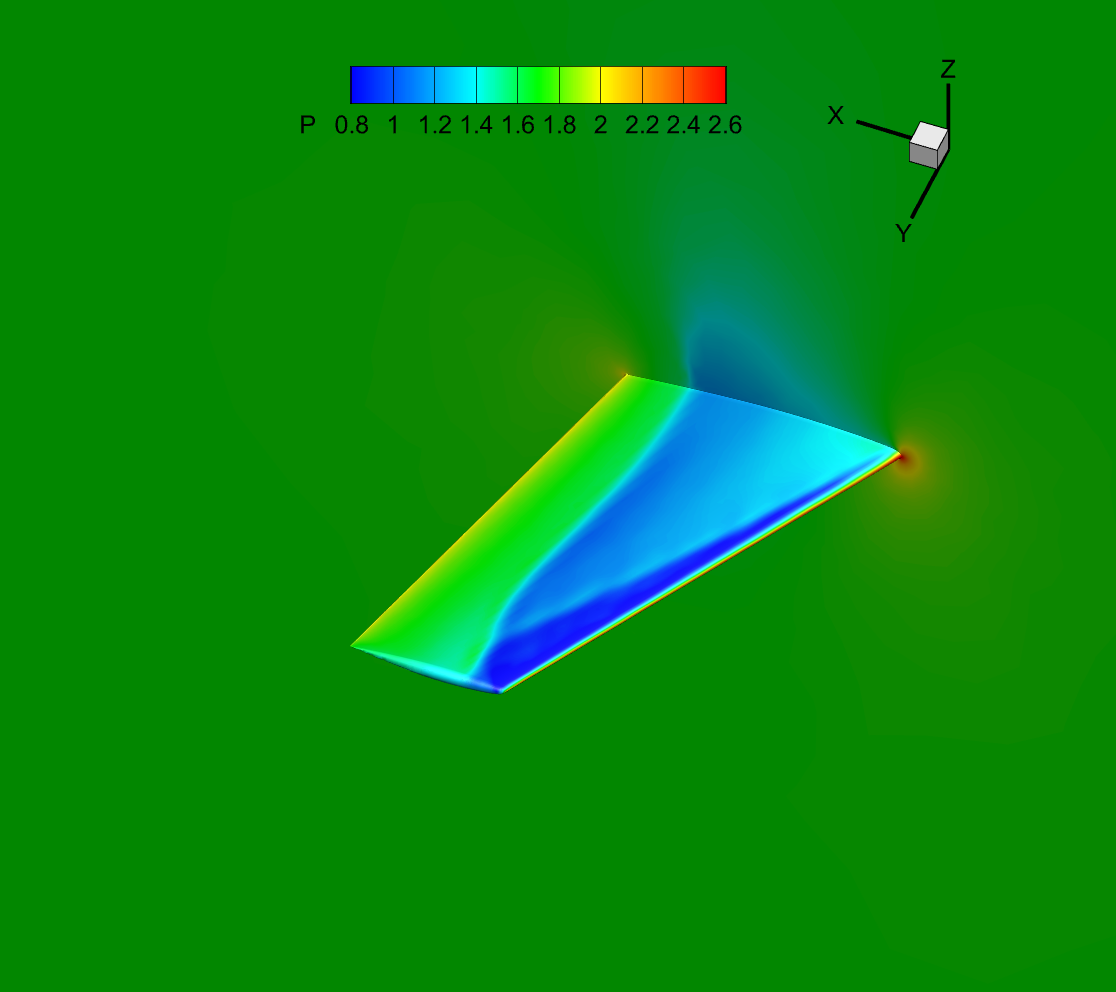}
	\caption{\label{M6 Contour}
		Transonic flow around a M6 wing. Left: Mesh. Right: Pressure contour.}
\end{figure*}

\begin{figure}[htp]	
	\centering	
	\includegraphics[height=0.35\textwidth]{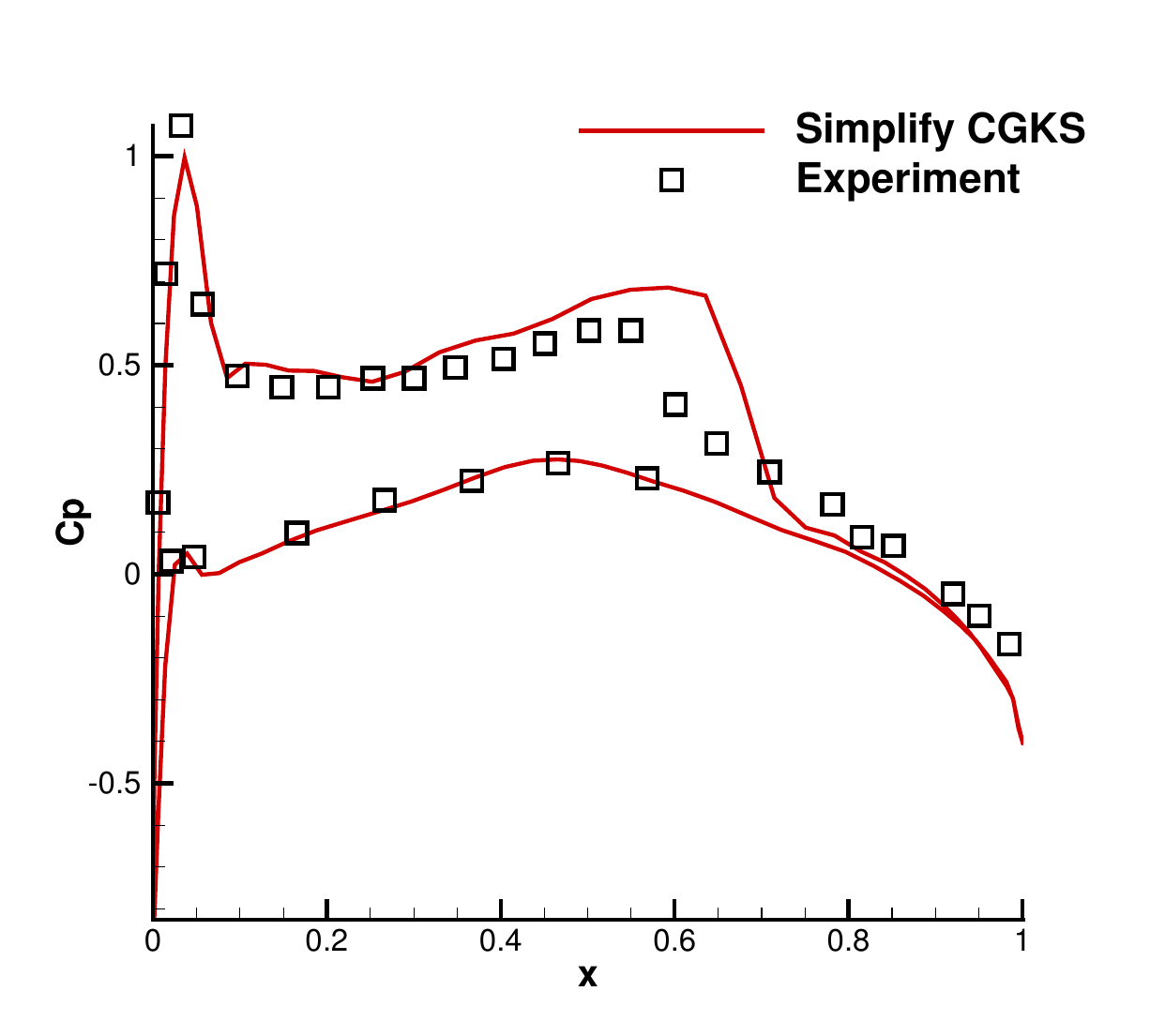}
	\includegraphics[height=0.35\textwidth]{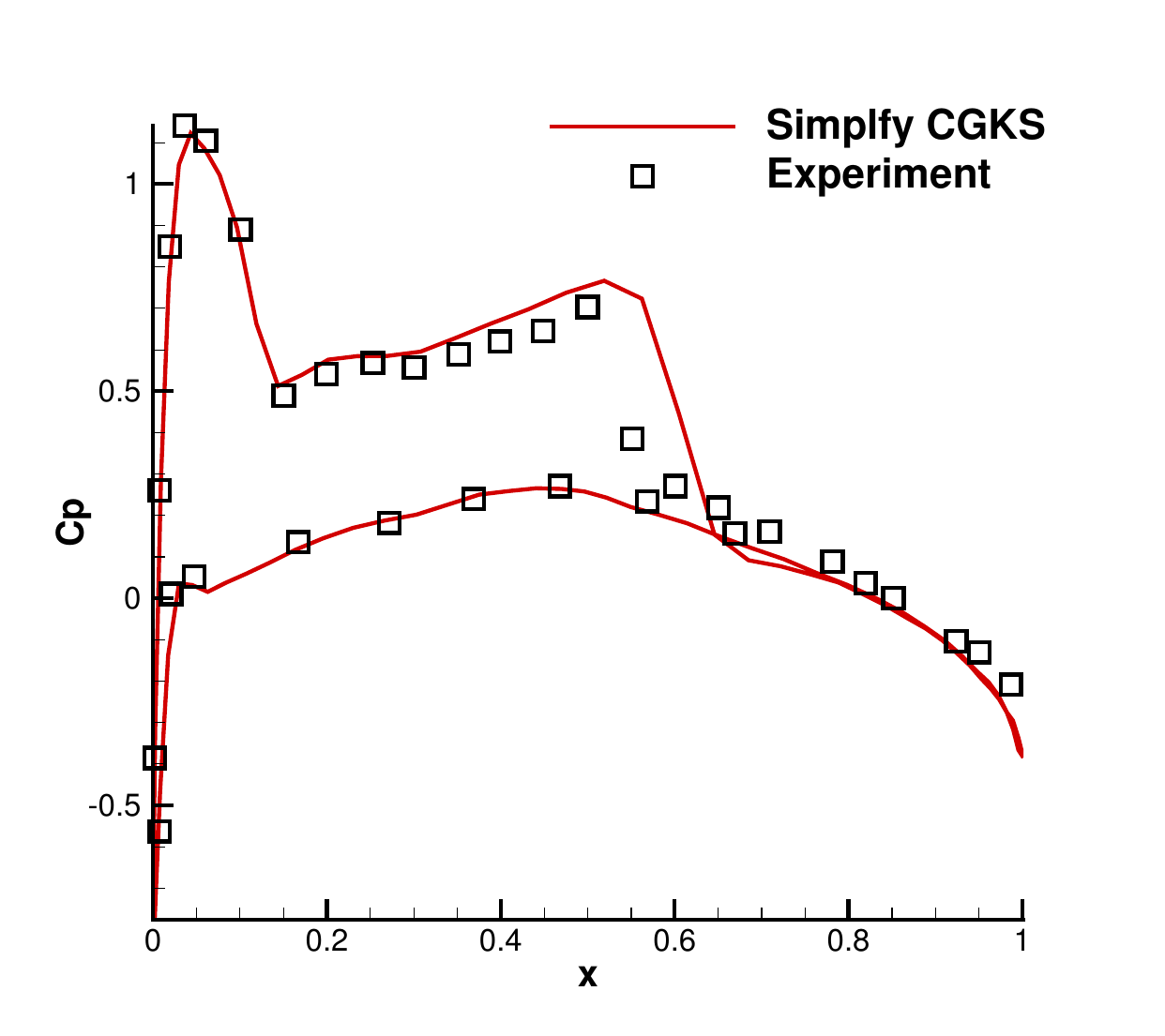}
	\includegraphics[height=0.35\textwidth]{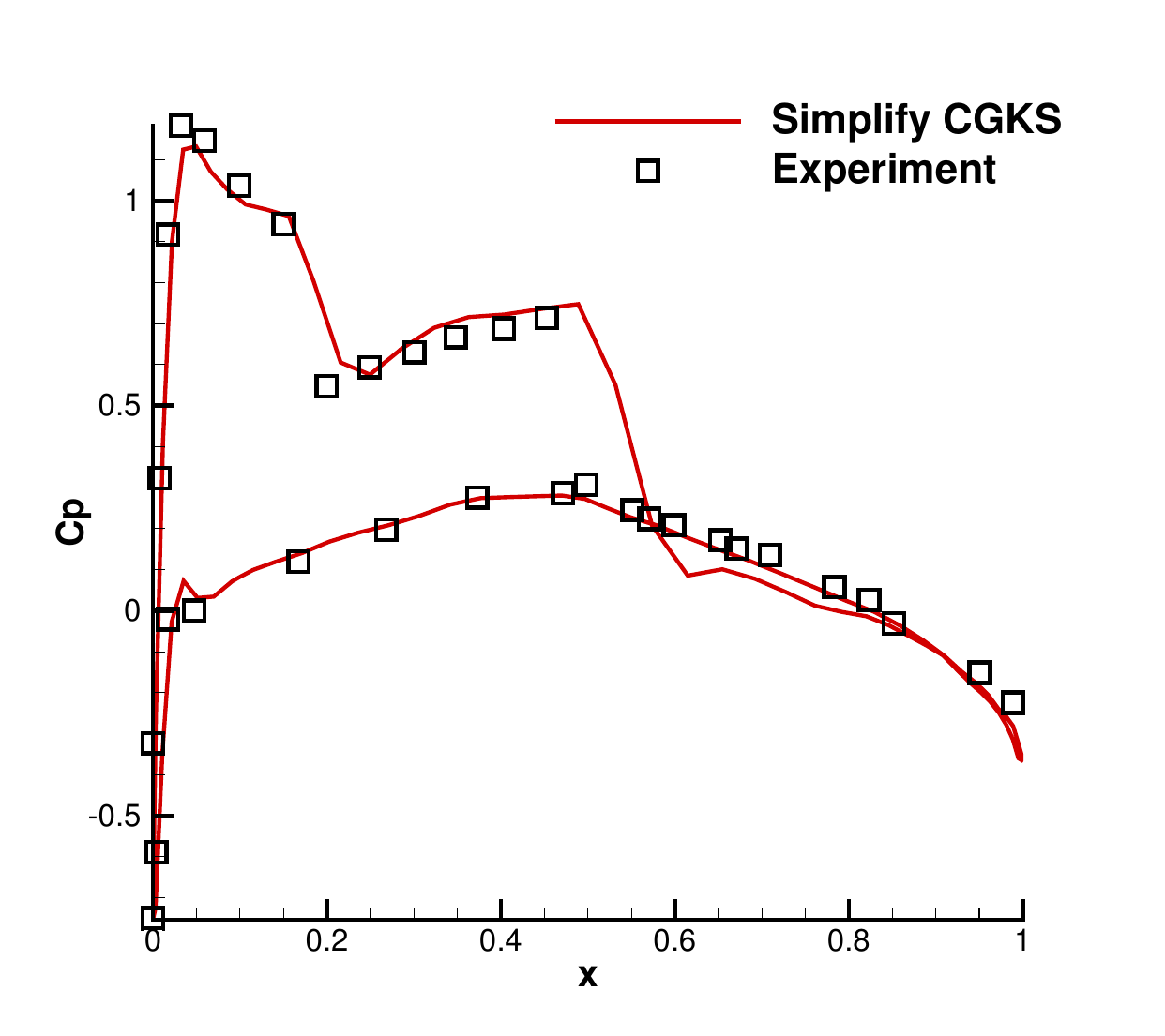}
	\includegraphics[height=0.35\textwidth]{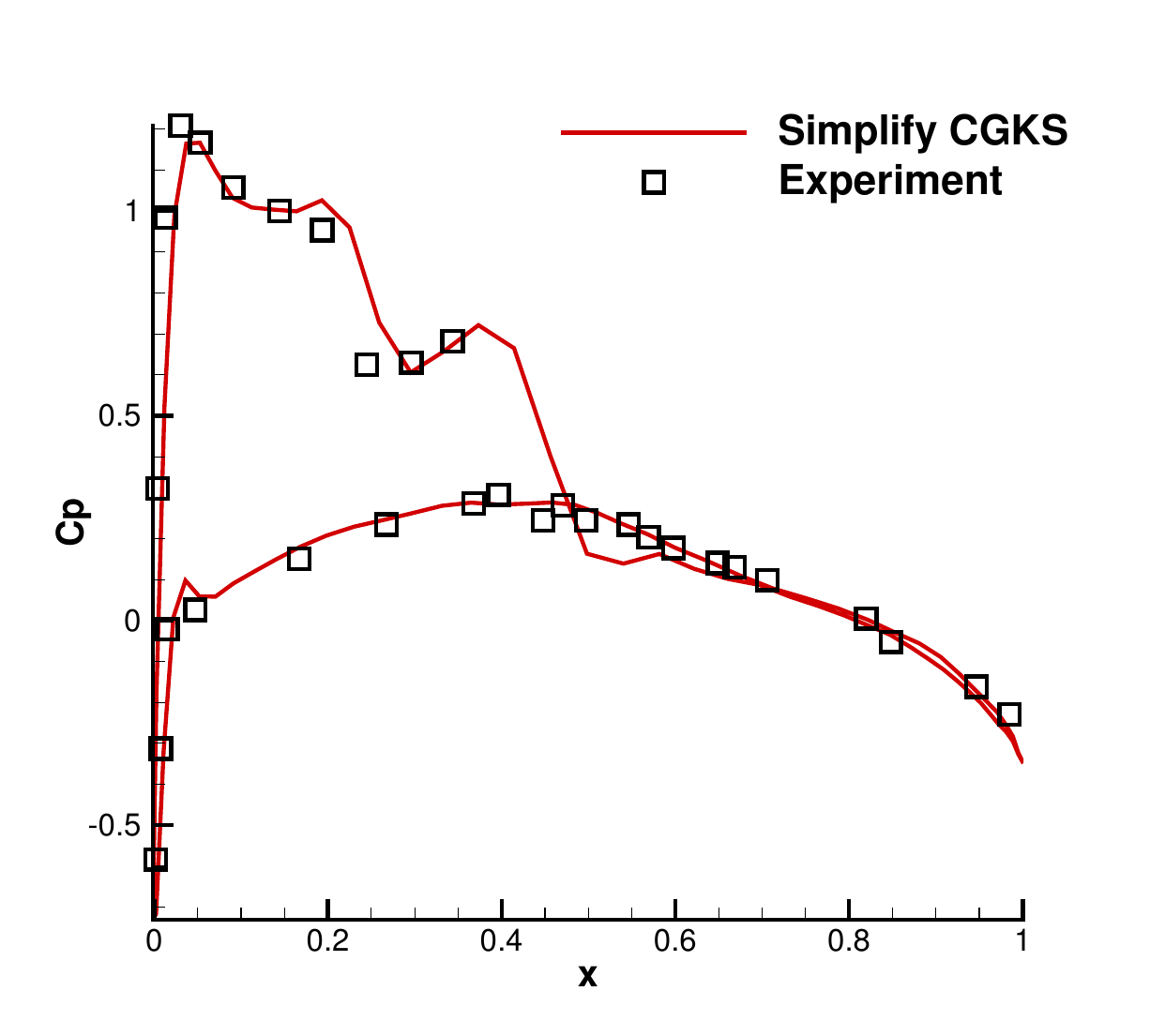}
	\includegraphics[height=0.35\textwidth]{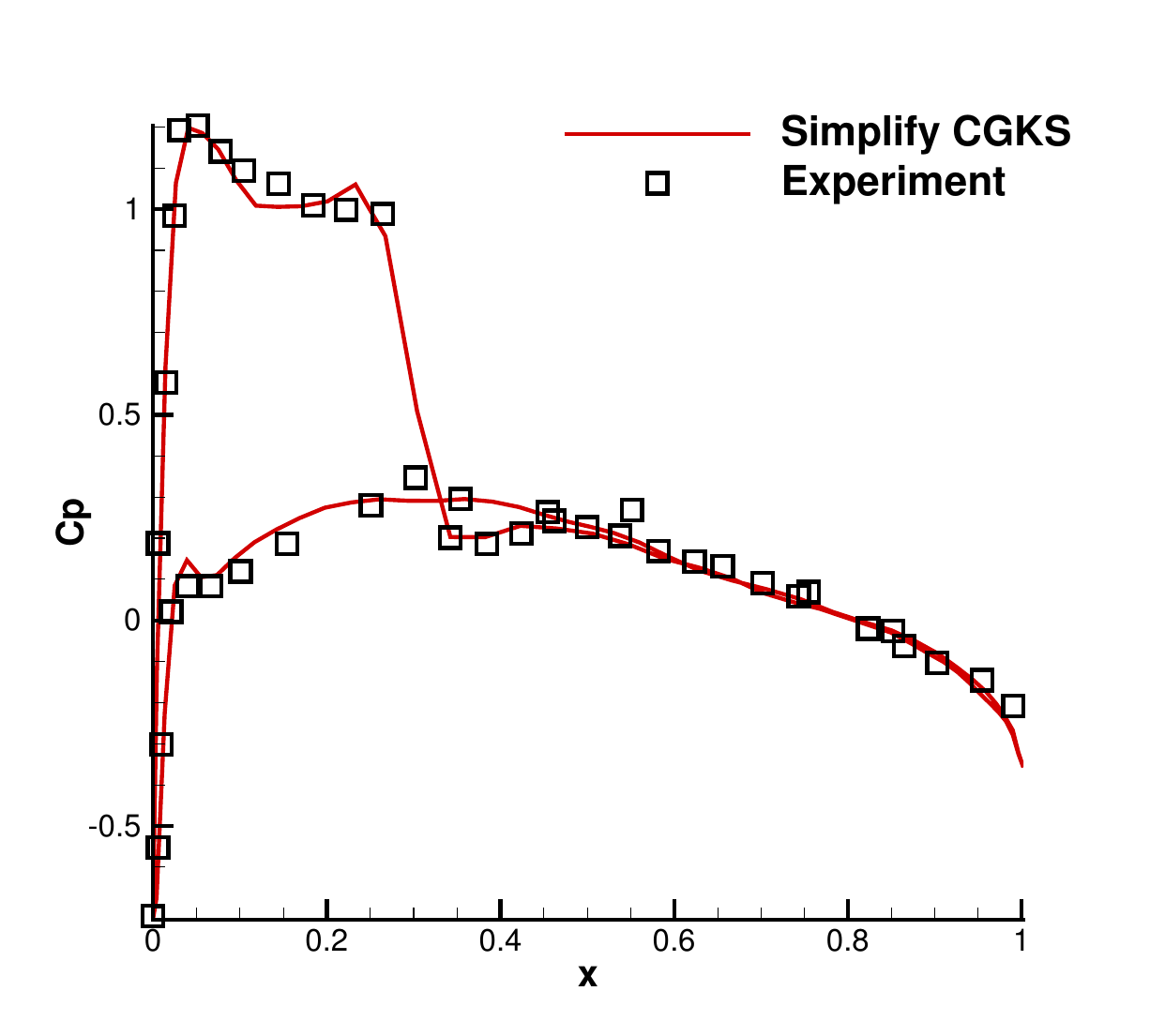}
	\includegraphics[height=0.35\textwidth]{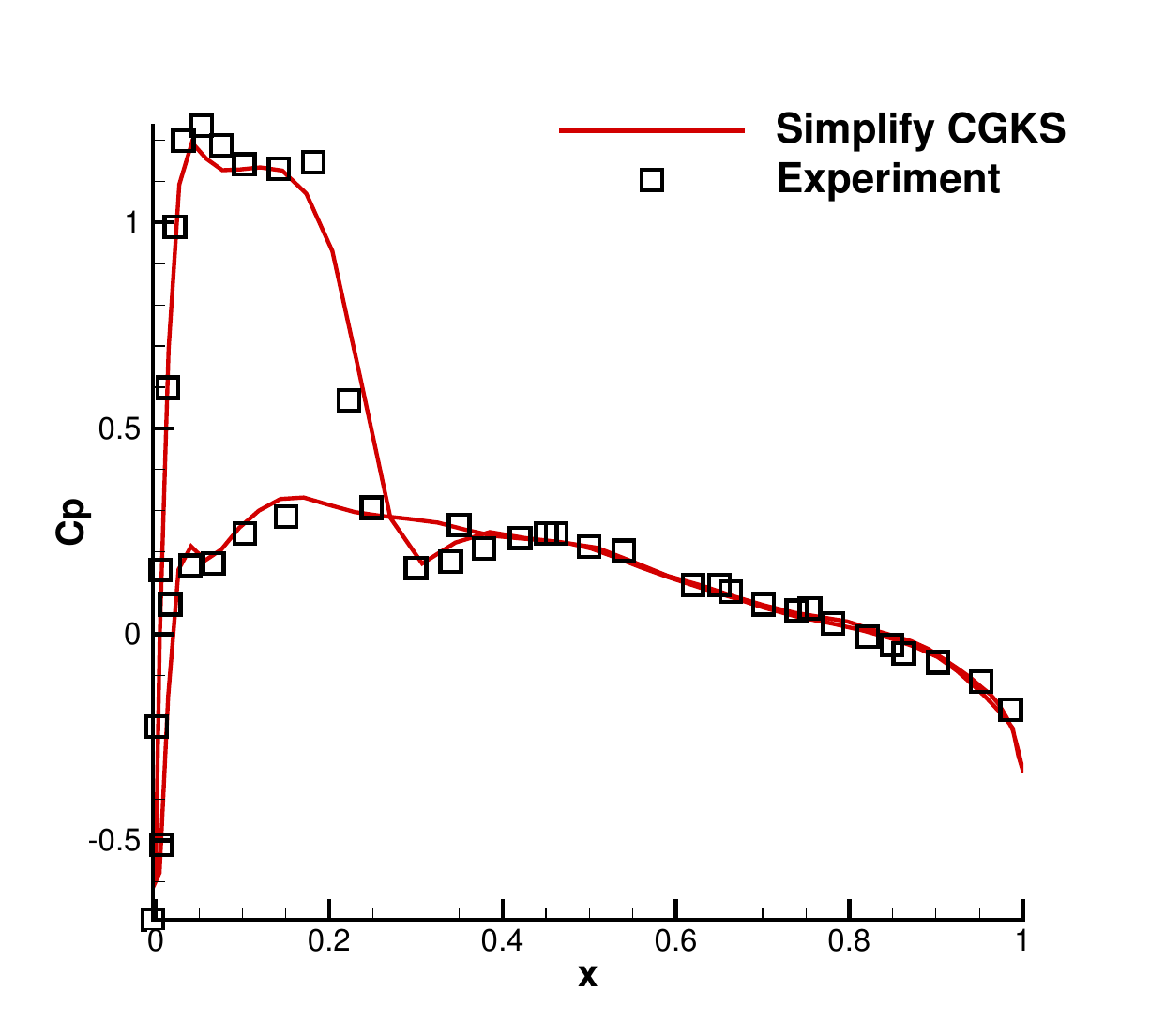}
	\caption{\label{M6 Contour-Cp}
		Transonic flow around a M6 wing. Cp distribution.}
\end{figure}

Transonic flow around an ONERA M6 \cite{eisfeld2006onera} wing is a widely used engineering case to verify the acceleration techniques used in CFD \cite{yang2023implicit}. 
The flow structure of it is complicated due to the interaction of shock and wall boundary. 
Moreover, three-dimensional mixed unstructured mesh is also a challenge to high-order schemes. 
Thus, it is an appropriate test case to verify the accuracy and robustness of the memory reduction CGKS. 
The far-field Mach number is set to be 0.8395 and the angle of attack is set to be 3.06$^{\circ}$. 
The adiabatic slip wall boundary is used on the surface of the ONERA M6 wing and the subsonic inflow boundary is set according to the local Riemann invariants.  
A hybrid unstructured mesh with a near-wall size $h\approx 2e^{-3}$ is used in the computation, as shown in Fig.~\ref{M6 Contour}. 
The pressure distribution on the wall surface is shown in Fig.~\ref{M6 Contour}. 
The pressure contour in Fig.~\ref{M6 Contour} indicates that the memory reduction CGKS has captured the shock accurately. 
The quantitative comparisons on the pressure distributions at the semi-span locations Y /B = 0.20, 0.44, 0.65, 0.80, 0.90, and 0.95 of the wing are given in Fig.~\ref{M6 Contour-Cp}. The numerical results quantitatively agree well with the experimental data.

\subsection{Supersonic flow around a rocket fairing}
\begin{figure}[htp]	
	\centering	
        \includegraphics[height=0.35\textwidth]{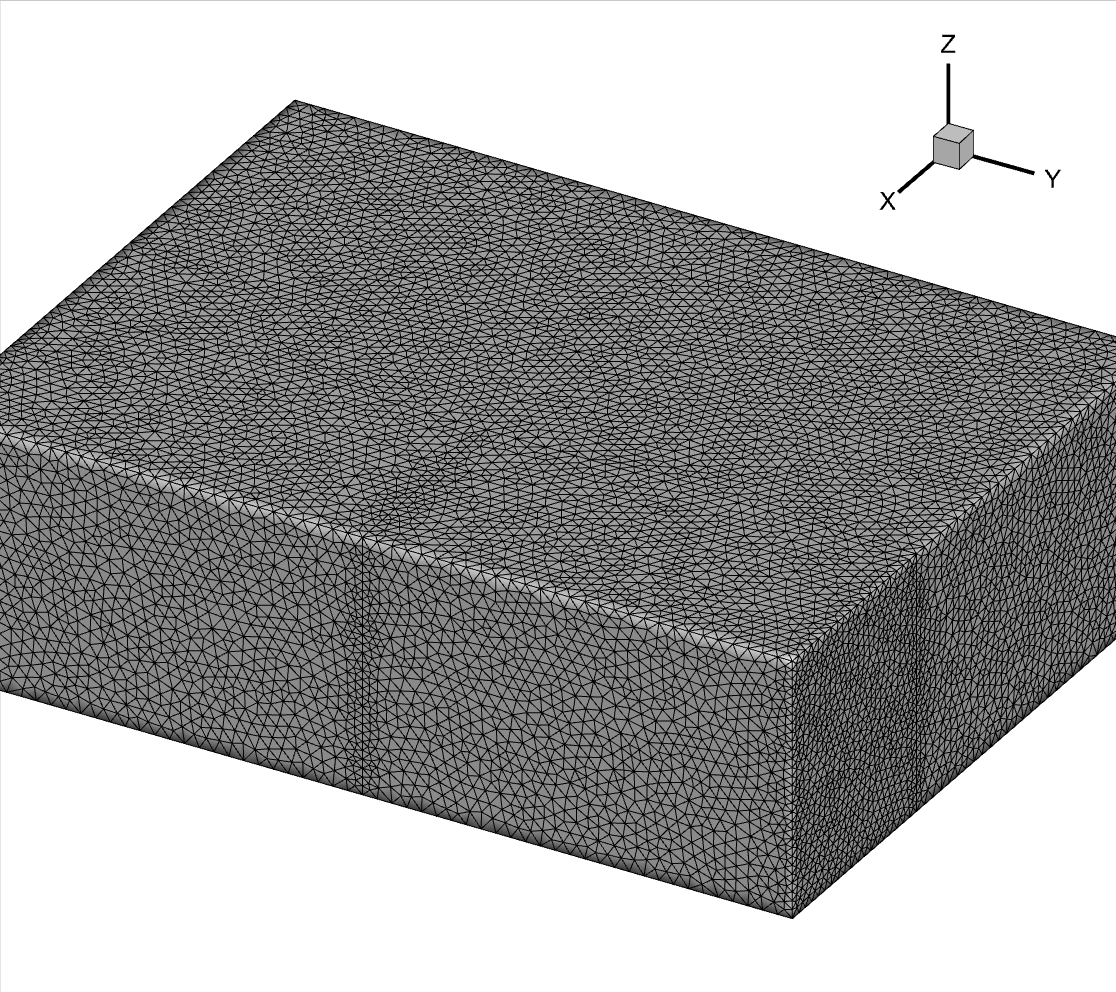}
        \includegraphics[height=0.35\textwidth]{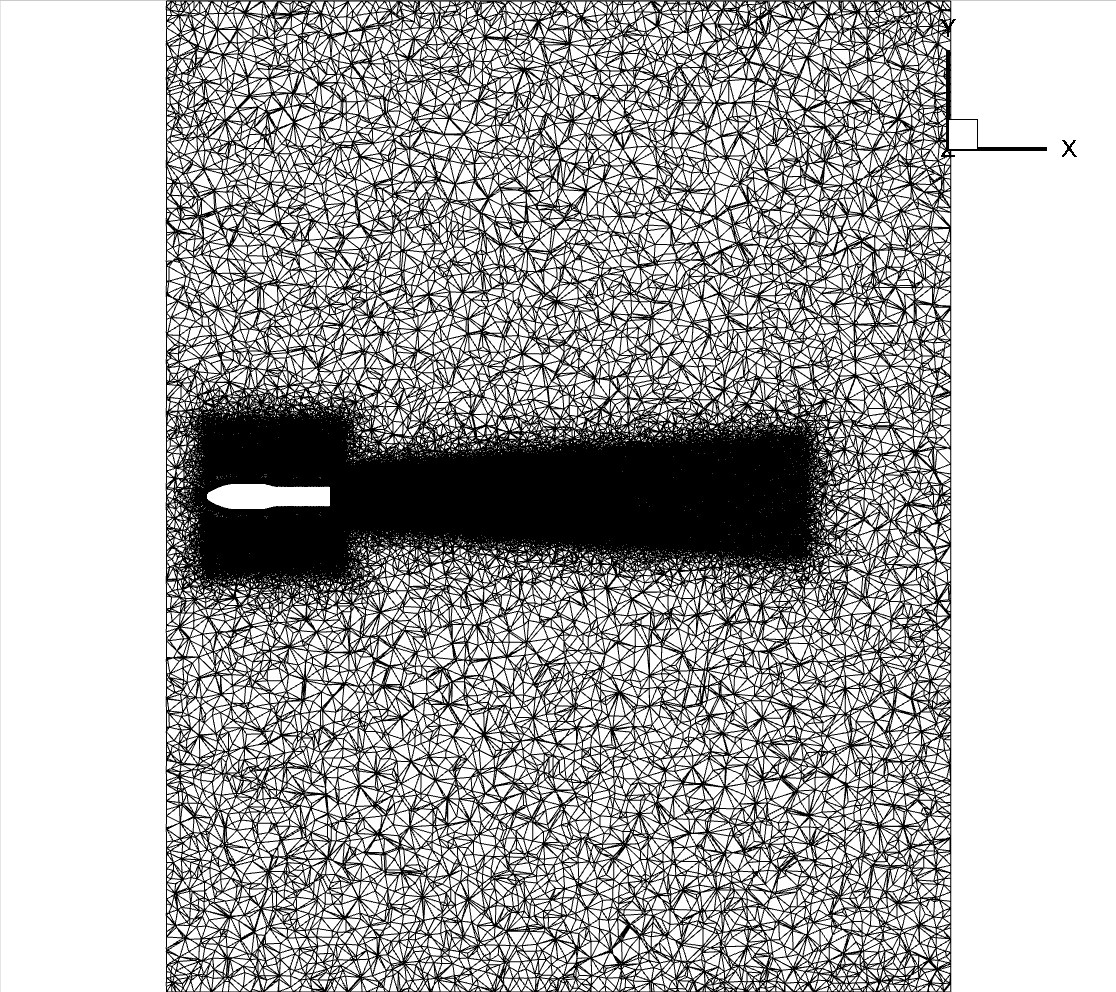}
	\caption{\label{rocket-mesh}
		The mesh of the supersonic flow around a rocket fairing. Left: Far-field mesh. Right: Slice of Z-plane.}
\end{figure}

In this case, a supersonic flow around a rocket fairing is simulated which can show the efficiency and robustness of the memory reduction CGKS.
The length of the rocket is 5 meters. 
To capture the unsteady vortex structure at the tail of the rocket and the shock at the head, the corresponding parts of the mesh are refined.
The element number of this mixed unstructured mesh is 11,948,652 and $y^+$ is set to be 10.
The mesh is shown in Fig.~\ref{rocket-mesh}

\begin{figure}[htp]	
	\centering	
        \includegraphics[height=0.4\textwidth]{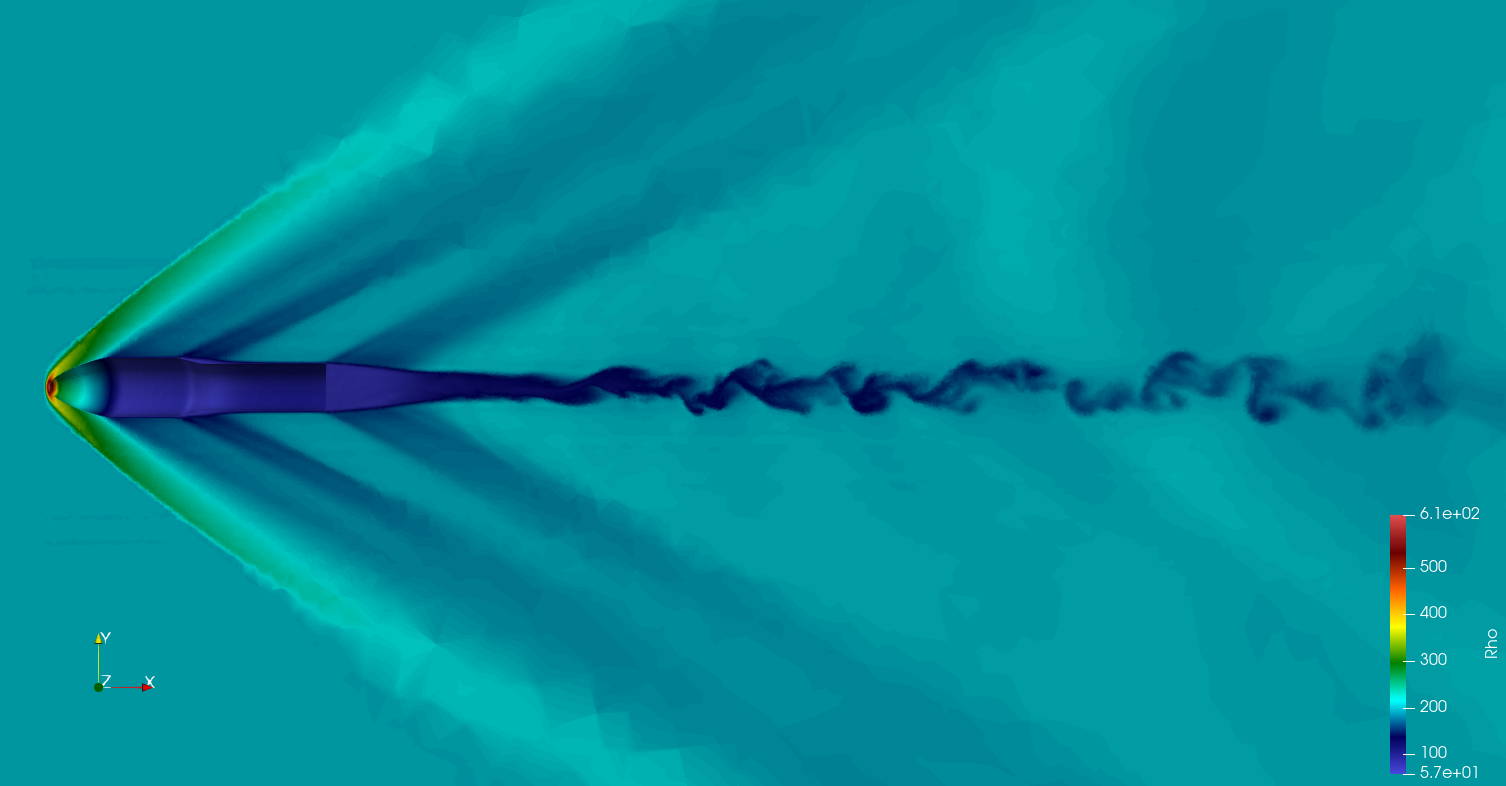}
        \includegraphics[height=0.4\textwidth]{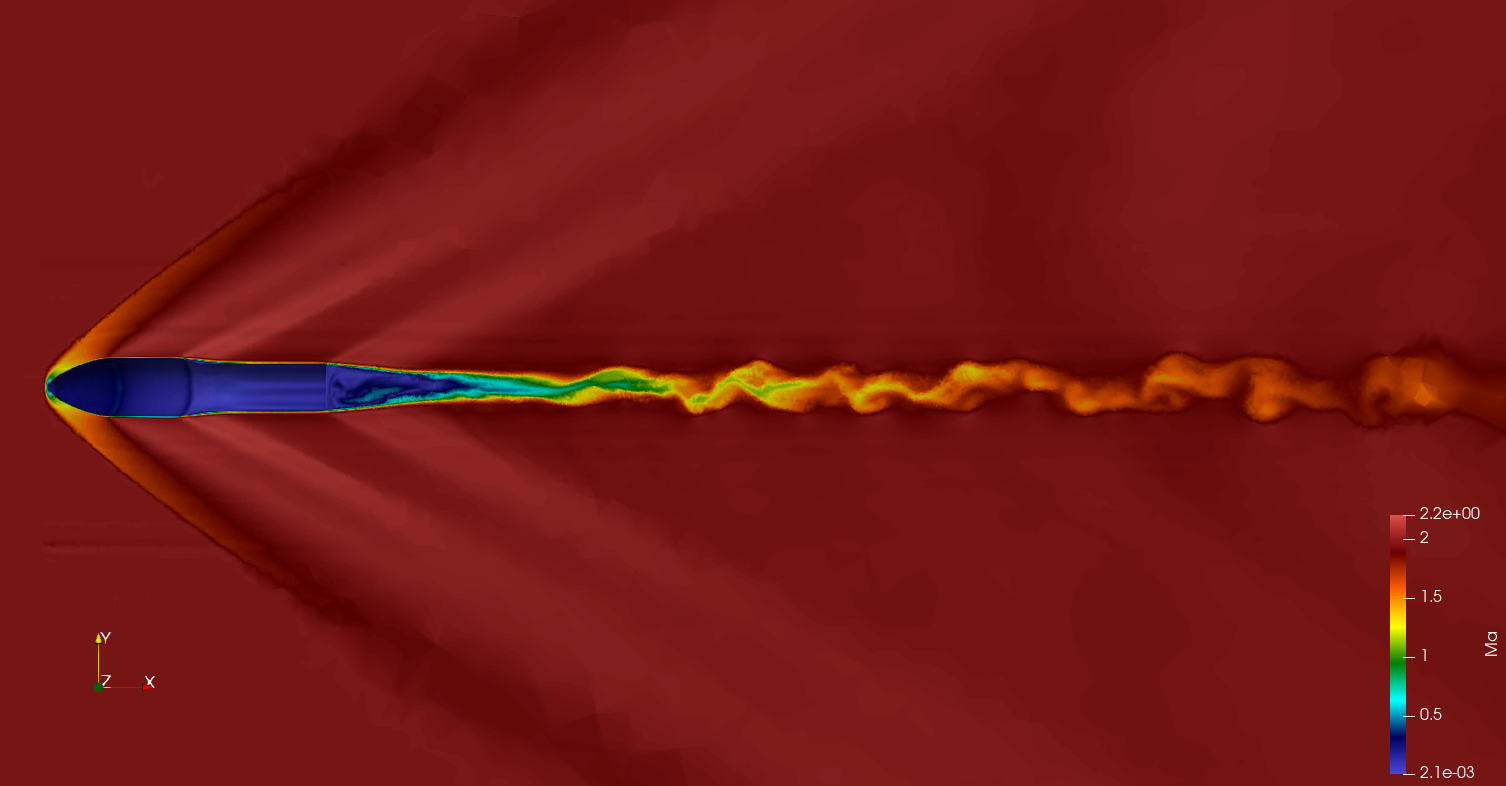}
         \includegraphics[height=0.4\textwidth]{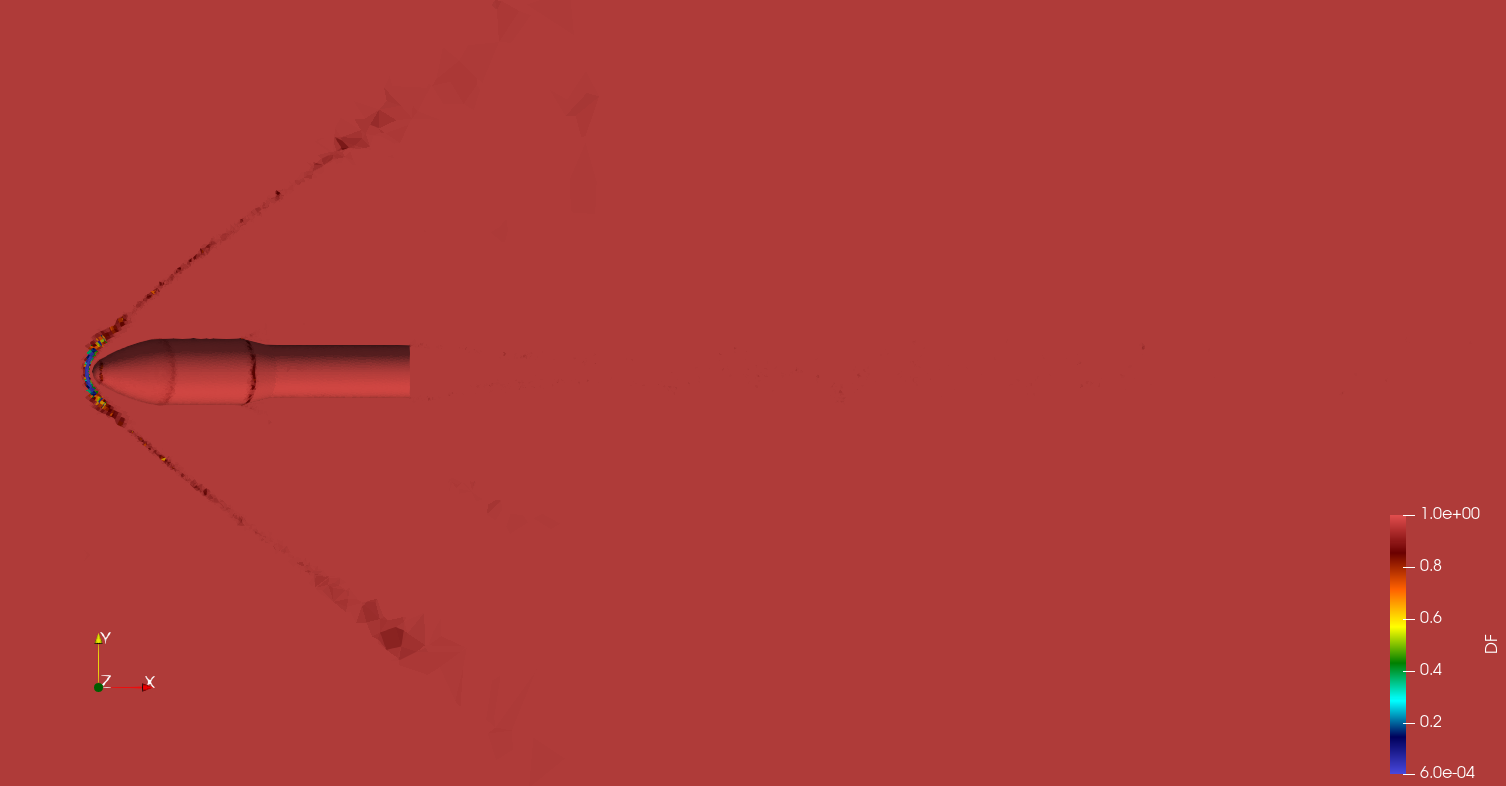}
	\caption{\label{rocket-contour}
		The contour of the supersonic flow around a rocket fairing. Top: Density contour. Middle: Mach number contour. Down: Discontinuity feedback distribution.}
\end{figure}

The incoming Mach number is set to be 2 and the Reynolds number is set to be 100,000.
The Mach number contour and density contour are shown in Fig.~\ref{rocket-contour}.
From the contour above, the current scheme can not only capture the shock sharply but also resolve the vortex at the tail of the rocket, which is a big challenge for high-order schemes. 
That means current memory reduction CGKS can deal well with both strong discontinuity and vortex in such a large-scale simulation.

\subsection{Hypersonic flow around a scramjet-powered lifting-body configuration}
\begin{figure}[htp]	
	\centering	
	\includegraphics[height=0.35\textwidth]{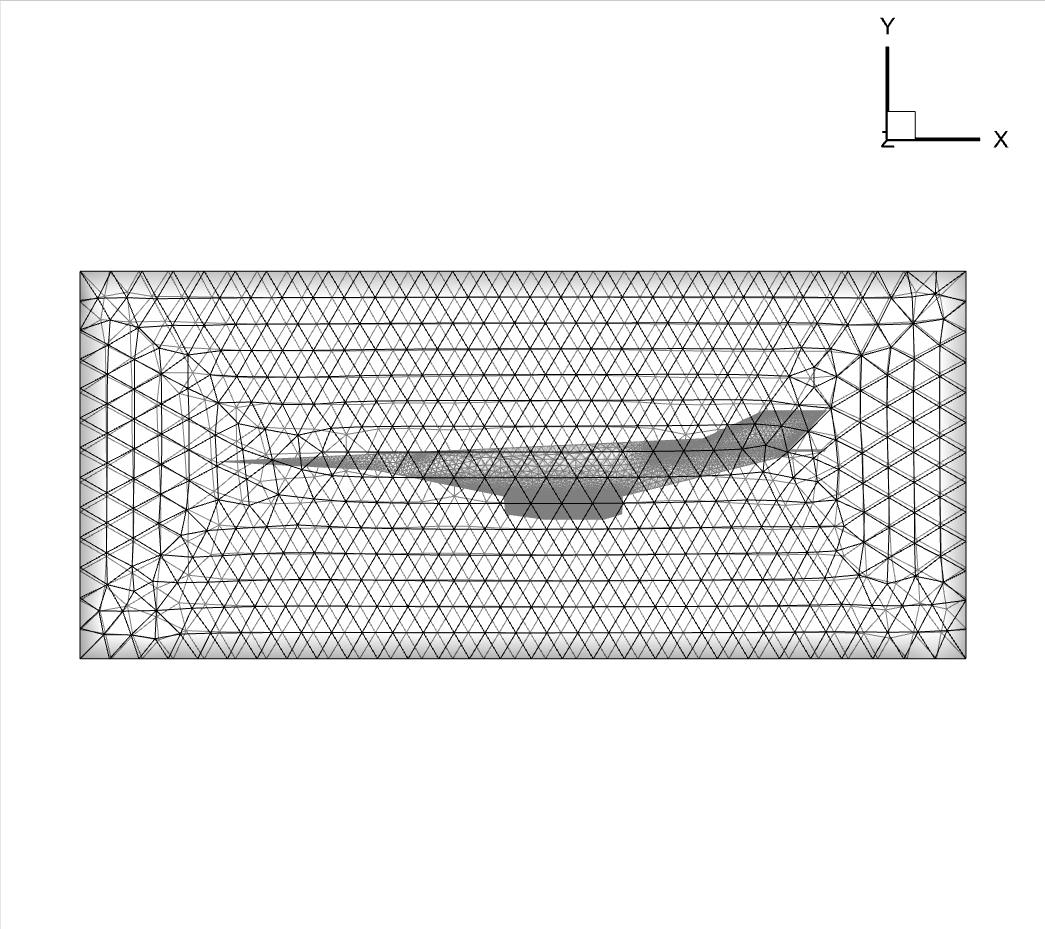}
	\includegraphics[height=0.35\textwidth]{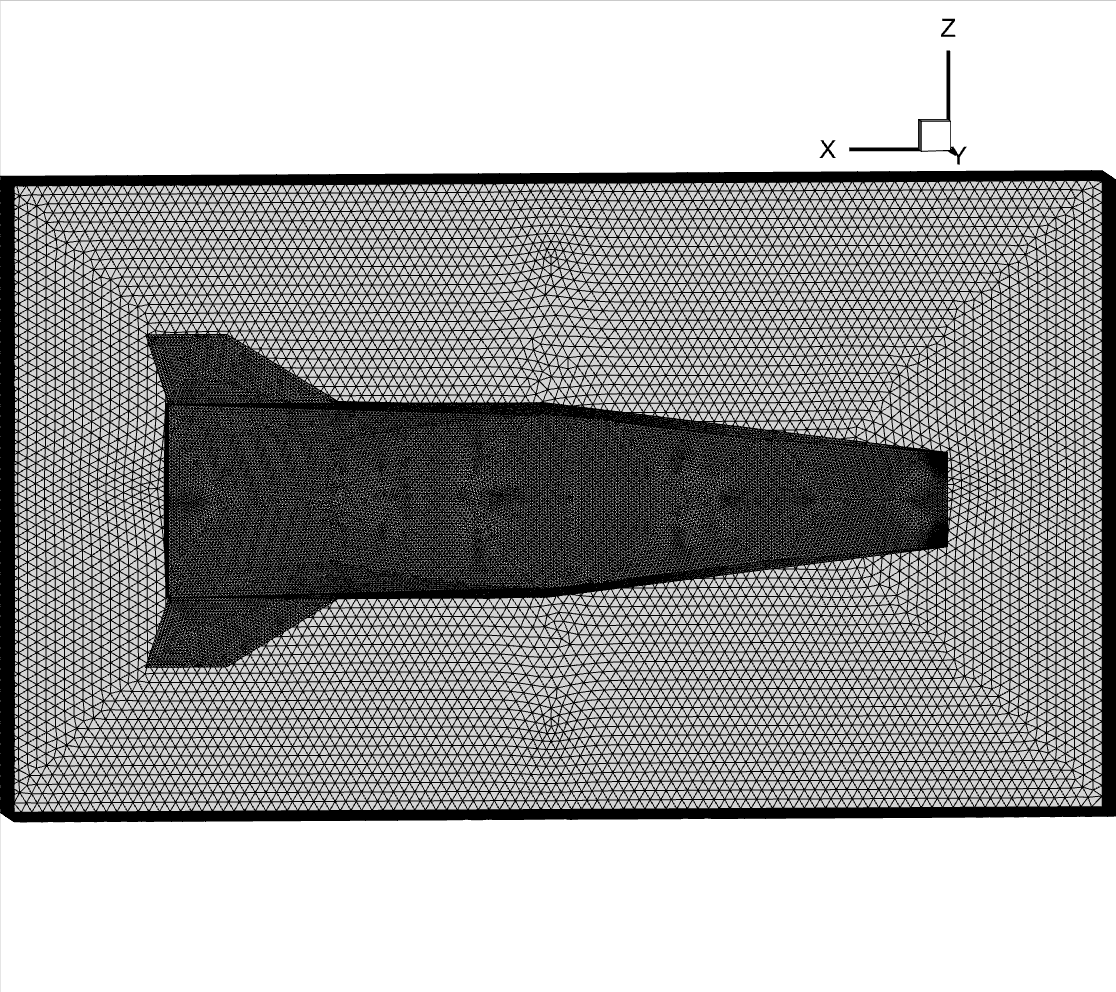}
	\caption{\label{x43-mesh}
		     Mesh of X-43A-like aircraft. Left: Global mesh. Right: Local mesh.}
\end{figure}

In this section, to highlight the advantages of high computational efficiency and strong robustness of the memory reduction CGKS, a hypersonic flow around an X-43A-like aircraft is simulated.
The X-43A aircraft is a scramjet-powered lifting-body vehicle with over 3.7m in length.
The nose of the X-43A aircraft is extremely sharp, and the transition in the middle section is not smooth. Additionally, the bottom of its fuselage features an air intake. In addition, the typical speed of X43-A is Mach number of 7.
These features make it difficult to simulate the X43-A aircraft using higher-order schemes.
In this test case, a mixed unstructured mesh is used and the total mesh number is about $5.6$ million. 
To capture the shock, the near-wall part of the mesh is refined. 
The mesh is shown in Fig.~\ref{x43-mesh}.

\begin{figure}[htp]	
	\centering	
	\includegraphics[height=0.40\textwidth]{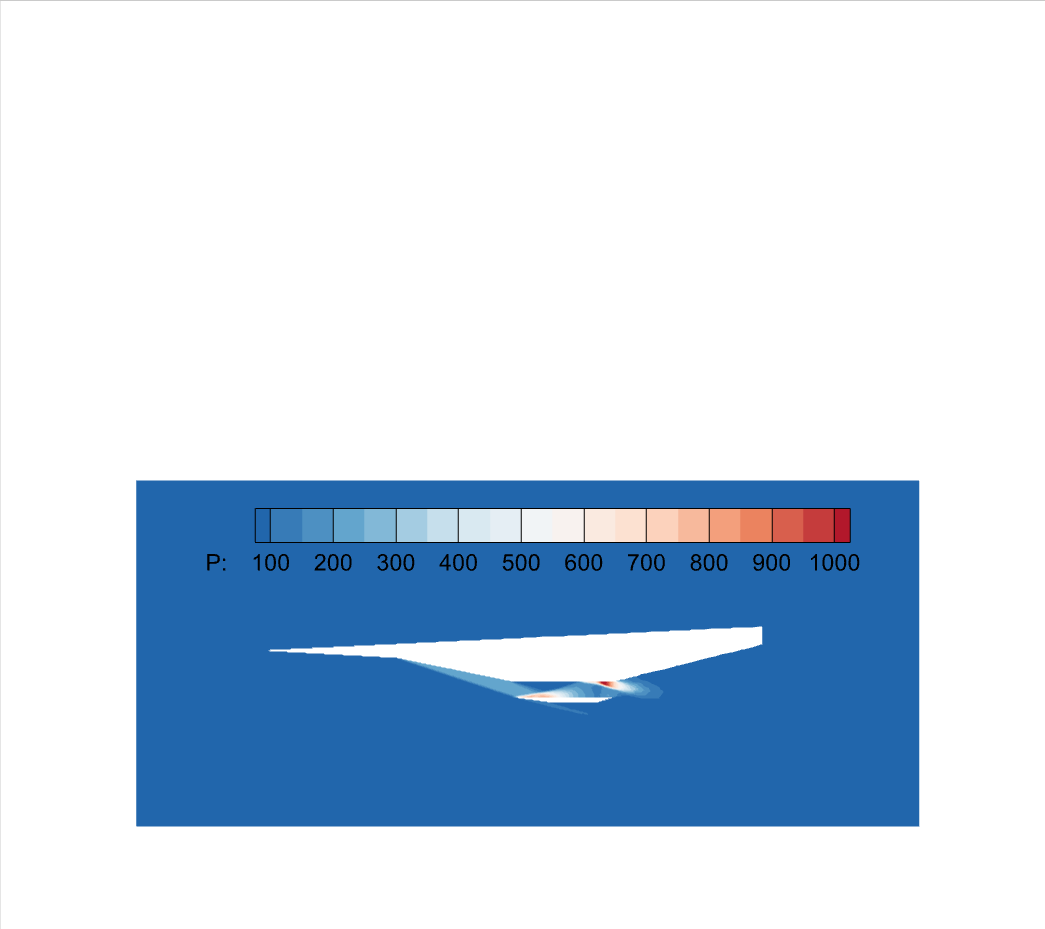}
	\includegraphics[height=0.40\textwidth]{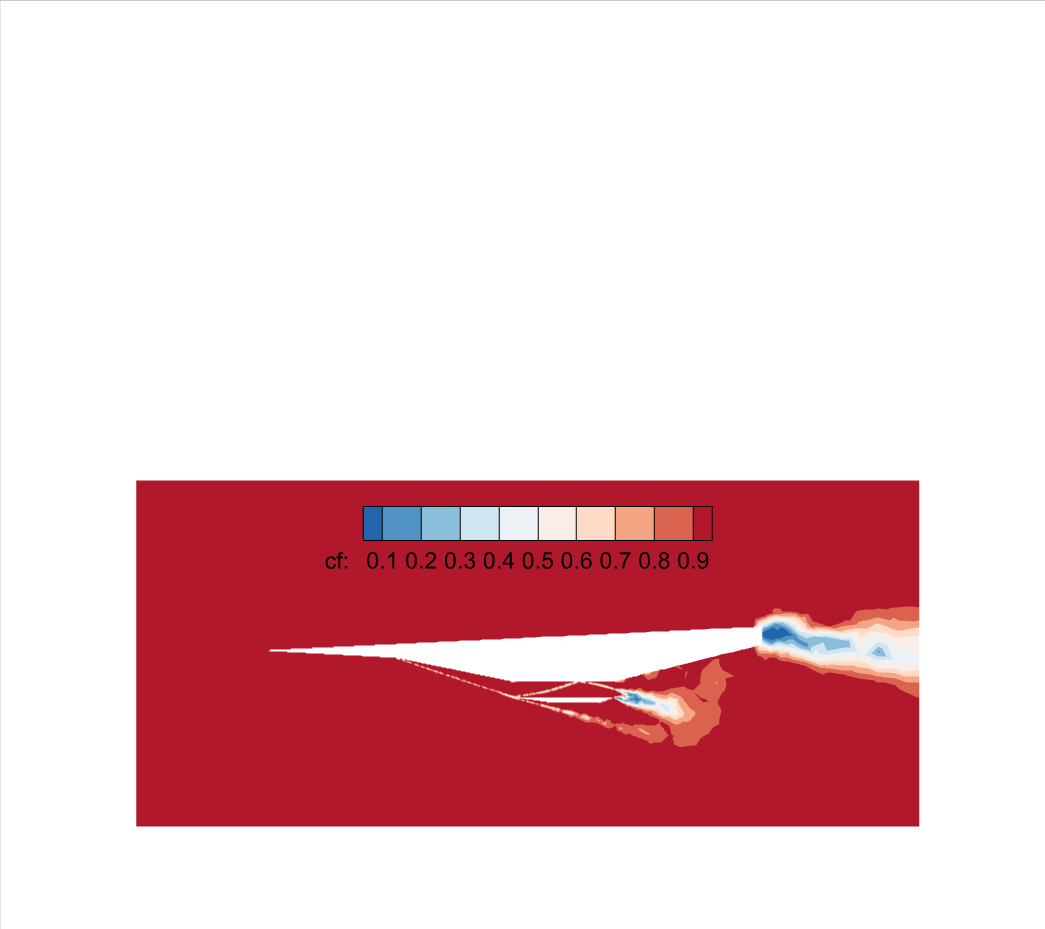}
	\caption{\label{x43-contour}
		Contour of hypersonic flow around X-43A-like aircraft. Left: Pressure contour. Right: DF distribution.}
\end{figure}

The inviscid flow condition is assumed in this case.
The incoming Mach number is 7, and the angle of attack is set to be $1.78^{\circ}$.
The pressure contour and DF distribution are shown in Fig.~\ref{x43-contour}.

\begin{figure}[htp]	
	\centering	
	\includegraphics[height=0.45\textwidth]{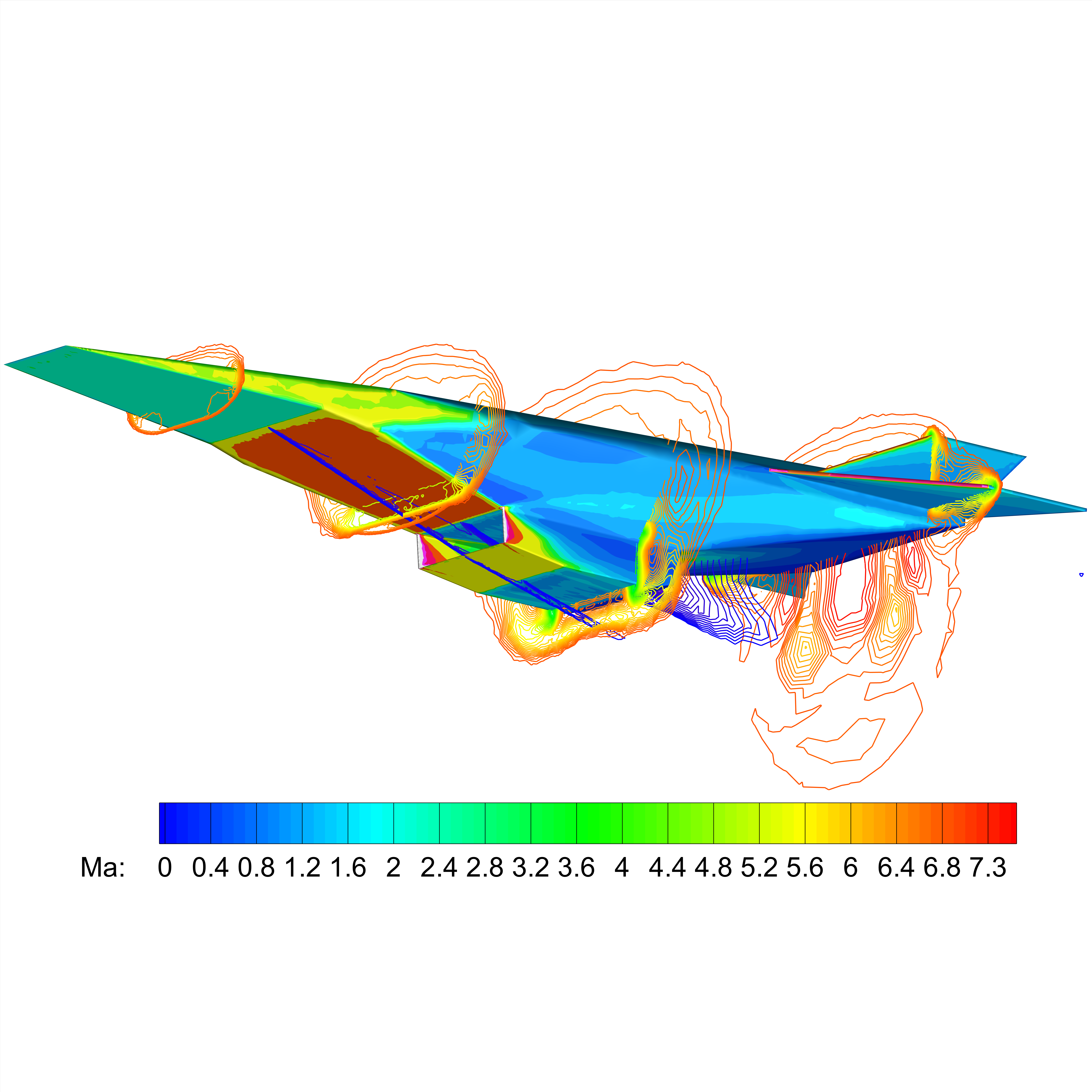}
	\caption{\label{x43-3d}
		3-D contour of hypersonic flow around X-43A-like aircraft.}
\end{figure}

The Mach number distribution of the wall surface with different slices of Mach number distribution in the space and the reference result is shown in Fig.~\ref{x43-3d}.

The pressure contour plot obtained from the memory reduction CGKS calculation shows that a strong shock wave precisely enters the scramjet engine inlet and continuously reflects within the engine, consistent with the design conditions of the aircraft. This demonstrates the accurate shock-capturing capability of the memory reduction CGKS method.

From the three-dimensional contour, we can see the complicated shock structures are captured with low oscillations.
From the DF distribution, only the strong shock and strong reflection wave parts are limited, demonstrating the high resolution of the memory reduction CGKS.

\subsection{Hypersonic flow around a gliding wave rider configuration}
\begin{figure}[htp]	
	\centering	
        \includegraphics[height=0.40\textwidth]{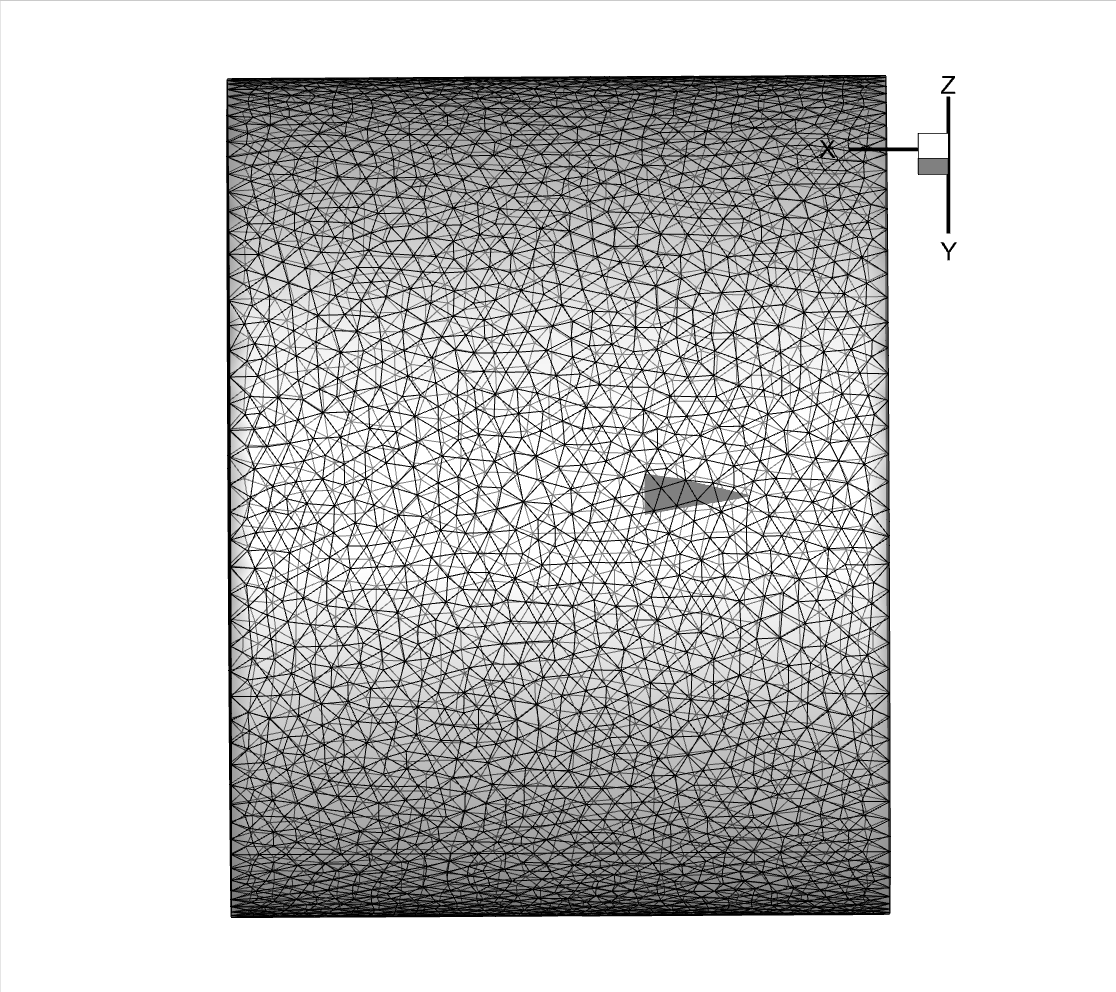}
        \includegraphics[height=0.40\textwidth]{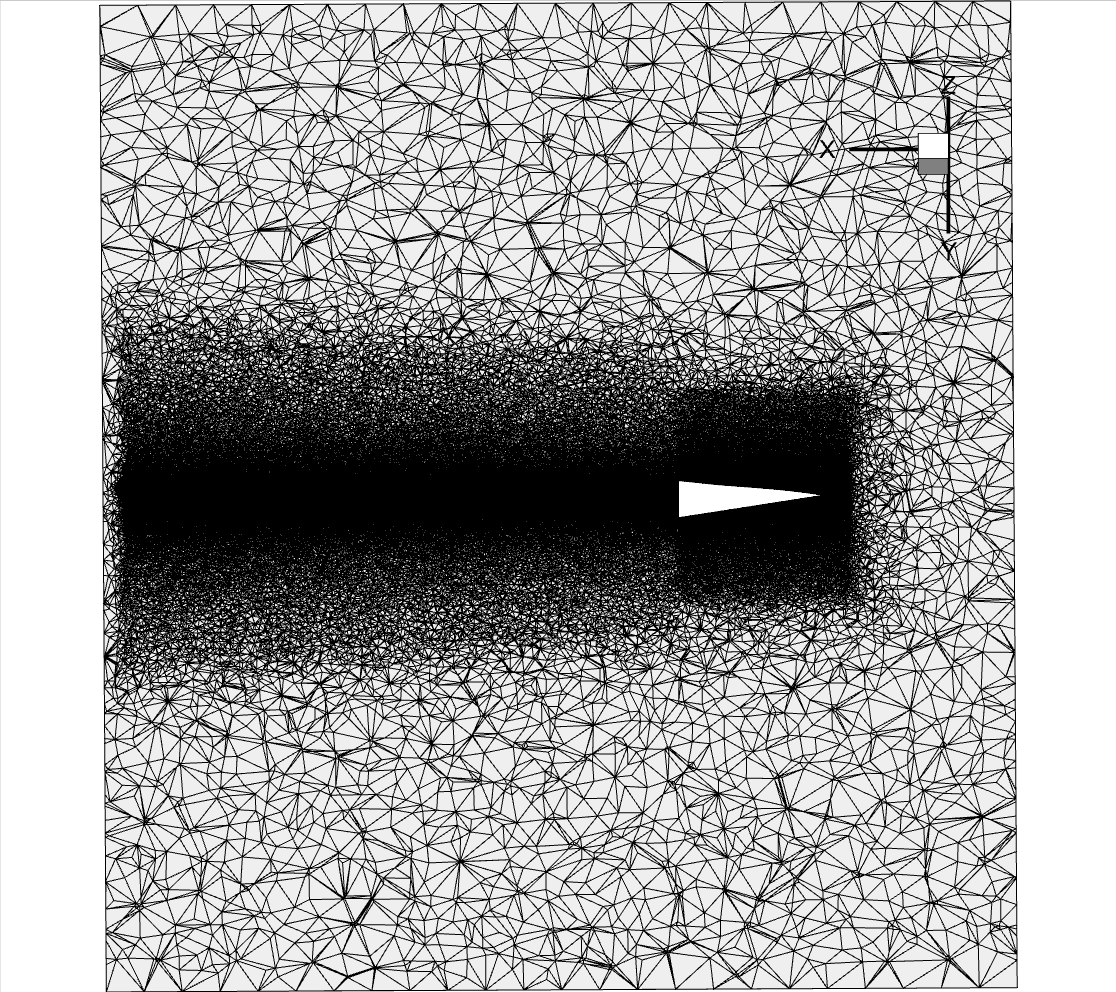}
	\caption{\label{htv2-mesh}
		Mesh of hypersonic flow around HTV-2-like aircraft. Left: Global mesh. Right: Local mesh.}
\end{figure}

In this section, a hypersonic flow around the HTV-2-like aircraft is simulated to further show the high resolution, high efficiency, and strong robustness of the current scheme.
HTV-2-like aircraft is 3 meters long and 1 meter wide. 
To capture the unsteady vortex structure at the tail of the aircraft and the shock at the head, the corresponding parts of the mesh are refined.
The element number of this mixed unstructured mesh is 19,482,823 and $y^+$ is set to be 10, as shown in Fig.~\ref{htv2-mesh}.

\begin{figure}[htp]	
	\centering	
        \includegraphics[height=0.40\textwidth]{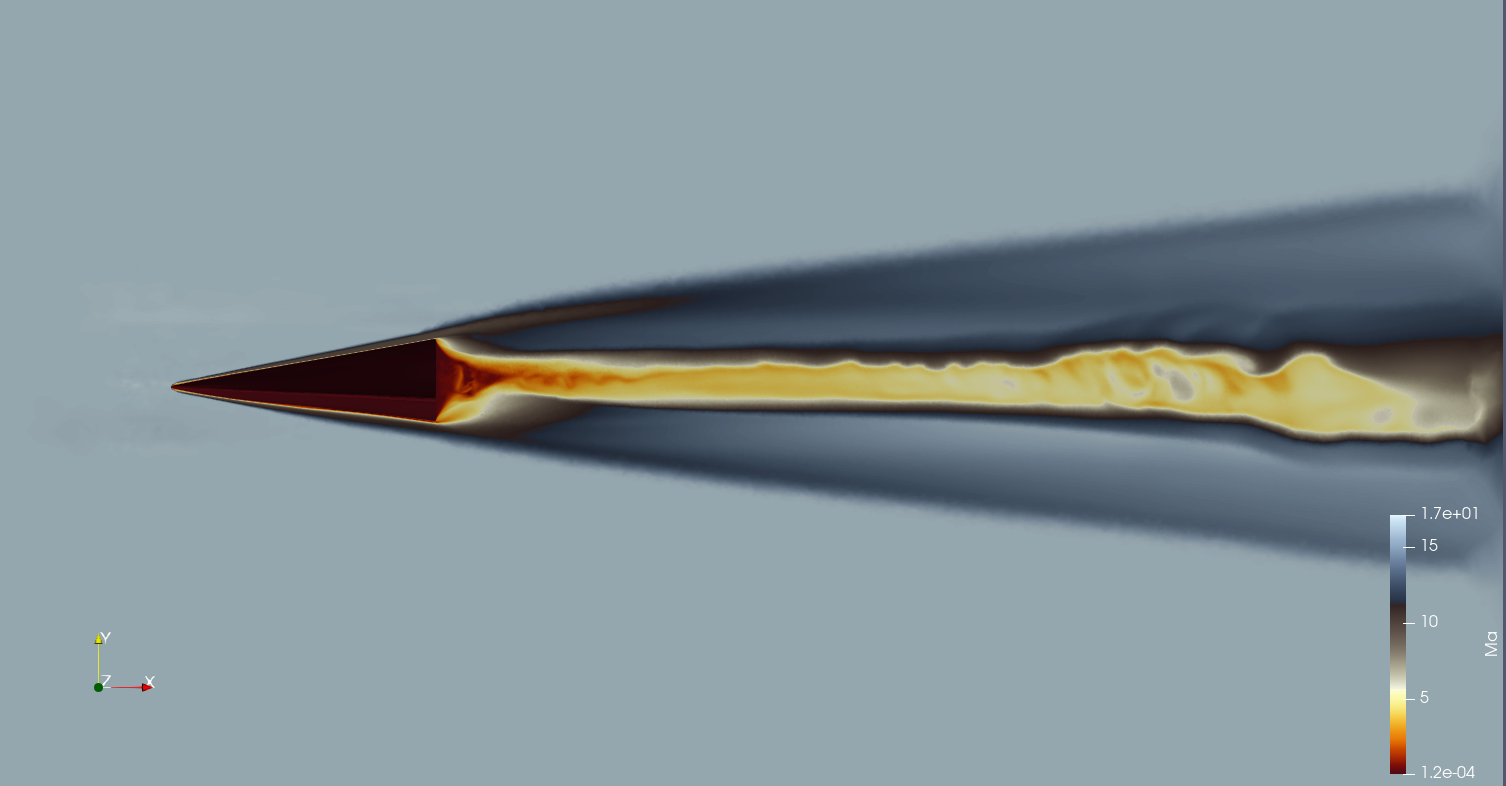}
        \includegraphics[height=0.40\textwidth]{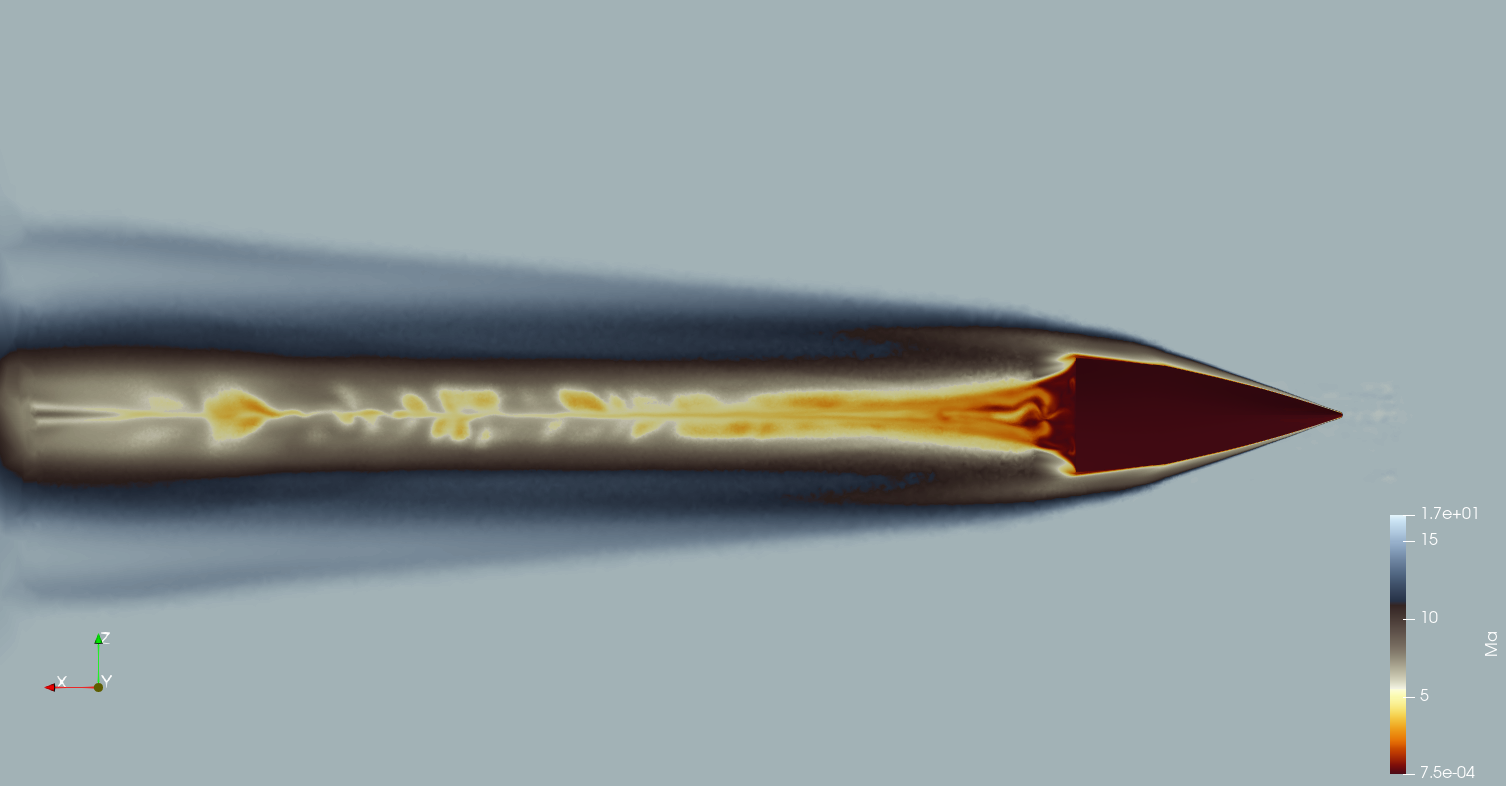}
	\caption{\label{htv2-ma}
		Mach number contour of hypersonic flow around HTV-2-like aircraft. Top: XOY-plane. Down: XOZ-plane.}
\end{figure}

\begin{figure}[htp]	
	\centering	
        \includegraphics[height=0.40\textwidth]{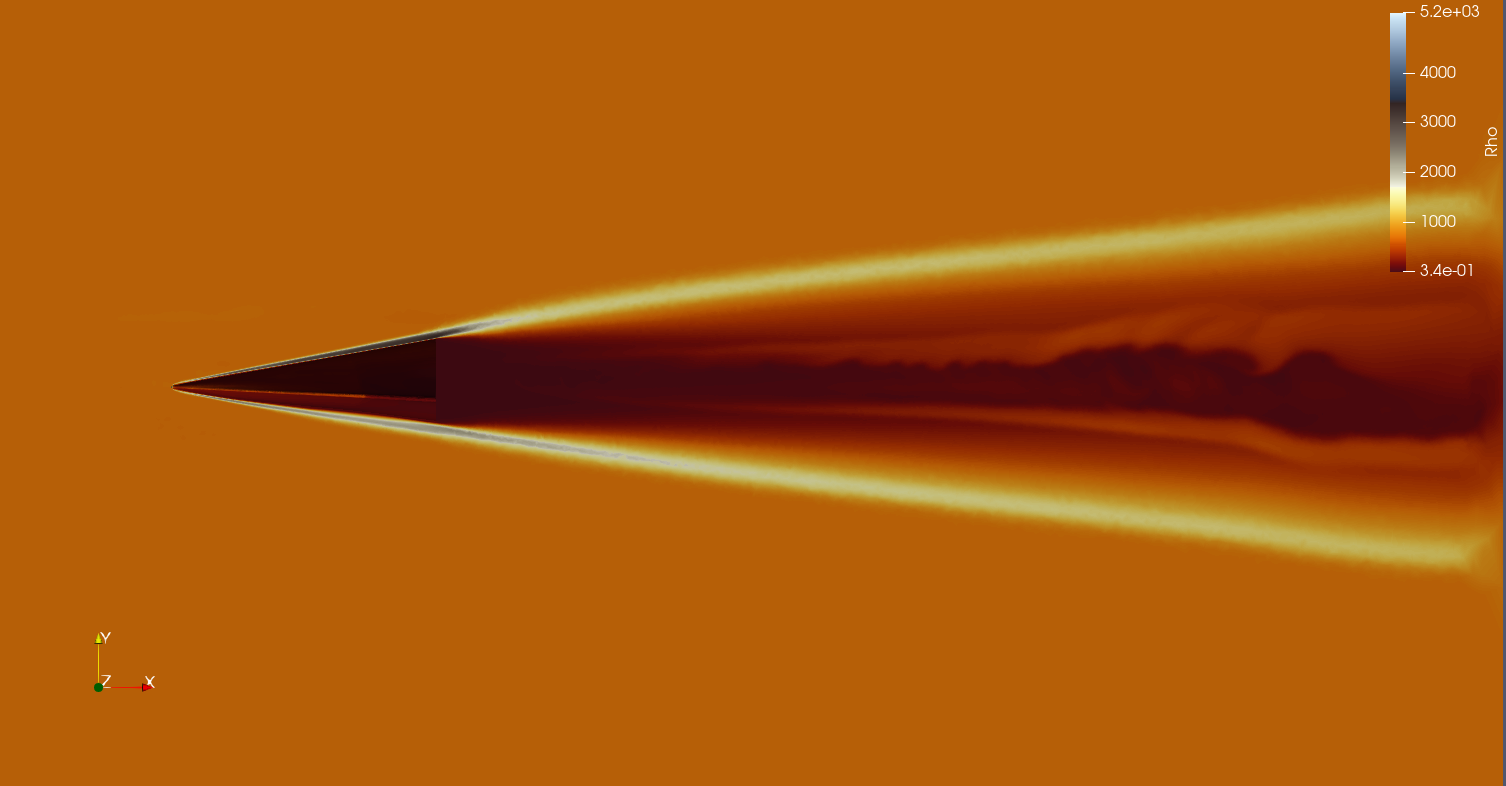}
        \includegraphics[height=0.40\textwidth]{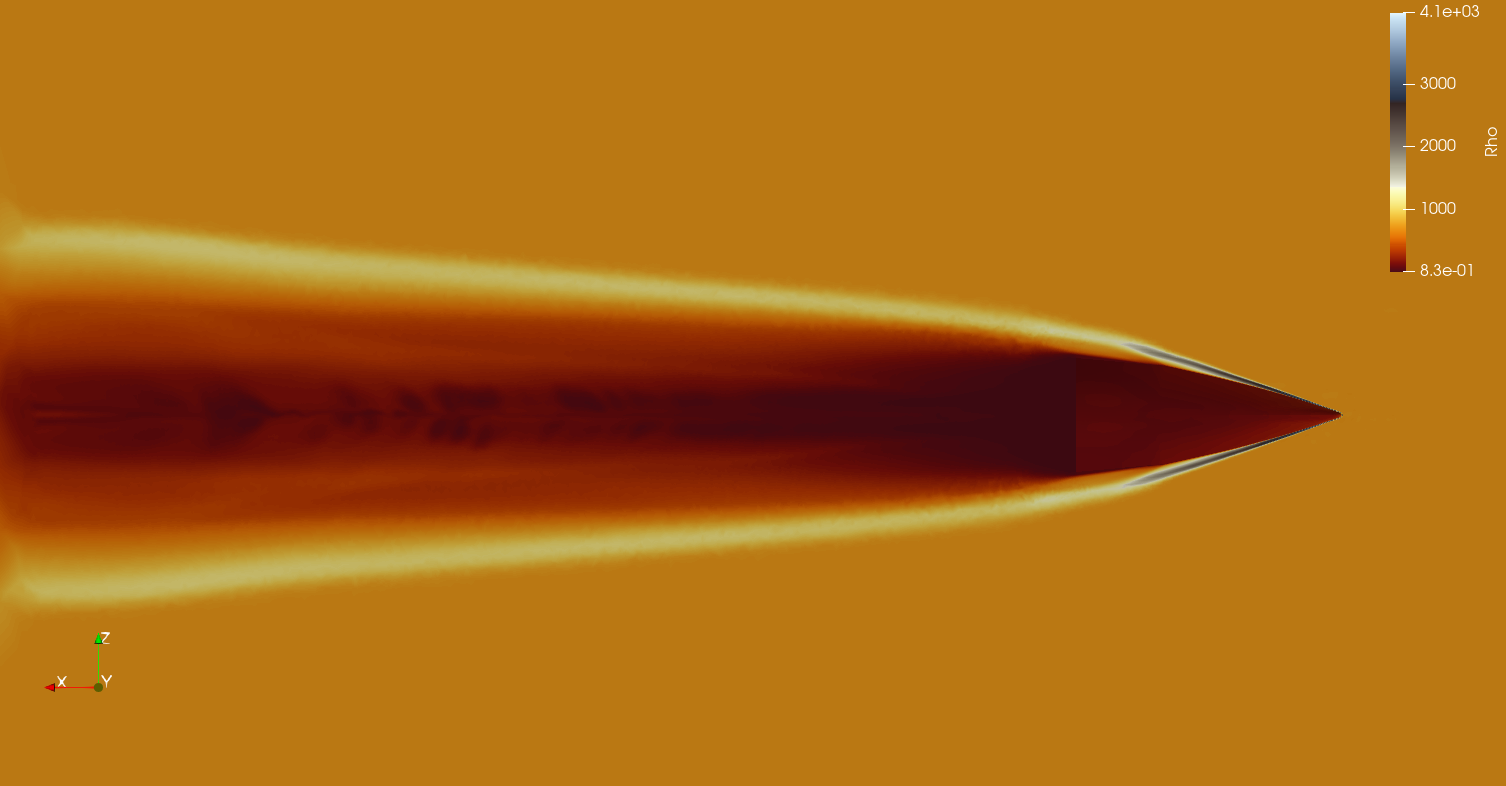}
	\caption{\label{htv2-density}
		Density contour of hypersonic flow around HTV-2-like aircraft. Top: XOY-plane. Down: XOZ-plane.}
\end{figure}

The incoming Mach number is set to be 16.38 and the Reynolds number is set to be 870,000.
The Mach number contour and density contour are shown in Fig.~\ref{htv2-ma} and Fig.~\ref{htv2-density}.

\begin{figure}[htp]	
	\centering	
        \includegraphics[height=0.40\textwidth]{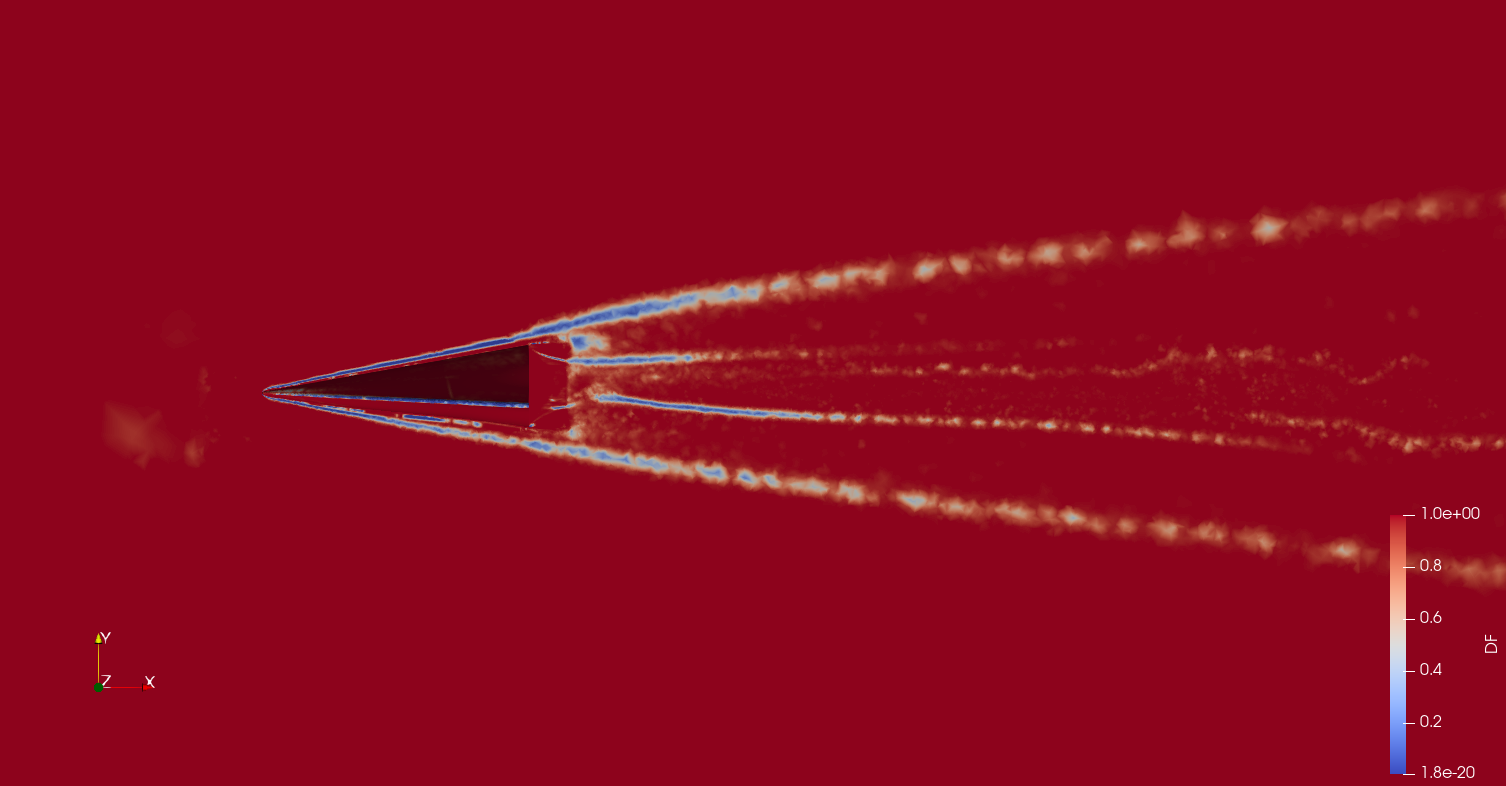}
	\caption{\label{htv2-DF}
		Discontinuity feedback distribution.}
\end{figure}

The above results show that even facing the extremely high Mach number the current scheme can also survive and capture the shock very sharply.
Meanwhile, we can see the high resolution of the vortex at the aircraft's tail, demonstrating the current scheme's low dissipation.
Finally, from the discontinuity feedback distribution, shown in Fig.~\ref{htv2-DF}, only the strong shock part is limited, indicating the accuracy of the shock-capturing and the low dissipation.

\subsection{Efficiency Comparison}
In this section, we use the test cases above to compare the computational efficiency of the current memory reduction CGKS and the original CGKS. The results are listed in Table \ref{cpu-time}.
The accuracy test case, flow around a sphere case and the flow around an M6 wing case are tested on a personal computer using 32 CPU cores. 
The flow around an X-43A-like aircraft, flow around a rocket and the flow around a HTV-2-like aircraft are tested on our clusters using 512 CPU cores. 
The CPU time listed in Table \ref{cpu-time} is the time consumption of every 100 steps.
From the results, improvements in computational efficiency were achieved in all test cases including large-scale simulations.

\begin{table}[htp]
	\small
	\begin{center}
		\def\temptablewidth{1.0\textwidth}
		{\rule{\temptablewidth}{1pt}}
		\begin{tabular*}{\temptablewidth}{@{\extracolsep{\fill}}c|c|c|c|c|c|c}
			Scheme & Accuracy Test & Sphere & M6 Wing  &  X-43A aircraft  & Rocket & HTV-2 aircraft\\
			\hline	
			Original CGKS & 793s & 78s & 321s  & 24s & 98s & 112s \\
                Current  & 571s & 59s & 248s  & 18s & 80s & 86s \\
                Improvement  & 28$$\%$$ & 24$$\%$$ & 23$$\%$$  & 25$$\%$$ & 18$$\%$$ & 23$$\%$$ \\
		\end{tabular*}
		{\rule{\temptablewidth}{0.1pt}}
	\end{center}
	\vspace{-4mm} \caption{\label{cpu-time} 
 CPU time comparison of the original CGKS and current scheme.
 }
\end{table}

\section{Conclusions}
In this paper, we develop a memory reduction third-order spatial reconstruction for large stencils to save memory consumption and enhance computational efficiency.
A two-step third-order linear reconstruction is employed. In the first step, we use the evolved point value on the integration points to obtain the cell-averaged slopes through the Green-Gauss theorem. 
Using cell-averaged slopes to do the least square reconstruction once.
Then, the coefficients of quadratic terms of the polynomial can be obtained.
In the second step, moving the quadratic terms to the right-hand side (RHS) of the original HWENO linear equations, means only linear terms need to be solved. 
Compared with the original reconstruction, the current reconstruction is matrix-free. 
As a result, computational efficiency has been improved. 
We demonstrate the method's performance on 3-D hybrid unstructured meshes, even using tens of millions of grids, suggesting its suitability for large-scale applications, including multi-GPU acceleration.
Our analysis indicates that the current memory reduction CGKS has high resolution and strong robustness from subsonic flow to hypersonic flow.
Future work will explore incorporating the ALE method into current memory reduction CGKS.
Due to the reason that the current scheme does not need to calculate the new reconstruction matrix for each deformation step, solving the problem of the   ALE method on traditional high-order FVM framework.
What's more multi-GPU acceleration technique is also a potential way to explore.

\section*{Acknowledgments}
The current research is supported by National Science Foundation of China (12172316, 12302378, 92371201, 92371107),
Hong Kong Research Grant Council (16208021,16301222).


\bibliographystyle{plain}%

\end{document}